\documentclass[11pt]{article}

\usepackage{amsmath,amsthm,amsfonts,latexsym}
\usepackage{amssymb}
\usepackage{mathabx}
\usepackage[latin1]{inputenc}
\usepackage{tikz-cd}

\begin{document}

\newtheorem*{theo}{Theorem}
\newtheorem*{pro} {Proposition}
\newtheorem*{cor} {Corollary}
\newtheorem*{lem} {Lemma}
\newtheorem{theorem}{Theorem}[section]
\newtheorem{corollary}[theorem]{Corollary}
\newtheorem{lemma}[theorem]{Lemma}
\newtheorem{proposition}[theorem]{Proposition}
\newtheorem{conjecture}[theorem]{Conjecture}

\theoremstyle{definition}
 \newtheorem{definition}[theorem]{Definition} 
  \newtheorem{example}[theorem]{Example}
   \newtheorem{remark}[theorem]{Remark}
   
\newcommand{\Naturali}{{\mathbb{N}}}
\newcommand{\Reali}{{\mathbb{R}}}
\newcommand{\Complessi}{{\mathbb{C}}}
\newcommand{\Toro}{{\mathbb{T}}}
\newcommand{\Relativi}{{\mathbb{Z}}}
\newcommand{\HH}{\mathfrak H}
\newcommand{\KK}{\mathfrak K}
\newcommand{\LL}{\mathfrak L}
\newcommand{\as}{\star_{\sigma}}
\newcommand{\tn}{\vert\hspace{-.3mm}\vert\hspace{-.3mm}\vert}
\def\A{{\mathcal A}}
\def\B{{\mathcal B}}
\def\E{{\mathcal E}}
\def\F{{\mathcal F}}
\def\H{{\mathcal H}}
\def\K{{\mathcal K}}
\def\L{{\mathcal L}}
\def\N{{\mathcal N}}
\def\M{{\cal M}}
\def\gM{{\frak M}}
\def\O{{\cal O}}
\def\P{{\cal P}}
\def\S{{\cal S}}
\def\T{{\cal T}}
\def\U{{\cal U}}
\def\C{{\mathcal C}}
\def\V{{\mathcal V}}
\def\X{{\mathcal X}}
\def\Y{{\mathcal Y}}
\def\qed{\hfill$\square$}
\def\veps{{\varepsilon}}

\title{Actions of Fell bundles }

\author{Erik B\'edos, Roberto Conti\\}
\date{\today}
\maketitle
\markboth{R. Conti, Erik B\'edos}{
}

\vspace{-3ex} \hspace{22ex}{\it  Dedicated to the memory of Sergio Doplicher}

\renewcommand{\sectionmark}[1]{}
\begin{abstract} We introduce and study  actions of Fell bundles over discrete groups on Hilbert bundles.
Many examples of such actions are presented. We discuss the connection with 
positive definite bundle maps between Fell bundles, culminating in the unital case in a Gelfand-Raikov type theorem. We also 
use these actions to construct C*-correspondences over cross-sectional $C^*$-algebras of Fell bundles.
\end{abstract}



\section{Introduction}

In many contexts it is useful to know that a given $C^*$-algebra $A$ is graded over a (nontrivial) discrete group $G$, meaning that there exists a linear independent family of closed subspaces $(A_g)_{g\in G}$ of $A$ satisfying that $A_gA_h\subseteq A_{gh}, \, (A_g)^* \subseteq A_{g^{-1}}$ for all $g,h \in G$ and that the algebraic direct sum of the $A_g$'s is norm-dense in $A$. The collection of all $A_g$'s form then a (concrete) Fell bundle $\A$ over $G$ and, as exposed in Exel's book \cite{Exel2}, understanding the structure of $\A$ is a key step in order to study certain structural properties of $A$. Examples of group-graded $C^*$-algebras are abundant in the literature, and frequently involve (partial) actions or coactions of discrete groups. 
Now, to any (not necessarily concrete) Fell bundle $\A= (A_g)_{g\in G}$, 
one may associate its full (resp.~reduced) cross-sectional $C^*$-algebra $C^*(\A)$ (resp.~$C_r^*(\A)$). Both these $C^*$-algebras are naturally graded over $G$, and their underlying concrete Fell bundles may both be identified with $\A$. Thus  
any Fell bundle over $G$ can be assumed to be concretely given. However, one should be aware that a $G$-graded $C^*$-algebra $A$ with underlying Fell bundle $\A$ may be ``exotic'' in the sense that it is neither isomorphic to $C^*(\A)$, nor to $C_r^*(\A)$.
In good cases, in particular if $\A$ is amenable \cite{Exel1, Exel2}, e.g., if $\A$ has the Exel approximation property, or the PD-approximation property \cite{BeCo7} (see also \cite{BF}), this does not happen. We will not discuss this point any further here and only consider full and reduced cross-sectional $C^*$-algebras associated to Fell bundles in the present paper.

Given two Fell bundles $\A=(A_g)_{g\in G}$ and $\B = (B_h)_{h\in H}$ over possibly different groups $G$ and $H$ connected by a homomorphism $\varphi:G\to H$, one would like to know how to produce various types of bounded linear maps between the associated $C^*$-algebras (full or reduced). At first one may look for maps that by restriction respect the linear structure of the fibers, i.e., arise from so-called $\A$-$\varphi$-$\B$ bundle maps.
Recently \cite{BeCo7}, we introduced a notion of positive definiteness for 
such bundle maps 
and showed that they induce completely positive maps from $C^*(\A)$ to $C^*(\B)$, and  from $C_r^*(\A)$ to $C_r^*(\B)$ whenever $\ker(\varphi)$ is an amenable group. 
When $\A$ is the group bundle over $G$ and $\B$ is the trivial bundle, i.e., $H$ is trivial and 
the unit fiber $B_e$ is equal to $\Complessi$, this notion reduces to the notion of positive definiteness for a complex function $f$ on $G$. As is well-known, this is equivalent to the fact that $f$ is a diagonal coefficient of a unitary representation $u$ of $G$ on some Hilbert space $\mathcal{H}$ over $\Complessi$. A natural question is therefore whether there exist some objects playing the role of $u$ and $\mathcal{H}$ in the general context of Fell bundles. A good candidate for replacing $\mathcal{H}$ is 
a so-called (right) Hilbert $\B$-bundle $\X=(X_h)_{h\in H}$ (see Section \ref{Hilb-B}).  This 
kind of bundles, which were introduced by Abadie and Ferraro in \cite{AF} and studied further in \cite{ABF1}, 
may be thought as the natural generalization of the concept of a (right) Hilbert $C^*$-module $X$ to the realm of Fell bundles. A (left) action of $\A$ on $\X$ should then be the analog of a (left) action of a $C^*$-algebra on $X$ by adjointable operators. This is basically what we propose 
in section \ref{action},
cf.~Definition \ref{BC-action-def}, which involves only a few natural axioms. 
 Such an action on $\X$  is called an $(\A, \varphi)$-action as it depends on $\varphi$.  
 
 Actions of Fell bundles arise in many contexts, as we illustrate in section \ref{action-ex}. 
When $\A=\B$ and $\varphi$ is the identity map on $G$, actions are simply called $\B$-actions.
 In subsection \ref{B-actions}, we exhibit some examples of 
$\B$-actions, such as the trivial and the regular $\B$-actions, which make sense for every Fell bundle $\B$. We also discuss a regularization process for $\B$-actions, producing a regularized $\B$-action from any given $\B$-action.
Subsection \ref{A-phi-actions}
deals with the general case of $(\A, \varphi)$-actions on Hilbert $\B$-bundles. The flexibility of the definition makes it possible to present examples from various sources.
For instance, if $H$ is trivial, 
then  $(\A, \varphi)$-actions correspond to $*$-representations of $\A$ on Hilbert $C^*$-modules. Moreover, one may enhance this basic example to see that $(\A, \varphi)$-actions always exist.
 In another direction, unitary representations of groups  
 and conditional expectations yield actions of  naturally associated Fell bundles. We also mention that the regularization process for $\B$-actions can be extended to $(\A, \varphi)$-actions.  
 In subsection \ref{action-ds}, we show that compatible actions of $C^*$-dynamical systems on $C^*$-correspondences, as studied in \cite{EKQR00, EKQR06},
 give also rise to actions of the Fell bundles associated to such systems.  

Given an $(\A, \varphi)$-action $\rho$ on a Hilbert $\B$-bundle $\X=(X_h)_{h\in H}$, it is elementary to see that
any diagonal coefficient map associated to $\rho$ and  some $x\in X_e$ is  a positive definite $\A$-$\varphi$-$\B$ bundle map in the sense of \cite{BeCo7}.  Section \ref{GR-PD} is devoted to proving that the converse holds when $\A$ and $\B$ are both unital: any positive definite $\A$-$\varphi$-$\B$ bundle map can then be written as a diagonal coefficient map for some suitable $\rho$ and  $x\in X_e$ as above, cf.~Theorem \ref{GR}.   This may be seen as an analog in our setting of the Gelfand-Raikov theorem for positive definite functions on discrete groups and  of the GNS-construction for positive linear functionals on unital $C^*$-algebras. This result also generalizes Paschke's theorem for completely positive maps between unital $C^*$-algebras \cite[Theorem 5.2]{Pas} (by letting both $G$ and $H$ be trivial). In fact, we first give a proof in subsection \ref{GR-1} for the case where $G=H$ and $\varphi$ is the identity map, by applying Paschke's theorem to the completely positive map $\Phi_T:C^*(\A)\to C^*(\B)$ associated to a given positive definite $\A$-$\varphi$-$\B$ bundle map $T$ in \cite[Theorem 3.11]{BeCo7}.  
Next, we sketch a proof of the general case in subsection \ref{GR-II}, following a similar pattern as in the proofs of the classical results mentioned above. 
 The unitality assumption is a reasonable price to pay at this stage in order to avoid introducing an appropriate concept of multiplier Hilbert $\B$-bundles and extending the notion of  $(\A, \varphi)$-action to such bundles (see Remark \ref{non-unital} for some further comments). 
  Theorem \ref{pd-etc}, which gives a sum-up for unital Fell bundles of the connections between actions,  positive definite bundle maps and completely positive maps on cross-sectional $C^*$-algebras, is in our opinion rewarding, and will expectedly be useful in further studies of unital Fell bundles.
  
  In the theory of $C^*$-crossed products associated with $C^*$-dynamical systems, an important observation going back to Combes \cite{Combes} and Kasparov \cite{Kasp, Kasp95} is that one may produce $C^*$-correspondences over the $C^*$-crossed products from compatible actions on Hilbert $C^*$-modules, cf.~\cite{EKQR00, EKQR06}.  An interesting application of actions of Fell bundles is that they can be used in a similar way to construct $C^*$-correspondences. Indeed, given an $(\A, \varphi)$-action $\rho$ on a Hilbert $\B$-bundle $\X$, there is a natural way of defining a (right) Hilbert $C^*(\B)$-module $\X\rtimes \B$ and a (left) action of $C^*(\A)$ on $\X\rtimes \B$, that is, a $C^*(\A)$-$C^*(\B)$ correspondence, cf.~Theorem \ref{correspond}. At the expense of amplification, we can also construct a $C_r^*(\A)$-$C_r^*(\B)$ correspondence. This amplification procedure appears to be unnecessary in certain circumstances (see subsection \ref{C*-c-red}), e.g., when $\A$ and $\B$ are unital and $\rho$ has a cyclic element in $X_e$ (in a suitable sense), as it is always the case when $\rho$ comes from a positive definite bundle map.  To illustrate the usefulness of our construction, we give a new proof of a result of Abadie and Ferraro \cite[Theorem 4.5]{AF}, saying that if $\A$ and $\B$ are (weakly) equivalent Fell bundles (over the same group), then the $C^*$-algebras $C^*(\A)$ and $C^*(\B)$ are Morita equivalent. The existence of the required  $C^*(\A)$-$C^*(\B)$ imprimitivity bimodule follows naturally by considering the action of $\A$ induced by the (weak) equivalence.  
  
 In the last section (Section \ref{final}), we briefly review some possible extensions of our present work.    In order to keep
the paper at a reasonable size, we decided to split our manuscript into two parts. In a forthcoming article \cite{bc25B} we will discuss in detail the notions of direct sum and tensor product of actions of Fell bundles as well as the Fourier-Stieltjes algebra $F(\B)$ of a unital Fell bundle $\B$. 
The latter topic is actually one of our main motivation for introducing actions of Fell bundles, in order to describe $F(\B)$ in terms of such actions, and the present paper provides the 
foundation for our study of $F(\B)$.

\section{Preliminaries} \label{prelim}

Let $A$ be a $C^*$-algebra. We will denote by $A^+=\{ a \in A: a\geq 0\}$ the cone of positive elements in $A$. Further,  
by a Hilbert $A$-module we will always mean a \emph{right} Hilbert $A$-module, as defined in \cite{La1}. If $X$ is such a module, then $\L_A(X)$ will denote the $C^*$-algebra consisting of all adjointable linear maps from $X$ into itself. When $B$ is a $C^*$-algebra and $A$ acts on the left of a Hilbert $B$-module $X$ by adjointable operators (i.e.,  there exists a $*$-representation of $A$ on $X$, meaning that there exists a $*$-homomorphism of $A$ into $\L_B(X)$), we will say that $X$ is a $C^*$-correspondence over $A$ and $B$, or an $A$-$B$ correspondence. Note that we do not assume that $A$ acts non-degenerately on $X$, as some authors do.  

Throughout this paper, 
$G$ and $H$ will always be discrete groups.
By a slight abuse of notation, we will denote the unit of any group by  $e$. If $X$ is a complex vector space,  then $C_c(G, X)$ will denote the vector space of all finitely supported functions from $G$ to $X$ (with pointwise operations).

We will be considering Fell bundles  $\A=(A_g)_{g \in G}$ and $\B=(B_h)_{h \in H}$ over $G$ and $H$, respectively. In particular, 
the unit fibres $A_e$ and $B_e$ are $C^*$-algebras. Our notation and terminology will essentially be as in Exel's book \cite{Exel2}, hence as in \cite{BeCo7}. In particular, $C_c(\A)$ denotes the $*$-algebra consisting of all finitely supported functions $f:G\to \A$ such that $ f(g) \in A_g \text{ for all } g \in G$, the product and involution being given by
\[ (f_1 \star f_2) (h) = \sum_{g\in G} f_1(g) f_2(g^{-1}h), \quad f^*(h) = f(h^{-1})^*\]
for all $h \in G$. Moreover,   $\ell^2(\A) $ will denote the Hilbert $A_e$-module 
consisting of the functions  $\xi:G\to \A$ satisfying that $ \xi(t) \in A_t \text{ for every } t$ in $G \text{ and } \sum_{t\in G} \xi(t)^*\xi(t)$ is (unconditionally) norm-convergent in $A_e$, with operations   
\[ (\xi \cdot a)(t) = \xi(t) a, \quad  \langle \xi, \eta\rangle_{A_e} = \sum_{t\in G}\,\xi(t)^*\eta(t)\quad \text{for every } t\in G,\]
 where $\xi, \eta \in \ell^2(\A)$ and $a \in A_e$. (As a Hilbert $A_e$-module, $\ell^2(\A)$ can also obtained as the completion of the natural 
 inner product  $A_e$-module $C_c(\A)$, cf.~\cite{Exel1, Exel2}.)
 
A $*$-representation $\pi=(\pi_g)_{g\in G}$ of $\A$ in some $*$-algebra $C$ is a family of linear maps $\pi_g:A_g\to C$ satisfying that \[\pi_g(a)\pi_{g'}(a')= \pi_{gg'}(aa') \quad \text{ and }\quad \pi_g(a)^*= \pi_{g^{-1}}(a^*)\] for all $g, g'\in G$ and $a \in A_g, a'\in A_{g'}$. It induces a $*$-homomorphism $\phi_\pi:C_c(\A) \to C$, often called the integrated form of $\pi$,  given by 
\[\phi_\pi(f) = \sum_{g\in G} \pi(f(g)) \quad \text{for all } f \in C_c(\A).\]
For $g \in G$, we let $j^\A_g: A_g \to C_c(\A)$ denote the canonical injection. If $a\in A_g$ we will sometimes write $a\odot g$ instead of $j^\A_g(a)$. We then have that $j^\A=(j^\A_g)_{g\in G}$ is a $*$-representation of $\A$ in $C_c(\A)$ whose integrated form is the identity map on $C_c(\A)$.

The (left) regular representation $\lambda^\A = (\lambda^\A_g)_{g\in G}$ of $\A$ in  $\L_{A_e}(\ell^2(\A))$ is the $*$-representation determined  for each $g\in G$
by
\[ (\lambda^\A_g(a)\xi)(h) = a\,\xi(g^{-1}h) \quad \text{whenever } a \in A_g, \xi \in \ell^2(\A) \text{ and } h \in G.\]
It induces a faithful $*$-representation $\iota^\A = \phi_{\lambda^\A}: C_c(\A)\to \L_{A_e}(\ell^2(\A))$ given  by
\[ \iota^\A(f) = \sum_{g\in G} \lambda^\A_g(f(g)) \quad \text{ for all } f\in C_c(\A),\]
see  \cite[Proposition 17.9, (ii)]{Exel2}.
 The reduced cross-sectional $C^*$-algebra $C_r^*(\A)$ is the $C^*$-subalgebra of $\L_{A_e}(\ell^2(\A))$ generated by 
 $\{ \lambda_g^\A(a): g\in G, a \in A_g\}$, i.e., by $\iota^\A(C_c(\A))$.

The full cross-sectional $C^*$-algebra $C^*(\A)$ is defined as the $C^*$-completion of $C_c(\A)$ with respect to the universal $C^*$-norm  $\|\cdot\|_{\rm u}$ given by \[\|f\|_{\rm u} 
=  \sup\big\{p(f): p \text{ is a } C^*\text{-seminorm on } C_c(\A)\big\} \leq \sum_{g\in G} \|f(g)\|.\]
 For $f\in C_c(\A)$, setting $\|f\|_r := \|\iota^\A(f)\|$, we have $\|f\|_{r} \leq \|f\|_{\rm u}$.

Further, we let $\kappa^\A:C_c(\A)\to C^*(\A)$ 
denote the canonical injection and let  $\widehat j^\A= (\widehat j^\A_g)_{g\in G}$ be the $*$-representation of $\A$ in $C^*(\A)$ given by $\widehat j^\A_g := \kappa^\A\circ j^\A_g: A_g \to C^*(\A)$  for each $g\in G$. As shown in \cite[Proposition 17.9, (iv)]{Exel2}, each $\widehat j^{\A}_g$ is an isometric map. Also, $\kappa^\A$ is the integrated form of $\widehat j^\A$, i.e.,
\[\kappa^\A(f) =  \sum_{g\in G} \widehat j^\A_g(f(g)) \quad \text{for all } f\in C_c(\A).\] 
The canonical $*$-homomorphism $\Lambda^\A$ from $C^*(\A)$ onto $C_r^*(\A)$ is determined by  $\iota^\A = \Lambda^\A \circ \kappa^\A$, i.e., by \[\Lambda^\A\big(\,\widehat{j}_g^\A(a)\big) = \lambda_g^\A(a)\] for all $g\in G$ and $a \in A_g$. 

We will make use of the matrix $C^*$-algebras associated to a Fell bundle $\A$ in \cite{AF}. We recall their definition. Let $g_1, \ldots, g_n \in G$ and set $\mathbf{g}:= (g_1, \ldots, g_n) \in G^n$. 
Then
\[ M_\mathbf{g}(\A) := \big\{R=[r_{ij}]\in M_n(\A) : r_{ij} \in A_{g_i^{-1}g_j} \text{ for all } i, j =1, \ldots, n\big\}\]
is a $*$-algebra with respect to the natural operations, 
which can be equipped with a norm turning it into a $C^*$-algebra, 
cf.~\cite[Lemma 2.8]{AF}.

As in \cite{BeCo7}, we will represent  $M_\mathbf{g}(\A)$ 
by adjointable operators on a  Hilbert $A_e$-module 
in the following way.
Consider $A_\mathbf{g}:= A_{g_1^{-1}}\oplus \cdots \oplus A_{g_n^{-1}}$ as the  Hilbert $A_e$-module obtained by taking the direct sum of the $A_{g_i}$'s (considered as  Hilbert $A_e$-modules), whose inner product is given by
\[ \big\langle (a_1, \ldots, a_n), (a'_1, \ldots, a'_n)\big\rangle_{A_e}  = \sum_{i=1}^n a_i^*a'_i\]
for $a_i, a'_i \in A_{g_i^{-1}}$, $i=1, \ldots, n$.
Then each $R=[r_{ij}] \in M_\mathbf{g}(\A)$ act on $A_\mathbf{g}$ by  
\[L_R(a_1, \ldots, a_n) = (a'_1, \ldots, a'_n), \ \text{ where } a'_i := \sum_{j=1}^n r_{ij}a_j \text { for each }i = 1, \ldots, n. \] 
It is easy to check that the operator $L_R : A_\mathbf{g}\to A_\mathbf{g}$ is adjointable with $(L_R)^*=L_{R^*}$, and that the map $R \mapsto L_R$ from $M_\mathbf{g}(\A)$ into  $\mathcal{L}_{A_e}(A_\mathbf{g})$ is an injective $*$-homomorphism. 
Thus, the norm on $M_\mathbf{g}(\A)$ satisfies that $\|R\|= \|L_R\|$.

We recall some other relevant definitions from \cite{BeCo7}.
Let $\varphi: G \to H$ be a group-homomorphism, i.e., $\varphi\in {\rm Hom}(G,H)$. An 
\emph{$\A$-$\varphi$-$\B$ bundle map}
is a family of maps $T = (T_g)_{g\in G}$  such that $T_g: A_g \to B_{\varphi(g)}$ is linear and bounded for every $g\in G$. 
When $G=H$ and $\varphi = {\rm id}_G$, so $\A$ and $\B$ are both Fell bundles over $G$, we simply say $\A$-$\B$ \emph{bundle map} instead of  $\A$-${\rm id}_G$-$\B$ bundle map. Moreover, if $\A=\B$, we just say \emph{$\B$-bundle map} instead of $\B$-$\B$ bundle map. 

Let $T = (T_g)_{g\in G}$ be an $\A$-$\varphi$-$\B$ bundle map. Then we can define
a linear map $\phi_T: C_c(\A) \to C_c(\B)$ by setting 
\[  \phi_T\big(f)= \sum_{g\in G} j^\B_{\varphi(g)}(T_g(f(g)))\quad \text{for all } f\in C_c(\A).\]
 We say that $T$ is \emph{reduced $($resp.~full\,$)$} if $\phi_T$ is bounded w.r.t.~the reduced (resp.~universal) $C^*$-norms on $C_c(\A)$ and $C_c(\B)$, and denote the extension of $\phi_T$ to a linear bounded map from $C_r^*(\A)$ into  $C_r^*(\B)$ (resp.~from $C^*(\A)$ into $C^*(\B)$) by $M_T$ (resp.~$\Phi_T$).

An $\A$-$\varphi$-$\B$ bundle map
  $T = (T_g)_{g\in G}$ is said to be \emph{positive definite} if
\begin{equation}\label{Bposdef}
\sum_{i, j=1}^n b_i \, T_{g_i^{-1}g_j} (a_i^* a_j) \, b_j^* \, \in  (B_e)^+
\end{equation}
for all $n\in \Naturali$, $g_1, \ldots, g_n \in G$ and all  $a_i \in A_{g_i}$, $b_i\in B_{\varphi(g_i)}$, $i=1, \ldots, n$.
Equivalently, cf.~\cite[Remark 3.4]{BeCo7}, 
  for every  $n\in \Naturali$, $\mathbf{g}= (g_1, \ldots, g_n) \in G^n$ and $a_i \in A_{g_i}$, $i=1, \ldots, n$, 
  the matrix 
\[ \Big[ T_{g_i^{-1}g_j} (a_i^* a_j)\Big]\] is positive in the $C^*$-algebra $M_{\varphi(\mathbf{g})}(\B)$, where $\varphi(\mathbf{g}):= (\varphi(g_1), \ldots, \varphi(g_n)) \in H^n$.

The following result is a combination of Theorem 3.11 and Theorem 3.14 in \cite{BeCo7}: 

\begin{theorem} \label{PDCP}
Let $\A=(A_g)_{g \in G}$ and $\B=(B_h)_{h \in H}$ be Fell bundles and let $T=(T_g)_{g\in G}$ be an $\A$-$\varphi$-$\B$ bundle map, where $\varphi \in {\rm Hom}(G,H)$.

Then  $T$ is positive definite
if and only if $T$ is full and  $\Phi_T: C^*(\A) \to C^*(\B)$ is completely positive.
Moreover, if $\ker(\varphi)$ is an amenable group, then this is also equivalent to $T$ being reduced with $M_T: C^*_r(\A) \to C^*_r(\B)$ being completely positive.
\end{theorem}

\section{ Hilbert bundles over Fell bundles} \label{Hilb-B}

The following definition is due to Abadie and Ferraro in \cite{AF}, except that we follow \cite{ABF1} and don't include the fullness condition. 
 
\begin{definition}\label{rHilbertBbundle} Let $\B=(B_g)_{g\in G}$ be a Fell bundle over $G$.
 A (\emph{right}) \emph{Hilbert $\B$-bundle}\footnote{With the exception of subsection \ref{Morita}, we will not consider \emph{left} Hilbert $\B$-bundles in this paper. 
 By a  
 Hilbert $\B$-bundle we otherwise always mean a \emph{(right)} Hilbert $\B$-bundle.}  is a family of 
complex 
Banach
spaces $\X = (X_r)_{r\in G}$ indexed over $G$
 such that if $\X$ also denotes the disjoint union of the $X_r$'s,  we have maps \[(x, b) \mapsto x b : \X \times \B \to \X, \, \text{ and }\, (x, y) \mapsto \langle x, y\rangle_\B : \X \times \X \to \B,\] satisfying the following conditions: 
\begin{itemize}
\item[(1)] For each $r, s\in G$, we have $X_rB_s \subseteq X_{rs}$, and the map  $\,(x,b)\mapsto xb$ from $X_r \times B_s$ into $X_{rs}$ is bilinear.
\item[(2)] For each $r, s\in G$, we have  $\langle X_r,X_s\rangle_\B \subseteq B_{r^{-1}s}$ and, for each $x\in X_r$,  the map $  \, y\mapsto \langle x,y\rangle_\B$ from $X_s$ into $B_{r^{-1}s}$ is linear. 

\smallskip \noindent (Hence, in particular, we have $\langle X_r,X_r\rangle_\B \,\subseteq \,B_e$ for every $r\in G$.)
\item[(3)] $\langle x,yb\rangle_\B = \langle x,y\rangle_\B\, b$  and $(\langle x,y\rangle_\B)^* = \langle y,x\rangle_ \B$ for all $x,y \in \X$ and $b \in \B$.
\item[(4)] $\langle x,x\rangle_\B \geq  0$ in $B_e$ for all $x \in \X$, and $\langle x,x\rangle_\B = 0$ implies $x = 0$. 
\item[(5)] 
For each $r\in G$, the norm $x\mapsto \|x\|$ on $X_r$ satisfies that $\|x\|=\| \langle x,x\rangle_\B \|^{1/2}$ for all $x\in X_r$.
\end{itemize} 
Note that for each $r \in G$, conditions (1)-(4) imply that $X_r$ is an  inner product $B_e$-module, hence that 
$x\mapsto \| \langle x,x\rangle_\B \|^{1/2}$ gives a norm on $X_r$ (cf.~\cite{La1}).  
Since $X_r$ is complete by assumption, it follows from 
condition (5) that $X_r$ is a  Hilbert $B_e$-module. 
Moreover,  it can be shown  (cf.~\cite[Lemma 2.7]{AF}) that for all $x, y\in \X$ and $b,c\in \B$, we  also have
\begin{itemize}
\item[(a)] \,$\langle xb,y\rangle_\B =b^*\langle x,y\rangle_\B$,
\item[(b)] \, $\|xb\| \leq \|x\| \,\|b\|$,
\item[(c)] \, $\|\langle x,y\rangle_\B\| \leq \|x\| \,\|y\|$, 
\item[(d)] \, $(xb)c = x(bc)$.
\end{itemize}
\end{definition}

\begin{remark} \label{pre_H} 
If $\X = (X_r)_{r\in G}$ is a family of 
(complex) vector 
spaces equipped with maps $\X \times \B \to \X$ and  $\X \times \X \to \B$  satisfying conditions (1)-(3) and 
$\langle x,x\rangle_\B \geq  0$ in $B_e$ for all $x \in \X$ (where $\X$ denotes the disjoint union of the $X_r$'s), then $\X$ will be called a 
\emph{semi-inner product $\B$-bundle}.
The map $x\mapsto \|x\|:= \| \langle x,x\rangle_\B \|^{1/2}$ is then a seminorm on each $X_r$, and the proof of \cite[Lemma 2.7]{AF} makes it clear that properties (a)-(d) above still hold. After a separation and completion procedure,  one obtains from $\X$ a  Hilbert $\B$-bundle $\widetilde{\X}$ (see \cite[Proposition 5.1]{AF} for a similar result). If the map $x\mapsto \|x\|$ is already a norm on each $X_r$, that is, if (4) also holds,  then $\X$ is called an
\emph{inner product $\B$-bundle}, and only completion
is needed to construct $\widetilde{\X}$.
\end{remark}

\begin{example}
The  \emph{trivial $\B$-bundle} is $\B$ itself (with $\langle b, c\rangle_\B := b^*c$). 
 \end{example}
 
  \begin{example} \label{reg-bundle}
Let $\B=(B_g)_{g\in G}$ be a Fell bundle over $G$. Assume  $\X:=(X_r)_{r\in G}$ is a Hilbert $\B$-bundle. We can define its regularization $\widecheck{\X}= (\widecheck X_r)_{r\in G}$ as follows. For each $r \in G$, set $\widecheck X_r^0 := C_c(G, X_r)$.
Then for $r,s\in G$, $\xi \in \widecheck X_r^0$, $\eta \in \widecheck X_s^0$ and $b \in B_s$, we define
 $\xi\cdot b \in \widecheck X_{rs}^0$ and $\langle \xi, \eta\rangle_\B \in B_{r^{-1}s}$ by
 \begin{align*} (\xi\cdot b)(t) &:= \xi(t) b \quad \text{for all} t \in G,\\
 \langle \xi, \eta\rangle_\B&:= \sum_{t\in G} \,\langle \xi(t), \eta(t)\rangle_\B \,. 
 \end{align*}
 Thus, letting $\widecheck{\X}^0$ denote the disjoint union of the $\widecheck{\X}^0_r$'s, we get maps $\widecheck{\X}^0 \times \B \to \widecheck{\X}^0$ and  $\widecheck{\X}^0 \times \widecheck{\X}^0 \to \B$.
 It is a routine exercise to check that $\widecheck {\X}^0=(\widecheck X^0_r)_{r\in G}$ becomes an inner product $\B$-bundle. We can now let $\widecheck \X$
  be the Hilbert $\B$-bundle obtained from $\widecheck \X^0$ after 
   completion. We will call it the \emph{regular Hilbert $\B$-bundle associated to $\X$}. 
 
  As a special case, we may choose $\X$ to be   the trivial Hilbert $\B$-bundle, namely $\B$ itself. We will refer to the regular Hilbert $\B$-bundle $\widecheck \B$ associated to $\B$ simply as \emph{the regular Hilbert $\B$-bundle}.  
  \end{example}
  
\begin{example} \label{ell2}
Let $\B=(B_g)_{g\in G}$ be a Fell bundle over $G$. In \cite[Section 3]{ABF1} the authors construct a certain  Hilbert $\B$-bundle $\Y=\mathcal{L}^2\B$ which they call the \emph{canonical $\ell^2$-bundle} of $\B$, cf.~\cite[Definition 3.2]{ABF1}. For completeness, we 
sketch this construction.
  For each $r\in G$,  consider $Y_r := \ell^2(\B) \times \{r\}$ as a Banach space w.r.t.~the norm 
 $\|(\xi, r)\| = \|\langle \xi, \xi\rangle\|^{1/2}$. Moreover, for $r, s\in G$, $ \xi, \eta\in \ell^2(\B)$ and $b \in B_s$, set
 \[ (\xi, r)\cdot b := (\xi\cdot b, rs), \quad \big\langle (\xi, r), (\eta, s)\big\rangle_\B:= \sum_{t\in G}\,\xi(tr)^*\eta(ts),\]
  the function $\xi\cdot b \in \ell^2(\B) $ being defined by
 \[ (\xi\cdot b)(t) := \xi(ts^{-1}) b \quad \text{ for each } t\in G.\] 
 Then $\Y:=(Y_r)_{r\in G}$ becomes a  Hilbert $\B$-bundle, 
  cf.~\cite[Proposition 3.1]{ABF1}. 
 \end{example}

\begin{example}\label{eqrep1}
Let $\B_\Sigma= (A\times\{g\})_{g\in G}$ be the canonical Fell bundle associated to a (discrete) $C^*$-dynamical system $\Sigma= (A, G, \alpha)$, where $\alpha: G \to {\rm Aut}(A)$ is an action of $G$ on a $C^*$-algebra $A$, cf.~\cite{Exel2}. 
We recall that the operations and the norms are given by \[(a, g) (a', g') = (a\alpha_g(a'), gg'),  (a, g)^* = (\alpha_{g^{-1}}(a^*), g^{-1}) \text{ and } \|(a, g)\| = \|a\|. \]  It is common to identify the unit fiber $A\times \{e\}$ with $A$ via the map $(a,e)\mapsto a$. 

Let $X$ be a  Hilbert $A$-module. We can then construct 
a Hilbert $\B_\Sigma$-bundle $\X$ associated to $X$ as  follows.
 For each $g\in G$, set $X_g:= X\times \{g\}$  and organize $X_g$ as a Banach space in the obvious way. 
 Moreover, for $x,y \in X, a\in A$ and $g, h \in G$, set
\begin{align*}
(x, g) \cdot (a, h) &:= \big(x\,\alpha_g(a), gh\big), \\
\big\langle (x, g), \, (y, h)\big\rangle_{\B_\Sigma} &:= \big(\alpha_g^{-1}(\langle x, y\rangle_A), g^{-1}h\big).
\end{align*}
 
 It is straightforward to check that conditions (1)-(5) hold. 
 As a sample, we check (3). For $x,y \in X, a\in A$ and $g, h, k \in G$, we get
\begin{align*}
 \big\langle (x, g),  \, (y, h) \cdot(a, k)\big\rangle_{\B_\Sigma} &=  \big\langle (x, g), \, (y\,\alpha_h(a), hk)\big\rangle_{\B_\Sigma} = \big(\alpha_g^{-1}(\langle x, y\,\alpha_h(a)\rangle), g^{-1}hk\big) \\
& = \big(\alpha_g^{-1}(\langle x, y\rangle) \, \alpha_{g^{-1}h}(a), g^{-1}hk\big) = \big(\alpha_g^{-1}(\langle x, y\rangle), g^{-1}h\big)\, (a,k) \\ 
& = \big\langle (x, g), \, (y, h)\big\rangle_{\B_\Sigma} \, (a,k), 
\end{align*}
so the first assertion in (3) holds. Further, 
we also get
\[ \big(\big\langle (x, g), (y, h)\big\rangle_{\B_\Sigma}\big)^* =  \big(\alpha_{(g^{-1}h)^{-1}}\big(\alpha_{g}^{-1}(\langle x, y\rangle_A)^*\big), (g^{-1}h)^{-1}\big) = \big(\alpha_{h}^{-1}(\langle x, y\rangle_A^*), h^{-1}g\big) \]
\[= \big(\alpha_{h}^{-1}(\langle y, x\rangle_A), h^{-1}g\big) = 
\big\langle (y, h), \, (x, g) \big\rangle_{\B_\Sigma},\]
showing that the second assertion in (3) holds. Thus, $\X=(X_g)_{g\in G}$ becomes a  Hilbert $\B_\Sigma$-bundle with these operations. We note that each $X_g$ can be identified as a  Hilbert $A$-module with the 
Hilbert $A$-module $X_{\alpha_g}$ obtained as the twisting of  $X$ by $\alpha_g \in {\rm Aut}(A)$.
\end{example}

\begin{example}\label{sub-b}
 Assume $\B = (B_g)_{g\in G}$ is a Fell sub-bundle of a Fell bundle $\A=(A_g)_{g\in G}$
 and $E= (E_g)_{g\in G}$ is a conditional expectation from $\A$ to $\B$, as defined in \cite[Definition 21.5]{Exel2} and  \cite[Definition 21.19]{Exel2}, respectively.
 It is then easy to check that $\A$ becomes a  semi-inner product 
 $\B$-bundle with respect to the map $\A\times \B \to \A, 
(a, b) \mapsto ab$, coming from the product map in $\A$,  and the map $\A \times \A \to \B, (a,a')\mapsto \langle a, a'\rangle_\B$,  given by
\[\langle a, a'\rangle_\B: = E_{g^{-1}g'}(a^*a')\] 
whenever $g, g'\in G, a\in A_g, a'\in A_{g'}$. Clearly, $\A$ is an inner product $\B$-bundle if  $E_e$ is faithful as a conditional expectation from $A_e$ onto $B_e$.
\end{example}

We will later need the following lemma, which generalizes the first part of \cite[Lemma 4.2]{La1}. 

\begin{lemma} \label{pos-mat} Let $\X = (X_g)_{g\in G}$ be a  Hilbert $\B$-bundle, $n\in \Naturali$ and $\mathbf{g}= (g_1, \ldots, g_n) \in G^n$.  
For $j =1, \ldots, n$, let $x_j \in X_{g_j}$. Set $R := \big[ \langle x_i, x_j\rangle_\B\big]$.
Then $R$ is positive in $M_{\mathbf{g}}(\B)$.
\end{lemma} 
\begin{proof}
Since $ \langle x_i, x_j\rangle_\B \in B_{g_i^{-1}g_j}$ for each $i, j\in \{1, \ldots, n\}$, $R  \in M_{\mathbf{g}}(\B)$. 
Note that for $\mathbf{b} = (b_1, \ldots , b_n) \in B_{\mathbf{g}}$, we have $x_i b_i \in X_{g_i}B_{g_i^{-1}} \subseteq X_e$ for each $i$, so $z:=\sum_{i=1}^n x_ib_i \in X_e$. Hence we get 
\[ \langle \mathbf{b}, L_R(\mathbf{b})\rangle_{B_e} = \sum_{i=1}^n  b_i^*\Big( \sum_{j=1}^n \langle x_i, x_j\rangle_\B\, b_j\Big) = \sum_{i, j=1}^n \langle x_i\,b_i, x_j\, b_j\rangle_\B = \langle z, z\rangle_\B \in (B_e)^+. \]
This shows that $L_R$ is positive in $\mathcal{L}_{B_e}(B_\mathbf{g})$, which implies that $R $ is positive in $M_{\mathbf{g}}(\B)$.
 \end{proof}

A natural notion of unitary equivalence for Hilbert $\B$-bundles is as follows. 
  \begin{definition} \label{unitary-eq} Let $\X=(X_g)_{g \in G}, \X' = (X'_g)_{g \in G}$ be Hilbert $\B$-bundles. We will denote the $\B$-valued inner product on $\X'$ by $\langle \cdot, \cdot\rangle_\B'$. Let $U=(U_g)_{g\in G}$ be a $\B$-bundle map from $\X$ to $\X'$. We will say that $U$ is \emph{unitary} if the following conditions hold:
 
 \begin{itemize}
 \item $U_g$ is bijective, for each $g \in G$. 
 \item $U_g(x) b = U_{gh}(xb)$
for all $g,h \in G$, $x \in X_g$ and $b \in B_h$.
 \item  $\langle U_g(x), U_h(y)\rangle_\B' = \langle x, y\rangle_\B$
   for all $g, h \in G,  x\in X_g, y\in X_h$. 
   \end{itemize}
Moreover, we will say that \emph{$\X$ and $\X'$ are unitarily equivalent} when such a unitary $\B$-bundle map $U$ exists.
 \end{definition}

 \begin{remark}
 Let $g\in G$. Considering $X_g$ and $X_g'$ as Hilbert $B_e$-modules, we get readily from 
 the previous definition that $U_g \in \L_{B_e}(X_g, X'_g)$ is unitary.
   In particular, every $U_g$ is isometric.
 \end{remark}

\section{Actions of Fell bundles on Hilbert bundles} \label{action}

In this section, we consider a Fell bundle   $\A=(A_g)_{g\in G}$ over a group $G$, a Fell bundle $\B=(B_h)_{h\in H}$ over a group $H$ and some $\varphi \in {\rm Hom}(G,H)$. Recall that $\A$ also denotes the disjoint union of the $A_g$'s, and the same with $\B$. 
\subsection{Definition and properties}
\begin{definition}\label{BC-action-def}
 Let $\X=(X_h)_{h\in H}$ be a  Hilbert $\B$-bundle and denote also by $\X$ the disjoint union of the $X_h$'s. Moreover, let $\F(\X) $ be the semigroup of maps from $\X$ into itself (w.r.t.~composition). A map $\rho: \A \to \F(\X)$ will be called an 
 \emph{$(\A, \varphi)$-action on $\X$} 
when it satisfies the following
conditions:
\begin{itemize}
\item[(i)] for each $g\in G$ and $h\in H$, we have $\rho(A_g)X_h \subseteq X_{\varphi(g)h}$, and the map
$(a, x) \mapsto \rho(a)x $ from $A_g \times X_h$ into $X_{\varphi(g)h}$ is bilinear;
\item[(ii)] $\rho(aa') = \rho(a)\rho(a')$ for all $a, a'\in \A$;
\item[(iii)] $\langle \rho(a) x, y\rangle_\B =   \langle x, \rho(a^*) y\rangle_\B$ for all 
$a \in \A$ and $x, y\in \X$;
\item[(iv)]  $(\rho(a)x)b = \rho(a)(x b)$ for all $a\in \A$, $b\in \B$ and $x \in \X$.
\end{itemize}
(Such a map $\rho$ could alternatively be called a representation of $(\A, \varphi)$ on $\X$.) 
When $G=H$ and $\varphi={\rm id}_G$, we will just say that $\rho$ is an 
\emph{$\A$-action} on the  Hilbert $\B$-bundle $\X$. 
\end{definition}

\begin{remark} \label{action-adjoint} 
We note that condition (iii) implies that if $a\in A_e$ and $h\in H$, then the restriction 
$\rho_h(a)$ 
of $\rho(a)$ to $X_h$ is a linear map from the Hilbert $B_e$-module $X_h$ into itself
which is adjointable, with $\rho_h(a)^* = \rho_h(a^*)$. One readily gets that the map $\rho_h:A_e \to \L_{B_e}(X_h)$ mapping each $a\in B_e$ to $\rho_h(a)$ is a $*$-homomorphism. 

More generally, let $g\in G$, $h\in H$, and $a\in A_g$. Then the restriction 
$\rho_{g,h}(a)$ 
of $\rho(a)$ to $X_h$ is a linear map from the Hilbert $B_e$-module $X_h$ into the Hilbert $B_e$-module $X_{\varphi(g)h}$
which is adjointable and satisfies $\rho_{g, h}(a)^* = \rho_{g^{-1}, \,\varphi(g)h}(a^*)$ for all $g\in G$ and $h\in G$. 
Note that for $a\in A_e$, $\rho_h(a)= \rho_{e,h}(a)$. 
Moreover, if $g, k\in G$,  $a\in A_g, a'\in A_{g'}$ and $h\in H$, then
$\rho_{g', \varphi(g)h}(a')\rho_{g, h}(a)= \rho_{g'g,h}(a'a)$.
 \end{remark} 

Actions of Fell bundles are automatically bounded in the following sense:
 \begin{lemma} \label{actionineq}
  Let $\X = (X_h)_{h\in H}$ be a  Hilbert $\B$-bundle and $\rho$ be an $(\A, \varphi)$-action on $\X$. Then for all $a\in \A$ and $x \in \X$ we have 
\[  \| \rho(a) x\| \, \leq \, \|a\| \, \|x\|
\, \text{ and }\,
\langle \rho(a)x, \rho(a)x\rangle_\B \, \leq \, \|a\|^2 \,\langle x, x\rangle_\B.
\]
 \end{lemma}
 \begin{proof} We use the same notation as in the remark above.
 Let $g\in G$, $h\in H$, $a\in A_g$ and $x \in X_h$. Then $a^*a \in A_e$ and
\[\rho_{g,h}(a)^*\rho_{g,h}(a) = \rho_{g^{-1}, \,\varphi(g)h}(a^*)\rho_{g,h}(a)
 =  \rho_{e, h}(a^*a)= \rho_{h}(a^*a).\]
 Hence, since $\rho_h$ is contractive (being a $*$-homomorphism between $C^*$-algebras), we get
\begin{align*} \| \rho_{g,h}(a) x\|^2 &= \| \langle \rho_{g,h}(a) x, \rho_{g,h}(a) x\rangle_\B\| =
\| \langle  x, \rho_{g,h}(a)^*\rho_{g,h}(a) x\rangle_\B\|\\
&= \| \langle  x, \rho_{h}(a^*a) x\rangle_\B\| \leq \, \|\rho_{h}(a^*a)\| \, \|x\|^2 \\
&\leq \|a^*a\| \, \|x\|^2 = \|a\|^2 \, \|x\|^2.
\end{align*}
This shows that $\| \rho(a)x\|=\| \rho_{g,h}(a)x\| \leq \|a\| \,\|x\| $, and also that $\| \rho_{g,h}(a)\|^2 \leq \|a\|^2$.

\medskip Moreover, as $\rho_{g,h}(a) \in \L_{B_e}(X_h, X_{\varphi(g)h})$,  \cite[Proposition 1.2]{La1} gives  that
\[ \langle \rho(a)x, \rho(a)x\rangle_\B=
  \langle \rho_{g,h}(a) x, \rho_{g,h}(a) x\rangle_\B \, \leq \, \|\rho_{g,h}(a)\|^2 \langle x, x\rangle_\B
  \leq \|a\|^2 \,  \langle x, x\rangle_\B , \]
  as desired.
 \end{proof}
 
 \begin{remark} \label{pre-act-def} One may define an $(\A, \varphi)$-\emph{pre-action} $\rho_0$ on a semi-inner product $\B$-bundle $\X^0= (X_h^0)_{h\in H}$ to be a map $\rho_0: \A\to \F(\X^0)$ satisfying the analogue of conditions  (i)-(iv) in Definition \ref{BC-action-def} and also the inequality
\[ \langle \rho_0(a)x, \rho_0(a) x\rangle_\B \, \leq \, \|a\|^2  \langle x, x\rangle_\B\] for all $a\in \A$ and $x \in \X^0$. Then it is not difficult to check that $\rho_0$ naturally induces an $(\A, \varphi)$-action $\rho$ on the Hilbert $\B$-bundle $\X$ obtained as the Hausdorff completion of $\X^0$.
\end{remark}

An important feature of actions of Fell bundles is that they give rise to positive definite bundle maps by choosing diagonal coefficients over the unit fibre. 

\begin{proposition} \label{action-pd} Let $\rho $ be an $(\A, \varphi)$-action  on a  Hilbert $\B$-bundle $\X$ and let $x \in X_e$. Then, for each $g\in G$, the map $T_g: A_g \to B_{\varphi(g)}$ given by \[T_g(a):= \langle x, \rho(a) x\rangle_\B\] for every $a\in A_g$ is well-defined, linear and bounded. Moreover,  
$T := (T_g)_{g\in G}$ is a positive definite $\A$-$\varphi$-$\B$ bundle map.
 \end{proposition} 

\begin{proof} 
Let $g \in G$ and $a\in A_g$. Then $\rho(a)x \in X_{\varphi(g)}$, so  $T_g(a) \in B_{e\varphi(g)} = B_{\varphi(g)}$.  
Hence each $T_g$ maps $A_g$ into $B_{\varphi(g)}$ and it is clearly linear and bounded. Thus $T$ is  an $\A$-$\varphi$-$\B$ bundle map.

\noindent Further, consider $g, g' \in G$, $a \in A_g, a' \in A_{g'}$ and $b \in B_{\varphi(g)}, c\in B_{\varphi(g')}$. Then we get that 
\[ b \,T_{g^{-1}g'}(a^* a')\, c^* = b\, \big\langle x, \rho(a^*a') x\big\rangle_\B \, c^* = b\, \big\langle x, \rho(a^*) \rho(a') x\big\rangle_\B \, c^* \] 
\[ = b\, \big\langle \rho(a) x, \rho(a') x\big\rangle_\B \, c^* =  \big\langle (\rho(a) x)\, b^*, (\rho(a') x)\,c^*\big\rangle_\B. \]
Let now $n\in \Naturali$, $g_1, \ldots, g_n \in G$ and $a_i \in A_{g_i}$, $b_i\in B_{\varphi(g_i)}$, $i=1, \ldots, n$. Then,
 setting  $z :=  \sum_{i =1}^n (\rho(a_i) x)\, b_i^* \in X_e$ and using the calculation above, we get that
\begin{equation*}
\sum_{i, j=1}^n b_i \, T_{g_i^{-1}g_j} (a_i^* a_j) \, b_j^* = \sum_{i, j =1}^n \big\langle (\rho(a_i) x)\, b_i^*, (\rho(a_j) x)\,b_j^*\big \rangle_\B = \langle z, z\rangle_\B
\end{equation*}
which is positive in $B_e$. This shows that $T$ is positive definite, as desired.
\end{proof} 

\begin{remark}\label{pre-act-pd} In the context of Proposition \ref{action-pd}, it is apparent from its proof that 
we also obtain a positive definite $\A$-$\varphi$-$\B$ bundle map $T$ if $\rho$ is 
assumed to be 
an $(\A, \varphi)$-pre-action on a semi-inner product $\B$-bundle $\X$. 

Also, if $x, y\in X_e$ and we define for each $g\in G$ the map $T^{\rho, x, y}_g: A_g \to B_{\varphi(g)}$  by \[T^{\rho, x, y}_g(a):= \langle x, \rho(a) y\rangle_\B\] for every $a\in A_g$, then $T^{\rho, x, y} := (T^{\rho, x, y}_g)_{g\in G}$ is an $\A$-$\varphi$-$\B$ bundle map, such that $T:=T^{\rho, x, y}$ is full with $\Phi_T$ being completely bounded, and such that $T$ is reduced with $M_T$ being completely bounded if ${\rm ker}(\varphi)$ is amenable. Indeed, by polarization of $T$, we can write $T$ as a linear combination of bundle maps of the form $T^{\rho, z_i, z_i}$ for some $z_i$'s in $X_e$ and apply Proposition \ref{action-pd} and Theorem \ref{PDCP} to each of these. A similar argument works in the reduced case if $\ker \varphi $ is amenable.
\end{remark}

\begin{proposition} Assume  $\X=(X_h)_{h \in H}$ and $ \X' = (X'_h)_{h \in H}$ are unitarily equivalent Hilbert $\B$-bundles and $\rho: \A \to \F(\X)$ is an $(\A, \varphi)$-action on $\X$. Let $U=(U_h)_{h \in H}$ be a unitary $\B$-bundle map from $\X$ to $\X'$. Then the formula
 \begin{equation} \label{rhoprime} \rho'(a)(U_h x) = U_{\varphi(g)h}(\rho(a)x), \ a \in A_g, x \in X_h\end{equation}
defines an  $(\A, \varphi)$-action $\rho'$ on $\X'$. 
\end{proposition}

\begin{proof} We check that the four conditions are satisfied. 
\begin{itemize}
\item[(i)] For each $g\in G$ and $h\in H$,  it is clear that $\rho'(A_g)X'_h \subseteq X'_{\varphi(g)h}$ and that the map
$(a, x') \mapsto \rho'(a)x' = U_{\varphi(g)h}\rho(a)U_h^{-1}x' $ from $A_g \times X'_h$ into $X'_{\varphi(g)h}$ is bilinear.
\item[(ii)]  For all $a \in A_g, a'\in A_{g'}$ and $x'\in X'_k$, we have
\begin{align*}
\rho'(aa')x'&= U_{\varphi(gg')k}\rho(aa')U_k^{-1}x' = U_{\varphi(gg')k}\rho(a)\rho(a')U_k^{-1}x'\\
&=U_{\varphi(g)\varphi(g')k}\rho(a)U_{\varphi(g')k}^{-1}U_{\varphi(g')k}\rho(a')U_k^{-1}x' =\rho'(a)(U_{\varphi(g')k}\rho(a')U_k^{-1}x')\\
&=\rho'(a)(\rho'(a)x') = \rho'(a)\rho'(a')x'
\end{align*}
\item[(iii)] For all 
$g\in G$, $h, k\in H$,  $a \in A_g$, $x' \in X'_h$  and $y'\in X'_k$, setting $x:=U_h^{-1}x' \in X_h$ and $y:= U_k^{-1}y' \in X_k$, we get 
 \begin{align*} 
 \langle \rho'(a) x', y'\rangle'_\B &=  \langle U_{\varphi(g)h}\rho(a)U_h^{-1}x' , U_ky\rangle'_\B
 =  \langle \rho(a)U_h^{-1}x' , y\rangle_\B \\ 
 &=  \langle U_h^{-1}x' , \rho(a^*)y\rangle_\B =  \langle U_h^{-1}x' , \rho(a^*)U_k^{-1}y'\rangle_\B\\ 
 &= \langle U_h^{-1}x', U_{\varphi(g^{-1})k}^{-1}U_{\varphi(g)^{-1}k}\rho(a^*)U_k^{-1}y'\rangle_\B\\
 &=\langle x', U_{\varphi(g^{-1})k}\rho(a^*)U_k^{-1}y'\rangle'_\B  =\langle x', \rho'(a^*) y'\rangle'_\B 
 \end{align*}
\item[(iv)] Let $g\in G$,  $h, k\in H$, $a \in B_g, x' \in X_h$ and $c\in B_k$. Set $x= U_h^{-1}x' \in X_h$.
Then we have 
\[U_{hk}^{-1}(x'c) = U_{hk}^{-1}(U_h(x)c) = U_{hk}^{-1}U_{hk}(xc) = xc = (U_h^{-1}x')c,\]
so we get 
\begin{align*}
(\rho'(a)x')c &= (U_{\varphi(g)h}\rho(a)U_h^{-1}x')c = U_{\varphi(g)hk}\big((\rho(a)U_h^{-1}x')c\big)\\
&=U_{\varphi(g)hk}\rho(a)((U_h^{-1}x')c) 
= U_{\varphi(g)hk}\rho(a)U_{hk}^{-1}(x'c) =\rho'(a)(x'c)
\end{align*} 
\end{itemize}

\vspace{-3ex} \end{proof}

\begin{definition}\label{unitary-eq2}
We will say that two $(\A, \varphi)$-actions $\rho$ and $\rho'$ on Hilbert $\B$-bundles $\X$ and $\X'$, respectively, are \emph{unitarily equivalent} when there is a unitary $\B$-bundle map $U$ from $\X$ to $\X'$ such that (\ref{rhoprime}) holds. 
When $G=H$ and $\varphi={\rm id}_G$, we will just say that the two $\A$-actions are unitarily equivalent.  
\end{definition}

We will later need the following result, which generalizes the second part of \cite[Lemma 4.2]{La1}. 
 
\begin{lemma}\label{Lance-gen}
  Let $\X = (X_h)_{h\in H}$ be a  Hilbert $\B$-bundle and $\rho$ be an $(\A, \varphi)$-action on $\X$.  
Let $a\in A_e$, $\mathbf{h}= (h_1, \ldots, h_n) \in H^n$  and $x_j \in X_{h_j}$ for $j =1, \ldots, n$.  

\smallskip \noindent Set $R := \big[ \langle x_i, x_j\rangle_\B\big] \in M_{\mathbf{h}}(\B)^+$, $S:= \big[ \langle \rho(a) x_i, \rho(a) x_j\rangle_\B\big] \in M_{\mathbf{h}}(\B)^+$, cf.~Lemma \ref{pos-mat}.

\medskip \noindent 
Then we have $S \, \leq \|a\|^2\, R.$
\end{lemma} 
\begin{proof}
Let $\mathbf{b} = (b_1, \ldots, b_n)\in \B_{\mathbf{h})} = B_{h_1^{-1}} \oplus \cdots \oplus B_{h_n^{-1}}$. Then  $y:=   \sum_{i=1}^n  \,x_i b_i$ belongs to $X_e$ and we have
\begin{align*}
\langle \mathbf{b}, L_S(\mathbf{b})\rangle_{B_e} &= \sum_{i =1}^n b_i^* \Big(\sum_{j=1}^n \langle \rho(a) x_i, \rho(a) x_j\rangle_\B\, b_j\Big)  =
\sum_{i, j=1}^n  \, \langle (\rho(a) x_i)b_i, (\rho(a) x_j)b_j\rangle_\B\\
& = \sum_{i, j=1}^n  \, \langle \rho(a) (x_ib_i), \rho(a) (x_jb_j)\rangle_\B = \langle \rho_e(a)y, \rho_e(a)y\rangle_{B_e} \\
& \leq \,\|\rho_e(a)\|^2\, \langle y, y\rangle_{B_e}  \,\leq \,\|a\|^2 \,\sum_{i, j=1}^n  \, \langle x_ib_i, x_jb_j\rangle_{B_e}\\  
&= \|a\|^2 \,\sum_{i=1}^n  \, b_i^* \Big(\sum_{j=1}^n\langle x_i, x_j\rangle_{\B}\,b_j\Big) = \|a\|^2\, \langle \mathbf{b}, L_R(\mathbf{b})\rangle_{B_e} 
\end{align*}
where the first inequality comes from \cite[Proposition 1.2]{La1} and the second from the fact that $\rho_e$ is a $*$-homomorphism of $A_e$ into $\L_{B_e}(X_e)$. This implies that  $ L_S \, \leq \|a\|^2\, L_R$, from which the stated inequality follows,  using that $T\mapsto L_T$ is an injective $*$-homomorphism. 
\end{proof}

\subsection{Examples of actions} \label{action-ex}
Our aim in this subsection is to convince the reader that there is an abundance of examples of actions of Fell bundles.

\subsubsection{$\B$-actions on $\B$-bundles} \label{B-actions}

 \begin{example} \label{trivialBact} The \emph{trivial $\B$-action} $\rho_{\B}$ 
on the trivial  Hilbert $\B$-bundle $\B$ is given by setting $\rho_{\B}(b) c := bc$ for all $b , c\in \B$.
Let $a\in B_e$. Then for $g\in G$ and $b\in B_g$, we have
\[ \big\langle a , \rho_\B(b) a \big\rangle_\B = a^*ba\,.\]
Hence we get from Proposition \ref{action-pd} that $T = (T_g)_{g\in G} $ given by 
$T_g(b) = a^*ba$ for every $g$ in $G$ and $b \in B_g$
is a positive definite $\B$-bundle map.
\end{example}

 \begin{example} \label{regB-action} Let $\rho: \B \to \F(\X)$ be a $\B$-action on a Hilbert $\B$-bundle $\X=(X_r)_{r\in G}$.  We may then form the \emph{regular $\B$-action $\widecheck \rho$ associated to $\rho$} as follows. Let  $\widecheck \X=(\widecheck X_r)_{r\in G}$ be the regular Hilbert $\B$-bundle associated to $\X$, cf.~Example \ref{reg-bundle}, and  $r, t \in G$, $b\in B_t$. Then we define $\widecheck\rho_0(b): \widecheck X_r^0 \to \widecheck X_{tr}^0$ by setting
 \[ (\widecheck\rho_0(b) \xi) (s) := \rho(b)\,\xi(t^{-1}s) \quad \text{for all } \xi \in \widecheck X_r^0 \text{ and } s\in G.\]
 It is obvious that the map $(\xi, b) \mapsto \widecheck \rho_0(b) \xi$ from 
 $ B_t \times \widecheck X_r^0 $ to $\widecheck X_{tr}^0$ is bilinear. Moreover, condition (ii) and (iv) in Definition \ref{BC-action-def} follow for $\widecheck\rho_0$ by routine computations. To see that condition (iii) also holds for $\widecheck\rho_0$, let $r, s, t \in G$,  $b\in B_t$, $\xi \in \widecheck X_r^0$ and $\eta \in \widecheck X_s^0$. Then we have
  \begin{align*}
 \big\langle \widecheck\rho_0(b) \xi , \eta \big\rangle_\B &= \sum_{g\in G} \langle \rho(b)\,\xi(t^{-1}g) ,\eta(g) \rangle_\B = \sum_{g\in G} \langle\,\xi(t^{-1}g) ,\rho(b^*)\eta(g) \rangle_\B \\
 & = \sum_{g\in G} \langle\,\xi(g) ,\rho(b^*)\eta(tg) \rangle_\B
 =  \big\langle  \xi , \widecheck\rho_0(b^*)\eta \big\rangle_\B.
 \end{align*}

Now, using Lemma \ref{actionineq}, we get
 \begin{align*}
 \langle \widecheck\rho_0(b) \xi , \widecheck\rho_0(b)\xi \rangle_\B 
 &=  \sum_{s\in G} \,\langle(\widecheck\rho_0(b)\xi)(s),  (\widecheck\rho_0(b)\xi)(s)\rangle_\B
 =\sum_{s\in G} \,\langle\rho(b)\xi(t^{-1}s),  \rho(b)\xi(t^{-1}s)\rangle_\B\\
 &= \sum_{s\in G} \,\langle\rho(b)\xi(s),  \rho(b)\xi(s)\rangle_\B\,
 \leq \sum_{s\in G} \|b\|^2 \langle \xi(s),\xi(s)\rangle_\B 
  =  \|b\|^2 \,\|\langle \xi, \xi\rangle_\B.
\end{align*}
Hence, $\widecheck\rho_0$ is a $\B$-pre-action on $\widecheck\X^0=(\widecheck X_r^0)_{r\in G}$, which can be extended to a $\B$-action $\widecheck\rho$ on $\X$, cf.~Remark \ref{pre-act-def}. 

As a special case, let us consider the trivial action $\rho_\B$ of $\B$ on itself. Then the regular $\B$-action 
$\widecheck{\rho_\B}$ associated to $\rho_\B$, which acts on the regular Hilbert $\B$-bundle $\widecheck\B$, will be called \emph{the regular $\B$-action}. On the underlying inner product $\B$-bundle $\big(C_c(G,B_r)\big)_{r\in G}$,  it is given by
 \[ (\widecheck{\rho_\B}(b) \xi) (s) := b\,\xi(t^{-1}s) \quad \text{for all } \xi \in C_c(G, B_r) \text{ and } s\in G.\]
for all $r, t \in G$ and $b\in B_t$. Let 
$\xi \in C_c(G, B_e)$. Then for $g\in G$ and $b\in B_g$, we have 
\[ \big\langle \xi , \widecheck{\rho_\B}(b) \xi \big\rangle_\B = \sum_{h\in G} \xi(h)^* \, b\, \xi(g^{-1}h) =  \sum_{h\in G} \xi(gh)^* \, b\, \xi(h).\]
Hence we get from Proposition \ref{action-pd} that $T = (T_g)_{g\in G} $ given by 
\[T_g(b) = \sum_{h\in G} \xi(gh)^* \, b\, \xi(h) \] 
$\text{for every } b \in B_g$
is a positive definite $\B$-bundle map. Note that $T$ is a $\B$-bundle map of Exel-type, as discussed in \cite[Example 4.5]{BeCo7}. 

\end{example}

\begin{example} \label{ell2-2}
One can also define a natural $\B$-action on the canonical $\ell^2$-bundle $\Y=(Y_r)_{r\in G}$ associated to $\B$, cf.~Example \ref{ell2}.
For  $r, s\in G$ and $b\in B_r$, let  $\rho(b): Y_s \to Y_{rs}$  be given by
\[ \rho(b) (\xi, s) := (\lambda^\B_r(b) \xi, rs)  \]
for every $(\xi, s)\in Y_s = \ell^2(\B) \times \{s\}$.
 We can now check that the resulting map $\rho: \B \to F(\Y)$ is a $\B$-action on $\Y$. Condition (i) is obviously true. Condition (ii)  follows readily from the fact that $\lambda^\B$ is a $*$-representation.
   To verify condition (iii), let $r, s, t \in G$,  $b\in B_r$ and $\xi , \eta \in \ell^2(\B)$. Then we have
 \begin{align*}
 \big\langle \rho(b) (\xi, s) , (\eta, t)\big\rangle_\B &=  \big\langle  (\lambda^\B_r(b) \xi, rs) , (\eta, t)\big\rangle = 
 \sum_{k \in G} ((\lambda^\B_r(b) \xi)(krs))^* \eta(kt)\\
 &= \sum_{k \in G} (b \xi(r^{-1}krs))^* \eta(kt) = \sum_{k \in G}  \xi(r^{-1}krs)^* b^*\eta(kt)
  \end{align*}
 On the other hand, we have
\begin{align*}
 \big\langle  (\xi, s) , \rho(b^*)(\eta, t)\big\rangle_\B &=\langle  (\xi, s) , (\lambda^\B_{r^{-1}}(b^*)\eta, r^{-1}t)\big\rangle
 =  \sum_{l \in G}  \xi(ls)^* (\lambda^\B_{r^{-1}}(b^*)\eta)(lr^{-1}t)\\
 &= \sum_{l \in G}  \xi(ls)^* b^*\eta(rlr^{-1}t) = \sum_{k \in G}  \xi(r^{-1}krs)^* b^*\eta(kt)
 \end{align*}
 where we have substituted $l = r^{-1}kr$ to get the last equality (this is allowed as the sums are unconditionally convergent). Hence we get that $\big\langle \rho(b) (\xi, s) , (\eta, t)\big\rangle_\B =  \big\langle  (\xi, s) , \rho(b^*)(\eta, t)\big\rangle_\B$, as desired. 
 
 Finally, condition (iv) is also satisfied. Indeed, let $r, s, t \in G$,  $b\in B_r$, $c \in B_t$ and $\xi \in \ell^2(\B)$. Then we have that 
 \[ ((\lambda^\B_r(b)\xi)\cdot c)(k) = (\lambda^\B_r(b)\xi)(kt^{-1})\, c = b \,\xi (r^{-1}kt^{-1}) \,c = b \,(\xi\cdot c) (r^{-1}k) = 
 (\lambda^\B_r(b)(\xi\cdot c))(k)  \]
 for all $k\in G$, i.e., $(\lambda^\B_r(b)\xi)\cdot c = \lambda^\B_r(b)(\xi\cdot c)$. Hence we get that
 \begin{align*} (\rho(b) (\xi, s))\cdot c &= (\lambda^\B_r(b) \xi, rs) \cdot c = ((\lambda^\B_r(b) \xi)\cdot c, rst)\\
 &= (\lambda^\B_r(b)( \xi\cdot c), rst) = \rho(b) (\xi\cdot c, st) = \rho(b) ((\xi, s) \cdot c),   
 \end{align*} as desired.
 
 Pick now $\xi \in \ell^2(\B)$. For $g\in G$ and $b\in B_g$ we get that
\begin{align*} \big\langle (\xi, e) , \rho(b) (\xi, e)\big\rangle_\B &=  \langle (\xi, e) , (\lambda^\B_g(b) \xi, g)\rangle_\B =
\sum_{h\in G} \xi(h)^* (\lambda^\B_g(b)\xi)(hg) \\
&= \sum_{h\in G} \xi(h)^* b\,\xi(g^{-1}hg) =  \sum_{h\in G} \xi(gh)^* b\,\xi(hg).
\end{align*}
Hence we get from Proposition \ref{action-pd} that $T = (T_g)_{g\in G} $ given by 
\[T_g(b) = \sum_{h\in G} \xi(gh)^* b\,\xi(hg) \] 
$\text{for every } g\in G $ and $b \in B_g$
is a $\B$-positive definite bundle map. Note that this expression  for $T$ differs from the one for Exel-type $\B$-bundle maps in Example \ref{regB-action}. In fact, it is not difficult to see that every $\B$-bundle map of Exel-type has finite support (in the sense that 
the set $\{g\in G: T_g\neq 0\}$ is finite). But if $\B$ is unital and we let $\xi\in \ell^2(\B)$ be defined by $\xi(g) = 1_{B_e}$ if $g=e$ and $\xi(g) =0$ otherwise, one gets that $T_g= {\rm id}_{B_g}$ for every $g\in G$, i.e., the $\B$-bundle map $T=(T_g)$ does not have finite support, hence is not of Exel-type.  

 \end{example}

\subsubsection{$(\A,\varphi)$-actions on Hilbert $\B$-bundles} \label{A-phi-actions}

\begin{example}
Assume $G =\{e\}$ and $H=\{e\}$ are both trivial, so $\A$ and $\B$ have only one fibre, say $A$ and $B$, respectively. Let $\varphi_0$
 denote the trivial homomorphism.   A  Hilbert $\B$-bundle $\X=\{X\}$ is then just a  Hilbert $B$-module $X$, and an $(\A, \varphi_0)$-action $\rho$ on $\X$ is nothing but a left action of $A$ on $X$ by adjointable operators, i.e., a representation of $A$ in $\L_B(X)$.  
Moreover, an $\A$-$\varphi_0$-$\B$-bundle map $T$ is simply a bounded linear map $T:A\to B$. If $T$  is a diagonal coefficient of an action $\rho$ as above, i.e., for all $a\in A$ we have $T(a) = \langle x, \rho(a) x\rangle_B$ for some $x\in X$, then it is well-known
 that  $T$ is completely positive. This is in accordance with the fact that in the present situation an $\A$-$\varphi_0$-$\B$-bundle map is positive definite if and only if it is completely positive, cf.~\cite[Example 4.2]{BeCo7}.    
\end{example}

\begin{example} 
Assume now that only $H$ is trivial, so $\B$ has only one fiber $B$ over $e$, and let  $\varphi_0$ denote the trivial homomorphism from $G$ into $\{e\}$. 
Then  an $(\A, \varphi_0)$-action corresponds to  a $*$-representation of $\A=(A_g)_{g\in G}$ on a Hilbert $B$-module. 

Indeed, let $\rho$ be an $(\A,\varphi_0)$-action  on a Hilbert $\B$-bundle $\X$.  Again, $\X=\{X\}$ for some Hilbert $B$-module $X$. For each $g \in G$, let $\pi_g$ denote the restriction of $\rho$ to $A_g$. Then each $\pi_g:A_g\to \mathcal{F}(X)$ is a map such that $\pi_g(a): X\to X$ is linear and  
$\pi_{gg'}(aa') = \pi_g(a)\pi_{g'}(a')$ for all $g'\in G$ and $a\in A_g, a'\in A_{g'}$. Moreover, we also get
$\langle \pi_g(a) x,  y\rangle_B = \langle x,  \pi_{g^{-1}}(a^*)  y\rangle_B$ for all $a\in A_g$ and $x, y \in X$, i.e., 
$\pi_g(a)$ is an adjointable operator on $X$ with 
$\pi_g(a)^* =  \pi_{g^{-1}}(a^*) $ for every $a \in A_g$. This means that $(\pi_g)_{g\in G}$ is a $*$-representation of $\A$ in $\L_B(X)$. 

Conversely, assume that  $\pi= (\pi_g)_{g\in G}$ is a $*$-representation of $\A$ in $\L_B(X)$ for some Hilbert $B$-module $X$. Thinking of $\X=\{X\}$ as a Hilbert $\B$-bundle, 
we can define $\rho_\pi: \A \to \F(\X)$ by $\rho_\pi (a):= \pi_g(a)$ for every $g\in G$ and $a\in A_g$. Then it is immediate that $\rho_\pi$ satisfies conditions (i)-(iii) for being an $(\A,\varphi_0)$-action on $\X$. Moreover, condition (iv) is also satisfied: one has to check that $(\pi_g(a)x)b = \pi_g(a)(xb)$ for all $g\in G$, $a\in A_g$, $x\in X$ and $b\in B$, i.e., each $\pi_g(a)$ is $B$-linear. This holds since  each $\pi_g(a)$ is an adjointable operator on $X$. Thus $\rho_\pi$ is an $(\A,\varphi_0)$-action on $\X$.

Note that for every $x\in X$, we get a positive definite $\A$-$\varphi_0$-$\B$ bundle map $T^{\pi, x}$ given by 
\[ T^{\pi, x}_g(a) = \langle x, \rho_\pi(a) x\rangle_\B = \langle x, \pi_g(a) x\rangle_B\]
for all $g\in G$ and $a\in A_g$. 
The associated completely positive map $\Phi^{\pi,x}:C^*(\A) \to B$ is then given by 
\[  \Phi^{\pi,x}(f) = \sum_{g\in G} \langle x, \pi_g(f(g)) x\rangle_B = \langle x, \tilde\pi(f) x\rangle_B\]
for all $f \in C_c(G, \A)$, where $\tilde\pi$ denotes the $*$-representation of $C^*(\A)$ on $X$ associated to $\pi$.
Clearly, when $B=\Complessi$, $X$ is a Hilbert space and $\Phi^{\pi,x}$ is a positive linear functional on $C^*(\A)$.
\end{example}

\begin{example}  \label{actionrep}
Let $\pi = (\pi_g)_{g\in G}$ be a $*$-representation  of $\A=(A_g)_{g\in G}$ on some  Hilbert $B$-module $X$ (i.e., in $\L_B(X)$), where $B$ is a $C^*$-algebra.  Assume $\varphi$ is a homomorphism from $G$ into a group $H$. 
Let then $\B = (B \times \{h\})_{h\in H}$ be the natural Fell bundle over $H$ associated to $B$, with its obvious linear structure and norm, multiplication and adjoint being given by \[ (b, h)\cdot (c, h') = (bc, hh'), \, (b, h)^* = (b^*, h^{-1}).\] 
Note that $\B$ is the Fell bundle associated to the system $(B, H, {\rm id})$, where ${\rm id}$ denotes the trivial action of $H$ on $B$, cf.~Example \ref{eqrep1}. Hence, $C^*(\B) \simeq B \otimes_{\rm max} C^*(H)$, cf.~\cite[Proposition 16.28]{Exel2}.
 
 Further, let $\X = 
  (X \times \{h\})_{h\in H}$ 
  be the natural Hilbert $\B$-bundle we get by letting the maps $\X\times \B \to \X$ and $\X\times \X \to \B$ be defined by
\[ (x, h) \cdot (b, h') = (x \cdot b, hh'), \quad \langle (x, h), (y, h')\rangle_\B = (\langle x, y\rangle_B, h^{-1}h')\,.\]
Then, as the reader will have no problem to check,  
we get an $(\A, \varphi)$-action $\rho_{\pi, \varphi}$ on  $\X$ by setting 
\[ \rho_{\pi, \varphi}(a)(x, h) := (\pi_g(a)x, \varphi(g)h) \]
for all $g \in G$, $a\in A_g$, $x \in X$ and $h\in H$.  

If $x \in X$, the diagonal coefficient of $\rho_{\pi,\varphi}$ associated to $(x, e) \in X\times \{e\}$ is given by
\[ \langle (x, e), \rho_{\pi, \varphi}(a) (x, e)\rangle_\B = \langle (x, e), (\pi_g(a)x, \varphi(g))\rangle_\B = \big(\langle x, \pi_g(a) x\rangle_B, \varphi(g)\big)\]
for all $g\in G$ and $a\in A_g$.  Thus, the $\A$-$\varphi$-$\B$-bundle map $T^{\pi, \varphi, x}=(T^{\pi, \varphi, x}_g)_{g\in G}$ given by 
\[ T_g^{\pi, \varphi, x}(a) = \big(\langle x, \pi_g(a) x\rangle_B, \varphi(g)\big)\]
for all $g\in G$ and $a\in A_g$ is positive definite, and the associated completely positive map $\Phi^{\pi, \varphi, x} : C^*(\A) \to C^*(\B) \simeq B \otimes_{\rm max} C^*(H)$ is determined by
\[ \Phi^{\pi, \varphi, x}\big(\,\widehat{j}_g^\A(a)\big) = \langle x, \pi_g(a) x\rangle_B \otimes u^H_{\varphi(g)}\]
for all $g\in G$ and $a\in A_g$ (where $u^H$ denotes the universal unitary representation of $H$ into $C^*(H)$). If ${\rm ker}(\varphi)$ is amenable, one can similarly obtain a completely positive map from $C_r^*(\A)$ into $B\otimes_{\rm min} C_r^*(H)$. 
\end{example}

\begin{example}\label{group-bundle} Assume $U$ is a unitary representation of $G$ on some Hilbert space $X$
 and $\varphi \in {\rm Hom}(G,H)$.  Let $\A_G=(\Complessi\times \{g\})_{g\in G}$  and $\B_H=(\Complessi\times \{h\})_{h\in H}$ denote the canonical group bundles associated to $G$ and $H$. Then $U$ gives rise to an $(\A_G, \varphi)$-action $\rho$ on  
the natural Hilbert $\B_H$-bundle $\X=   (X \times \{h\})_{h\in H}$ 
given by
\[ [\rho((z, g))](x, h) :=  (zU_g(x), \varphi(g)h)\]
for all $(z, g) \in \A_g$ and $(x, h) \in \X$. Note that $\rho((1, e))$ is the identity operator on each fiber $X\times\{h\}$, i.e., $\rho$ is unital.  If $x \in X$, we get
\[ \big\langle (x, e) , [\rho((z,g))](x, e)\big\rangle_{\B_H} 
= ( \langle x, zU_g(x)\rangle, \varphi(g)) = z\,( \langle x,  U_g(x)\rangle, \varphi(g))
\]
for all $z \in \Complessi, g\in G$,
so Proposition \ref{action-pd} gives that $T=(T_g)_{g\in G}$ defined by 
\[T_g((z, g)) = z\,( \langle x,  U_g(x)\rangle, \varphi(g))\] is a positive definite $\A_G$-$\varphi$-$\B_H$ bundle map. 
The associated completely positive map $\Phi_T: C^*(G)\to C^*(H)$ satisfies that
$\Phi(u^G(g))=  \langle x,  U_g(x)\rangle \, u^H(\varphi(g))$
for all $g\in G$. If ${\rm ker}(\varphi)$ is amenable, one can similarly obtain a completely positive map $M$ from $C_r^*(G)$ into $C_r^*(H)$.

We also note that if $P$ is an orthogonal projection in $\B(X)$ commuting with the range of $U$, then we can  set
\[ [\rho_P((z, g))](x, h) :=  (zPU_g(x), \varphi(g)h)\]
for all $(z, g) \in \A_g$ and $(x, h) \in \X$ and obtain an $(\A_G, \varphi)$-action $\rho_P$, which is  unital only when $P= {\rm id}_X$.
\end{example}

\begin{example} \label{regA-action} Let $\A=(A_g)_{g\in G}$ and $\B=(B_h)_{h\in H}$ be Fell bundles over $G$ and $H$, respectively, and $\varphi \in {\rm Hom}(G, H)$. Assume that $\rho$ is an $(\A, \varphi)$-action on a  Hilbert $\B$-bundle $\X=(X_h)_{h\in H}$.   We may then form the \emph{regular $(\A,\varphi)$-action $\widecheck \rho$ associated to $\rho$} as follows. 

Let  $\widecheck \X=(\widecheck X_h)_{h\in H}$ be the regular Hilbert $\B$-bundle associated to $\X$, cf.~Example \ref{reg-bundle}.
Recall that $\widecheck \X$ is obtained by completing the inner product $\B$-bundle $\widecheck \X^0=(\widecheck X_h^0)_{h\in H}$, where $\widecheck X_h^0= C_c(H, X_h)$ for each $h\in H$. Let  $g \in G$ and $a\in A_g$, and define $\widecheck\rho_0(a): \widecheck \X^0 \to \widecheck \X^0$ by setting for each $h \in H$ and $\xi \in \widecheck X_h^0$
 \[ (\widecheck\rho_0(a) \xi) (h') := \rho(a)\,\xi\big(\varphi(g)^{-1}h'\big) \quad \text{for all } h'\in H.\]
 Then $\widecheck\rho_0$ is an $(\A, \varphi)$-pre-action on $\widecheck\X^0$. Indeed, conditions (i)--(iv) in Definition \ref{BC-action-def} follow for $\widecheck\rho_0$ by routine computations. As a sample, we check that condition (iii) holds for $\widecheck\rho_0$. Let 
  $g\in G$, $a\in A_g$ and $\xi, \eta \in \widecheck \X^0$. Then we have
  \begin{align*}
 \big\langle \widecheck\rho_0(a) \xi , \eta \big\rangle_\B &= \sum_{h\in H} \langle \rho(a)\,\xi(\varphi(g)^{-1}h) ,\eta(h) \rangle_\B = \sum_{h\in H} \langle\,\xi(\varphi(g)^{-1}h) ,\rho(a^*)\eta(h) \rangle_\B \\
& = \sum_{h\in H} \langle\,\xi(h) ,\rho(a^*)\eta(\varphi(g)h)\rangle_\B
= \sum_{h\in H} \langle\,\xi(h) ,\rho(a^*)\eta(\varphi(g^{-1})^{-1}h)\rangle_\B\\
 &=  \big\langle  \xi , \widecheck\rho_0(a^*)\eta \big\rangle_\B.
 \end{align*}
 Moreover, using  Lemma \ref{actionineq}, we also have
 \begin{align*}
 \langle \widecheck\rho_0(a) \xi , \widecheck\rho_0(a)\xi \rangle_\B &
 = \sum_{h'\in H} \,\langle(\widecheck\rho_0(a)\xi)(h'),  (\widecheck\rho_0(a)\xi)(h')\rangle_\B
 = \sum_{h'\in H} \,\langle\rho(a)\xi(\varphi(g)^{-1}h'),  \rho(a)\xi(\varphi(g)^{-1}h')\rangle_\B\\ 
 &= \sum_{k\in H} \,\langle\rho(a)\xi(k),  \rho(a)\xi(k)\rangle_\B
 \leq \sum_{k\in H} \|a\|^2 \langle \xi(k),\xi(k)\rangle_\B = \|a\|^2 \,\langle \xi, \xi\rangle_\B.
\end{align*}
Hence, $\widecheck\rho_0$ induces an $\A$-$\varphi$-$\B$ action $\widecheck\rho$ on the completion $\widecheck\X$ of $\widecheck\X^0$, cf.~Remark \ref{pre-act-def}.  
\end{example}

 \begin{example} \label{condexp-Hilbund}
 Consider   a Fell sub-bundle $\B = (B_g)_{g\in G}$  of a Fell bundle $\A=(A_g)_{g\in G}$ and assume $E= (E_g)_{g\in G}$ is a conditional expectation from $\A$ to $\B$. We may then regard $\A$ as a semi-inner product $\B$-bundle, as in Example \ref{sub-b}. It is easy to check that the map $\rho_0:\A \to \F(\A)$ given by $\rho_0(a) a':= aa'$ for $a, a'\in \A$ is then an $\A$-pre-action on  $\A$. 
 Let  $x \in A_e$. Then for $g\in G$ and $a\in A_g$ we have \[ \langle x, \rho_0(a) x\rangle_\B =E_g (x^*ax).\]
 Thus Remark \ref{pre-act-pd} gives that $T=(T_g)_{g\in G}$ defined by $T_g(a) = E_g(x^*ax)$ for all $g\in G$ and $a\in A_g$ is a positive definite $\A$-$\B$ bundle map. This fits with \cite[Example 4.7]{BeCo7}, where this was observed for $T=E$. 
\end{example}

\subsection{Actions and $C^*$-dynamical systems} \label{action-ds}
In the study of imprimitivity bimodules over $C^*$-crossed products and of Morita equivalence for $C^*$-dynamical systems, the notion of compatible action of a group has played a useful role, see e.g.~\cite{EKQR00, EKQR06} and references therein. Also, the closely related concept of equivariant representations of discrete unital $C^*$-dynamical systems is an essential tool in the study of Fourier-Stieltjes algebras associated to such systems \cite{BeCo6}. We introduce below a more general notion of compatibility of actions for $C^*$-dynamical systems and show that it yields actions of the associated Fell bundles.

\medskip Let $\A_\Sigma =(A_g)_{g\in G}$ and $\B_\Omega=(B_h)_{h\in H}$ denote the canonical Fell bundles associated to two discrete $C^*$-dynamical systems $\Sigma= (A,  G, \alpha)$ and $\Omega=(B, H, \beta)$, and  let $\varphi \in {\rm Hom}(G, H)$. 
\begin{definition}
By an \emph{$(\alpha,\beta, \varphi)$-compatible action of  $G$ on an $A$-$B$ correspondence $X$}, we will mean a homomorphism $\gamma$  from $G$ into the group $\mathcal{I}(X)$ of $\Complessi$-linear invertible isometries of $X$ satisfying
\begin{itemize}
\item[(i)] $\gamma_g(a\cdot x) = \alpha_g(a)\cdot\gamma_g(x)$
\item[(ii)] $\gamma_g(x\cdot b) = \gamma_g(x)\cdot \beta_{\varphi(g)}(b)$
\item[(iii)] $\big\langle\gamma_g(x),\gamma_g(y)\big\rangle_B = \beta_{\varphi(g)}\big(\langle x, y\rangle_B\big)$
\end{itemize}
for all $g\in G, a\in A, b\in B$ and $x, y \in X$. 
\end{definition}

Note that when $G=H$ and $\varphi={\rm id}_G$, then $\gamma$  is called an $\alpha$-$\beta$ compatible action of $G$ on the  Hilbert $A$-$B$-bimodule $X$ in \cite{EKQR06} (where it is also required that  the left action of $A$ on $X$ is non-degenerate). 
Note also that if $A$ is unital, $\Omega=\Sigma$, and $(\varrho, v)$ is an equivariant representation of $\Sigma$ on a  Hilbert $A$-module in the sense of \cite{BeCo6}, then $v$ can be considered as an $\alpha$-$\alpha$ compatible action of $G$ on $X$ by setting $a\cdot x= \varrho(a)x$ for $a\in A$ and $x \in X$. 

\begin{proposition} \label{compatible} Let $\gamma$ be an $(\alpha,\beta, \varphi)$-compatible action of  $G$ on an $A$-$B$ correspondence $X$. 
Let  $\X=(X_h)_{h\in H}$ be the natural Hilbert $\B_\Omega$-bundle associated to $X$ and $\Omega$, cf.~Example \ref{eqrep1}. 

Then the map $\rho_\gamma: \A_\Sigma \to \F(\X)$, defined by 
\[ [\rho_\gamma((a, g))](x,h) = (a\cdot\gamma_g(x), \varphi(g)h)\]  
for all $ a\in A, g\in G, x\in X$ and $h\in H$,
is an $(\A_\Sigma, \varphi)$-action on the Hilbert $\B_\Omega$-bundle $\X$. 

Hence, if $x\in X$ and $T=(T_g)_{g\in G}$ is defined by 
\[T_g((a, g)) = (\langle x, a\cdot\gamma_g(x)\rangle_B , \varphi(g))\]
for all $a\in A$ and $g\in G$, then $T$ is a positive definite $\A_\Sigma$-$\varphi$-$\B_\Omega$ bundle map.

\end{proposition}

\begin{proof}
Property (i) in Definition \ref{BC-action-def} is clearly satisfied.
For completeness, we check properties (ii)-(iv).
First, we have
\begin{align*}
[\rho_\gamma((a, g)(a',g'))](x,h) & = [\rho_\gamma((a\alpha_g(a'),gg'))](x,h) \\
& =  (a \alpha_g(a') \cdot \gamma_{gg'}(x),\varphi(gg')h) \\
& =  (a \alpha_g(a') \cdot \gamma_{g}(\gamma_{g'}(x)),\varphi(gg')h) \\
& =  (a \cdot \gamma_{g}(a' \cdot \gamma_{g'}(x)),\varphi(g)\varphi(g')h) \\
& = \rho_\gamma((a, g))\Big( (a'\cdot\gamma_{g'}(x), \varphi(g')h)\ \Big) \\
& = \rho_\gamma((a, g))\Big([\rho_\gamma((a', g'))](x,h)\Big)
\end{align*}
for all $a,a' \in A$, $g,g' \in G$, $h \in H$ and $x \in X$, which is property (ii).
Next we show that property (iii) holds. Let $(x, h), (x', h') \in \X$ and $(a, g) \in \A_\Sigma$. Then we have
\begin{align*}
\big\langle [\rho_\gamma((a,g))](x,h), (x',h')\big\rangle_{\B_\Omega} & = \big\langle(a\cdot\gamma_g(x), \varphi(g)h), (x',h')\big\rangle_{\B_\Omega}\\
&= \Big( \beta_{\varphi(g)h}^{-1}\big(\langle a\cdot\gamma_g(x), \,x'\rangle_B\big), (\varphi(g)h)^{-1}h' \Big)\\
 &= \Big( \beta_{h}^{-1} \beta_{\varphi(g^{-1})}\big(\langle a\cdot\gamma_g(x), \,x'\rangle_B\big), (\varphi(g)h)^{-1}h' \Big)\\
 &= \Big( \beta_{h}^{-1} \big(\langle \gamma_{g^{-1}}(a\cdot\gamma_g(x)), \,\gamma_{g^{-1}}(x')\rangle_B\big), (\varphi(g)h)^{-1}h' \Big)\\
 &= \Big( \beta_{h}^{-1} \big(\langle \alpha_{g^{-1}}(a)\cdot x, \,\gamma_{g^{-1}}(x')\rangle_B\big), (\varphi(g)h)^{-1}h' \Big)\\
&=\Big( \beta_h^{-1} \big(\langle x, \alpha_{g^{-1}}(a^*)\cdot \gamma_{g^{-1}}(x')\rangle_B\big), h^{-1}\varphi(g)^{-1}h'\Big)\\
&= \big\langle (x,h), (\alpha_{g^{-1}}(a^*)\cdot \gamma_{g^{-1}}(x'), \varphi(g^{-1})h')\big\rangle_{\B_\Omega}\\
&= \big\langle (x,h), [\rho_\gamma((\alpha_{g^{-1}}(a^*), g^{-1})](x',h')\big\rangle_{\B_\Omega}\\
& = \big\langle (x,h), [\rho_\gamma((a,g)^*)](x',h')\big\rangle_{\B_\Omega}
\end{align*} 
as expected.
Finally, for $(b, k) \in \B_\Omega$, we also have
\begin{align*}
\Big([\rho_\gamma((a,g))](x,h)\Big)\cdot(b, k)& = (a\cdot\gamma_g(x), \varphi(g)h)\cdot(b,k)\\
& = \big(a\cdot\gamma_g(x)\cdot\beta_{\varphi(g)h}(b), \varphi(g)h k\big)\\
&= \big(a \cdot \gamma_g(x)\cdot \beta_{\varphi(g)}(\beta_h(b)), \varphi(g)hk\big)\\
&= \big(a \cdot \gamma_g(x\cdot\beta_h(b)), \varphi(g)hk\big)\\
&= [\rho_\gamma((a,g))]\big((x\cdot\beta_h(b), hk)\big)\\
&=[\rho_\gamma((a,g))]\big((x,h)\cdot(b, k)\big) \ , 
\end{align*}
which means that property (iv)  is satisfied as well.

Finally, if $x \in X$, then we have
\[ \big\langle (x, e) , [\rho_\gamma((a,g))](x, e)\big\rangle_{\B_\Omega} 
= ( \langle x, a\cdot\gamma_g(x)\rangle_B , \varphi(g))
\]
for all $a \in A$ and $g\in G$, so the last assertion follows form Proposition \ref{action-pd}.
\end{proof}

\bigskip 
\begin{example} \label{eqrep-Baction} 
Let us have a closer look at the case
where $G=H$, $\varphi={\rm id}_G$, $\Sigma =  \Omega = (A, G, \alpha)$ and $A$ is unital, 
and $(\varrho, v)$ is an equivariant representation of $\Sigma$ on a  Hilbert $A$-module $X$ (as defined in \cite{BeCo6}, see also \cite{BeCo3, BeCo4}).
By Proposition \ref{compatible}, 
 we get a $\B_\Sigma$-action  $\rho_{(\varrho, v)}$ on 
 the  Hilbert $\B_\Sigma$-bundle  $\X=(X_g)_{g\in G}$ associated to $X$, given by
 \[\rho_{(\varrho,v)}((a, g)) (x, h) = \big(\varrho(a)v(g) x, gh\big)\] for every $(a, g) \in \B_\Sigma$ and $(x, h) \in \X$. 
Further, if $x \in X$,  the associated positive-definite $\B_\Sigma$-bundle map $T:= (T_g)_{g\in G}$ takes the form 
\[ T_g ((a, g)) = 
\big(\langle x, \varrho(a)v(g)x\rangle_A, g\big)\]
for every $g\in G$ and $a\in A$. 
 When $(\varrho, v)=(\ell,\alpha)$
is the trivial equivariant representation of $\Sigma$ on $A$ (so $\ell: A \to \L_A(A)$ is given by left multiplication), then it is immediate that $\rho_{(\ell, \alpha)}$ coincides with the trivial $\B_\Sigma$-action $\rho_{\B_\Sigma}$.

Assume now that $(\varrho, v)=(\widecheck\ell,\widecheck\alpha)$ is the  regular equivariant representation of $\Sigma$ on the Hilbert $A$-module $X= A^G$ (cf.~\cite[Example 4.7]{BeCo3}). We recall that \[A^G = \big\{ \xi:G\to A: \sum_{g\in G} \xi(t)^*\xi(t) \text{ is norm-convergent  in } A\big\}\] is equipped with operations given by $(\xi\cdot a)(t) = \xi(t)a$ and  $\langle \xi, \eta\rangle_A = \sum_{t\in G}\xi(t)^*\eta(t)$ for $\xi, \eta \in A^G, a\in A$ and  $t\in G$. Moreover,  $\widecheck\ell$ and $\widecheck\alpha$ are given by \[ (\widecheck\ell(a)\xi)(s) =  a \,\xi(s) \text{  and } (\widecheck\alpha(r) \xi)(s) = \alpha_r(\xi(r^{-1}s))\] for $a\in A, \xi \in A^G$ and $r, s \in G$, so we have
\[(\widecheck\ell(a)\widecheck\alpha(r)\xi)(s) = a (\widecheck\alpha(r)\xi)(s) = a  \alpha_r(\xi(r^{-1}s)).\]
The   Hilbert $\B_\Sigma$-bundle $\X=(X_r)_{g\in G}$ associated to $X=A^G$ is then given by for each $r\in G$ by $X_r= A^G \times \{r\}$, with operations given by
\[ (\xi, r)\cdot (a, s) = (\xi \cdot\alpha_r(a), rs), \quad \langle (\xi, r), (\eta, s)\rangle_{\B_\Sigma} = \big(\alpha_{r^{-1}}(\langle \xi, \eta\rangle_A), r^{-1}s\big),\] 
for $\xi, \eta \in A^G, a\in A$ and $s\in G$.
The corresponding action $\rho_{(\widecheck\ell,\widecheck\alpha)}$ of $\B_\Sigma$ on $\X$ is then given by 
\[ \rho_{(\widecheck\ell,\widecheck\alpha)}((a, g)) (\xi, r) = ( \varrho(a)v(g)\xi, gr).\]

On the other hand,  according to Example \ref{regB-action}, we can also form the regular $\B_\Sigma$-action 
$\widecheck{\rho_{\B_\Sigma}}$ on the regular  Hilbert $\B_\Sigma$-bundle $\widecheck{\B_\Sigma}$,
which is the completion of the inner product 
$\B_\Sigma$-bundle $ \big(C_c(G, A \times \{r\})\big)_{r \in G}$. Identifying $C_c(G, A \times \{r\})$
with $C_c(G, A) \times \{r\}$ (using the identification given by $(\xi, r) (t) = (\xi(t), r)$ for all $t\in G$), the operations are given by
\[ (\xi, r) \cdot (a, s) = (\xi \cdot \alpha_r(a), rs)\] 
for $\xi \in C_c(G,A), a\in A$ and $r, s\in G$, 
where $(\xi \cdot a')(t) = \xi(t) a'$ for all $a' \in A$ and $t\in G$, and 
\begin{align*} \langle (\xi, r), (\eta, s) \rangle_{\B_\Sigma} &= \sum_{t\in G} (\xi(t), r)^*(\eta(t), s) =  \sum_{t\in G} (\alpha_{r^{-1}}(\xi(t)^*), r^{-1})(\eta(t), s)\\
&=  \sum_{t\in G} \Big(\alpha_{r^{-1}}\big(\xi(t)^*)\eta(t)\big), r^{-1}s\Big) = \Big(\alpha_{r^{-1}}\Big(\sum_{t\in G}\xi(t)^*\eta(t)\Big),  r^{-1}s\Big)\\ &= \big( \alpha_{r^{-1}}\big( \langle \xi, \eta \rangle_{A}\big), r^{-1}s\big)
\end{align*}
for $\xi, \eta \in C_c(G,A)$ and $r, s\in G$. It is then clear that 
for each $r \in G$ the fiber of $\widecheck{\B_\Sigma}$ over $r$ can be identified with $A^G \times \{r\}$, with operations given as above, but now with $\xi, \eta \in A^G$. Thus we get that 
$\X = \widecheck{\B_\Sigma}$
(up to the above identification). To be on the formal side, 
$\X$ and  $\widecheck{\B_\Sigma}$ are unitarily equivalent in the sense of Definition  \ref{unitary-eq}).

Now, for all $a \in A$, $g,r,s \in G$ and $\xi \in A^G$, we have
\begin{align*} \big(\widecheck{\rho_{\B_\Sigma}}
((a,g)) (\xi,r) \big)(s)&=  (a,g) ((\xi, r)(g^{-1}s)) = (a,g) (\xi(g^{-1}s),r) \\
&= (a \alpha_g(\xi(g^{-1}s)),gr) =  ((\varrho(a)v(g)\xi)(s), gr)  \\
&=  (\varrho(a)v(g)\xi, gr)(s)\\
&= \big(\rho_{(\widecheck\ell,\widecheck\alpha)}((a,g)) (\xi, r)\big)(s)
\end{align*}
Thus it follows that the regular action $\widecheck{\rho_{\B_\Sigma}}$ associated to $\rho_{\B_\Sigma}=\rho_{(\ell, \alpha)}$ coincides with the action $\rho_{(\widecheck\ell,\widecheck\alpha)}$ associated to the regular equivariant representation $(\widecheck\ell,\widecheck\alpha)$. Formally, these two actions are unitarily equivalent, in the sense of Definition \ref{unitary-eq2}.
Actually, by similar arguments,
one can show that for any equivariant representation $(\varrho,v)$ of $\Sigma$ on $X$, 
the regular $\B_\Sigma$-action associated to $\rho_{(\varrho,v)}$ is unitarily equivalent to the $\B_\Sigma$-action $\rho_{(\widecheck{\varrho},\widecheck{v})}$
associated to the {induced regular equivariant representation} $(\widecheck{\varrho},\widecheck{v})$ of $\Sigma$ on the Hilbert $A$-module $X^G$
(as defined in \cite[Example 4.8]{BeCo3}).
 \end{example}

\begin{example} Let $\Sigma=(A,G,\alpha)$ be a unital $C^*$-dynamical system and $X_1$ and $X_2$ be Hilbert $A$-modules. Let then $\X_1=(X_1\times \{g\})_{g\in G}$ (resp.~$\X_2=(X_2\times \{g\})_{g\in G}$) denote the canonical Hilbert $\B_\Sigma$-bundle associated to $X_1$ (resp.~$X_2$). If $X_1$ and $X_2$ are unitarily equivalent, then it is straightforward to check that 
$\X_1$ and $\X_2$ are unitarily equivalent.  

Moreover, assume  $(\varrho_1,v_1)$ (resp.~$(\varrho_2,v_2)$) is an equivariant representation of $\Sigma$ on $X_1$ (resp.~$X_2$), 
and $(\varrho_1,v_1)$ is unitarily equivalent to $(\varrho_2,v_2)$, i.e., there exists a unitary operator $u:X_1\to X_2$ such that $\rho_2(a) = u\rho_1(a)u^*$ and $v_2(g)=uv_1(g)u^*$ for all $a\in $ and $g\in G$.  Then the $\B_\Sigma$-action $\rho_{(\varrho_1,v_1)}$ associated to $(\varrho_1,v_1)$ on 
$\X_1$ is easily seen to be unitarily equivalent to the $\B_\Sigma$-action $\rho_{(\varrho_2,v_2)}$ associated to $(\varrho_2,v_2)$ on 
$\X_2$ via the $\B_\Sigma$-bundle map $U=(U_g)_{g\in G}$ from $\X_1$ to $\X_2$ given by $U_g(x_1, g)= (ux_1, g)$ for all $g\in G$ and $x_1\in X_1$. 
\end{example}

\section{Actions on Hilbert bundles from positive definite bundle maps}\label{GR-PD}
In this section, $\A=(A_g)_{g \in G}$ and $\B=(B_h)_{h \in H}$ are Fell bundles over discrete groups $G$ and $H$, respectively, and $\varphi \in {\rm Hom}(G, H)$. 
Our main goal  is to prove a Gelfand-Raikov type theorem for positive definite $\A$-$\varphi$-$\B$ bundle maps 
in the case where
 $\A$ and $\B$ are assumed to be unital, i.e., $A_e$ and $B_e$ are unital. 
 We show that any such map $T$ gives rise to an $(\A, \varphi)$-action $\rho$ on a Hilbert $\B$-bundle $\X$ such that $T$ can be written as a diagonal coefficient of $\rho$ over the unit fibre of $\X$, thereby showing that the converse of Proposition \ref{action-pd} holds. 
 
\begin{theorem} \label{GR}
 Let $\A=(A_g)_{g \in G}$ and $\B=(B_h)_{h \in H}$ be unital Fell bundles,  $\varphi \in {\rm Hom}(G, H)$ and 
 $T=(T_g)_{g\in G}$ be a positive definite $\A$-$\varphi$-$\B$ bundle map. Then there exist  a  Hilbert $\B$-bundle $\X=(X_h)_{h\in H}$, an $(\A, \varphi)$-action $\rho$ on $\X$, and $\xi\in X_e$, such that 
\[T_g(a) = \langle \xi, \rho(a) \xi \rangle_\B\]
for all $g \in G$ and $a \in A_g$. 

\end{theorem}
\begin{remark}\label{cyclic1}
By construction, it will become apparent that the vector $\xi\in X_e$ is \emph{cyclic} (for $\rho$) in the sense that for each $r \in H$, ${\rm span}\{(\rho(a)\xi)b: g\in G, a \in A_g, h\in H, b\in B_h, \varphi(g)h=r\}$ is dense in $X_r$.  
\end{remark}

We will present two distinct proofs of Theorem \ref{GR}. The first one deals only with the case where $G=H$ and $\varphi = {\rm id}_G$.
The idea is to use Theorem \ref{PDCP} 
to promote a given positive definite $\A$-$\B$ bundle map $T$ to a completely positive map $\Phi_T: C^*(\A)\to C^*(\B)$, apply Paschke's GNS-type result in \cite{Pas} to $\Phi_T$, and construct $\X$ and $\rho$ from the resulting GNS-triple $(X, \Psi, \xi)$.     
The second one addresses the general case.
 We only sketch the proof, as it follows the same abstract pattern as in the proofs of other analogous results in the literature.

\subsection{Proof of Theorem \ref{GR} for a  positive definite $\A$-$\B$ bundle map} \label{GR-1} 
We assume here that $H=G$, $\varphi={\rm id}_G$,  $\A=(A_g)_{g \in G}$ and $\B=(B_g)_{g \in G}$ are unital, and  $T$ is a positive definite $\A$-$\B$ bundle map.
We recall that Paschke's GNS-type theorem for completely positive maps, cf.~\cite[Theorem 5.2]{Pas}, says that
 that if $C$ and $ D$ are unital $C^*$-algebras and $\Phi:C\to D$ is a completely positive linear map, then there exists a  Hilbert module $X$ over $D$, a unital $*$-homomorphism $\Psi: C \to \L_D(X)$ and  $\xi \in X$ such that $\Phi(c) = \langle \xi, \Psi(c) \xi\rangle_D$ for all $c\in C$ and ${\rm span}\{ (\Psi(c)\xi)d: c\in C, d \in D\}$ is dense in $X$.  

\medskip 
Now, using Theorem \ref{PDCP},
we can form the (unique) completely positive linear map $\Phi_T: C^*(\A) \to C^*(\B)$ satisfying that \[\Phi_T(\widehat{j}^\A_g(a)) = \widehat{j}^\B_g(T_g(a))\] for all $g \in G, a \in A_g$.  
Since $C^*(\A)$ and $C^*(\B)$ are unital (with units 
$1_\A:= \widehat{j}^\A_e(1_{A_e})$ and $1_\B:=  \widehat{j}^\B_e(1_{B_e})$), 
we may apply Paschke's theorem to $\Phi_T$,
giving a  Hilbert $C^*(\B)$-module $X$, a unital $*$-homomorphism $\Psi: C^*(\A) \to {\mathcal L}_{C^*(\B)}(X)$ 
and $\xi \in X$ such that 
\[\Phi_T(y) = \big\langle \xi, \Psi(y) \xi \big\rangle_{C^*(\B)}, \quad y \in C^*(\A). \]
We will denote the norm on $X$ by $\|\cdot\|_X$, which is given by $\|\eta\|_X=\|\langle \eta, \eta\rangle_{C^*(\B)}\|^{1/2}$.

For $g \in G$, define 
\[ X^0_g  = {\rm span} \big\{ (\Psi(\widehat{j}^\A_h(a))\xi) \cdot \widehat{j}^\B_k(b) \ | \ (a,b) \in A_h \times B_k, h, k \in G, hk=g  \big\},  \]
and let $X_g$ denote the closure of $X^0_g$ in $X$, so $X_g$ is a closed subspace of $X$, hence a Banach space. 

We set $\X := (X_g)_{g\in G}$ and denote also the disjoint union of the $X_g$'s by $\X$. We can  now define maps 
$\X\times \B \to \X, (z,c)\mapsto z\cdot c$, and $\X\times \X \to \B, (z, z')\mapsto \langle z, z'\rangle_{\B}$, as follows.

First, for $g, l \in G$, $z \in X_g$ and $c \in B_l$, we set
\[ z \cdot c := z\cdot\widehat{j}^\B_{l}(c)\,.\] 
Note that if $z \in X^0_g$, then $z \cdot c \in X^0_{gl}$, because
\[ \Big(\Psi(\widehat{j}^\A_h(a))\xi) \cdot \widehat{j}^\B_k(b) \Big)\cdot \widehat{j}^B_{l}(c) = (\Psi(\widehat{j}^\A_h(a))\xi) \cdot \widehat{j}^\B_{kl}(bc) \in X^0_{hkl} = X^0_{gl} \]
whenever $(a,b) \in A_h \times B_k, h, k \in G$ and $ hk=g$. 
Since the map $\eta\mapsto \eta\cdot\widehat{j}^\B_{l}(c)$ is continuous on $X$, it follows readily that $z \cdot c \in X_{gl} $ whenever $z \in X_g$ and $c \in B_l$.

\smallskip 
Next, for $g, g' \in G$ and $z \in X_g, z'\in X_{g'}$, we set 
\[ \langle z, z'\rangle_\B:= (\widehat{j}^\B_{g^{-1}g'})^{-1}\big( \langle z, z'\rangle_{C^*(\B)}\big),\] 
which is well-defined and belongs to $B_{g^{-1}g'}$. 
Indeed, consider 
$(a,b) \in A_h \times B_k, h, k \in G, hk=g$ and $(a',b') \in A_{h'} \times B_{k'}, h', k' \in G, h'k'=g'$.   
Note that 
\begin{align*}
\big\langle (\Psi(\widehat{j}^\A_h(a))\xi) \cdot \widehat{j}^\B_k(b), (\Psi(\widehat{j}^\A_{h'}(a'))\xi) \cdot \widehat{j}^\B_{k'}(b') \big\rangle_{C^*(\B)}
 & = \widehat{j}^\B_k(b)^* \big\langle \xi, \Psi(\widehat{j}^\A_h(a))^* \Psi(\widehat{j}^\A_{h'}(a')) \xi \big\rangle_{C^*(\B)} \,\widehat{j}^\B_{k'}(b') \\
& = \widehat{j}^\B_{k^{-1}}(b^*) \big\langle \xi, \Psi(\widehat{j}^\A_{h^{-1}h'}(a^* a')) \xi \big\rangle_{C^*(\B)} \,\widehat{j}^\B_{k'}(b') \\
& = \widehat{j}^\B_{k^{-1}}(b^*) \,\Phi_T\big(\widehat{j}^\A_{h^{-1}h'}(a^* a')\big) \,\widehat{j}^\B_{k'}(b') \\
& = \widehat{j}^\B_{k^{-1}}(b^*) \,\widehat{j}^\A_{h^{-1}h'}\big(T_{h^{-1}h'}(a^* a')\big)\, \widehat{j}^\B_{k'}(b') \\
& = \widehat{j}^\B_{(hk)^{-1}h'k'}\big(b^*T_{h^{-1}h'}(a^* a')b'\big),
 \end{align*}
 which gives that \[\big\langle (\Psi(\widehat{j}^\A_h(a))\xi) \cdot \widehat{j}^\B_k(b), (\Psi(\widehat{j}^\A_{h'}(a')\xi) \cdot \widehat{j}^\B_{k'}(b') \big\rangle_{C^*(\B)} \in \widehat{j}^\B_{(hk)^{-1}h'k'}\big(B_{(hk)^{-1}h'k'}\big) = \widehat{j}^\B_{g^{-1}g'}\big(B_{g^{-1}g'}\big).\] 
 By density of $X^0_g$ in $X_g$, continuity of $\langle \cdot, \cdot\rangle_{C^*(\B)}$ in each variable and the fact that $\widehat{j}^\B_{g^{-1}g'}$ is isometric, one readily gets that   $\langle z, z'\rangle_{C^*(\B)} \in \widehat{j}^\B_{g^{-1}g'}\big(B_{g^{-1}g'}\big)$, hence that 
 $\langle z, z'\rangle_{\B}$ is well-defined and belongs to $B_{g^{-1}g'}$, as claimed above. 
 We also record that  this computation shows that
\begin{equation} \label{T-eq}
\big\langle (\Psi(\widehat{j}^\A_h(a))\xi) \cdot \widehat{j}^\B_k(b), (\Psi(\widehat{j}^\A_{h'}(a'))\xi) \cdot \widehat{j}^\B_{k'}(b') \big\rangle_\B
=  b^* T_{h^{-1}h'}(a^*a') b' \in B_{g^{-1}g'}.
\end{equation}
It is now a routine matter to check that $\X = (X_g)_{g\in G}$ is a  Hilbert $\B$-bundle w.r.t.~these operations, after noticing  
that for $z \in \X$ we have
\[ \big(\|\langle z, z\rangle_\B\|_{B_e}\big)^{1/2} =  \big(\|(\widehat{j}^\B_e)^{-1}\big(\langle z, z\rangle_{C^*(\B)}\big)\|_{B_e}\big)^{1/2}
=  \|\langle z, z\rangle_{C^*(\B)}\|^{1/2} = \|z\|_X\]
(since $\widehat{j}^\B_e: B_e \to C^*(\B)$ is isometric).

Further, we can define an $\A$-action 
$\rho$ on $\X$ as follows. 
Let $g, g' \in G$ and $a \in A_g$. Then, since  
\[\Psi(\widehat{j}^\A_g(a)) \Big( \Psi(\widehat{j}^\A_{h'}(a'))\xi) \cdot \widehat{j}^\B_{k'}(b') \Big) = \Big( \Psi(\widehat{j}^\A_{gh'}(aa'))\xi \Big) \cdot \widehat{j}^\B_{k'}(b') \in X^0_{gg'}\]
whenever 
$(a',b') \in A_{h'} \times B_{k'}, h', k' \in G, h'k'=g'$, 
we readily get that  $\Psi(\widehat{j}^\A_g(a))$ maps $X_g$ into $X_{gg'}$. We can therefore set
\[\rho(a)z' := \Psi(\widehat{j}^\A_g(a))z'\] for all $z'\in X_{g'}$. 

It is now straightforward to verify that $\rho: \A \to \F(\X)$ is an $\A$-action on $\X$. 
For example, if $g\in G$, $a \in A_g$
 $z \in X_r$ and $w\in X_s$, then $\rho(a)z \in X_{gr}$ and $\rho(a^*)w \in X_{g^{-1}s}$, so
\begin{align*}
\langle \rho(a)z, w\rangle_\B & = (\widehat{j}^\B_{(gr)^{-1}s})^{-1}\big(\langle  \Psi(\widehat{j}^\A_g(a))z, w\rangle_{C^*(\B)}\big) 
= (\widehat{j}^\B_{(gr)^{-1}s})^{-1}\big(\langle  z, \Psi(\widehat{j}^\A_g(a))^*w\rangle_{C^*(\B)}\big) \\
& = (\widehat{j}^\B_{r^{-1}g^{-1}s})^{-1}\big(\langle  z, \Psi(\widehat{j}^\A_{g^{-1}}(a^*))w\rangle_{C^*(\B)}\big) =  \langle z,  \rho(a^*)w\rangle_\B.
\end{align*}
This shows that property (iii) for $\rho$ to be an $\A$-action holds. 

\smallskip Finally, since $A_e$ and $B_e$ are unital, we have $\xi =  (\Psi(\widehat{j}^\A_e(1_{A_e})) \xi)\cdot \widehat{j}^\B_e(1_{B_e})  \in X^0_e$. 
Using (\ref{T-eq}), we get that
\[\langle \xi, \rho(a) \xi \rangle_\B = \langle \xi, \Psi(\widehat{j}^\A_g(a)) \xi \rangle_\B = T_g(a)\]
for all $g \in G$ and $a \in A_g$, as desired.

\begin{remark}\label{non-unital}
 Note that if $\B$ is not unital,
 then we may 
 consider $\Phi_T$ as a completely positive map from $C^*(\A)$ into the unitization of $C^*(\B)$. Up to some minor modifications of the proof 
 presented above, it is not difficult to see that the conclusion of Theorem \ref{GR} (when $G=H$ and $\varphi={\rm id}_G$) still holds if we only assume $\A$ to be unital.   
Now, let us assume that $\A$ is not unital.
Then we may 
apply Paschke's theorem to the ``natural'' completely positive extension $\widetilde{\Phi_T}$ of $\Phi_T$, defined on the unitization of $C^*(\A)$, 
taking values in the unitization of $C^*(\B)$ (whenever needed), and satisfying that $\|\widetilde{\Phi_T}\| = \| \Phi_T\|$ 
(cf.~\cite[Proposition 2.2.1 and Remark 2.2.2]{BrOz}). Arguing 
(mutatis mutandis) as above, we may still construct a  Hilbert module $X$ over the unitization of $C^*(\B)$,  a vector $\xi \in X$,  a  Hilbert $\B$-bundle $\X=(X_g)_{g\in G}$ such that each $X_g$ is a closed subspace of $X$, and an $\A$-action $\rho$ on $\X$.
 However, it is not clear that we can recover $T$ as a diagonal coefficient of $\rho$,  as  $\xi$ 
 may not belong to $X_e$ in general.
  To circumvent this problem, 
it might be that some strictness assumption on $T$ should be introduced, as well as some notion of multiplier for Hilbert $\B$-bundles, analogous to the approach developed in \cite{BKMS}. Alternatively, Lin's version of Paschke's result for non-unital $C^*$-algebras (cf.~\cite[Theorem 2.1]{LinH}) might be helpful.
\end{remark} 

\subsection{Proof of Theorem \ref{GR} in the general case}\label{GR-II}
Consider now two unital Fell bundles $\A=(A_g)_{g \in G}$ and $\B=(B_h)_{h \in H}$,  $\varphi \in {\rm Hom}(G, H)$, and 
   a positive definite $\A$-$\varphi$-$\B$ bundle map $T=(T_g)_{g\in G}$. Although it is possible to generalize our arguments in the previous subsection to prove Theorem \ref{GR} also in this case, we briefly sketch below a more direct way of achieving this, in line with proofs of earlier results with a similar flavor. 

 For $r\in H$, set $X^0_r:=\big\{ \xi \in \prod_{k\in G} A_k \odot B_{\varphi(k)^{-1}r}: \xi \text{ is finitely supported}\big\}$. Then set
$\X^0 = \coprod_{r\in G} X_r^0$ (disjoint union).
We  define maps
\[\text{$\X^0\times \B \to \X^0, (\xi,b)\mapsto \xi\cdot b$, and $\X^0\times \X^0 \to \B, (\xi, \eta)\mapsto \langle \xi, \eta\rangle_{\B}$}\] as follows.
 Let $r, h\in H$, $\xi \in X_0^r$ and $b\in B_h$. Then, for each $k \in G$, set
\[ (\xi\cdot b)(k) := \xi(k)\cdot b \, \in A_k \odot B_{\varphi(k)^{-1}rh},\]
where  $y \mapsto y\cdot b$ denotes the unique linear map from $A_k \odot B_{\varphi(k)^{-1}r}$ into $A_k \odot B_{\varphi(k)^{-1}rh}$ satisfying
\[ \Big(\sum_{i=1}^n a_i \odot b_i\Big)\cdot b = \sum_{i=1}^n a_i \odot (b_ib)\]
whenever $a_1, \ldots, a_n \in A_k$ and $b_1, \ldots, b_n \in B_{\varphi(k)^{-1}r}$. We then get that $\xi\cdot b \in X^0_{rh}$ and the map $(\xi,b)\mapsto \xi\cdot b$ is clearly bilinear. 

Next, let $k, k' \in G$ and $r, r'\in H$. For $a\in A_k$, $a'\in A_{k'}, b \in B_{\varphi(k)^{-1}r}$ and $b' \in B_{\varphi(k')^{-1}r'}$, set
\[ [a\odot b, a'\odot b']_T := b^*T_{k^{-1}k'}(a^*a')b' \, \,\in \,B_{r^{-1}\varphi(k)}B_{\varphi(k^{-1}k')}B_{\varphi(k')^{-1}r'} \, \subseteq \,B_{r^{-1}r'}.\] 
Using the universal properties of algebraic tensor products, 
we can extend $[\cdot, \cdot]_T$ to \\ $(A_k \odot B_{\varphi(k)^{-1}r}) \times  (A_{k'} \odot B_{\varphi(k')^{-1}r'})$. Now, for $\xi \in X^0_r, \eta \in X^0_{r'}$, we set 
\[ \langle \xi, \eta\rangle_\B^0 := \sum_{k, k'\in G} [\xi(k), \eta(k')]_T \, \in B_{r^{-1}r'} .\]
 Using that $T$ is positive definite, it is a routine matter to check that $\X^0$ is a semi-inner product $\B$-bundle.
Hence 
we can let $\X=(X_r)_{r\in G}$ be the Hilbert $\B$-bundle obtained from $\X^0$ after separation and completion.
 
Next we will construct an $(\A, \varphi)$-action $\rho$ on $\X$. Let $g \in G$, $a\in A_g$ and $r\in H$. For any $k\in G, l \in H$, we let  $y \mapsto a\cdot y$ denote the unique linear map from $A_k \odot B_l$ into $A_{gk} \odot B_l$ satisfying
\[ a\cdot \Big(\sum_{i=1}^n a_i \odot b_i\Big) = \sum_{i=1}^n (aa_i) \odot b_i\]
whenever $a_1, \ldots, a_n \in A_k$ and $b_1, \ldots, b_n \in B_l$.
Then we define $\rho_0(a): X_r^0\to X^0_{\varphi(g)r}$ by
\[ (\rho_0(a)\xi)(k)= a\cdot \xi(g^{-1}k) \quad \text{for all } \xi \in X_r^0 \text{ and } k \in G.\]
This is meaningful since $\xi(g^{-1}k) \in A_{g^{-1}k}\odot B_{\varphi(g^{-1}k)^{-1}r}$, so $a\cdot \xi(g^{-1}k) \in A_{k}\odot B_{k^{-1}\varphi(g)r}$.
We note that 
\[ (\rho_0(a^*)\xi)(k)= a^*\cdot \xi((g^{-1})^{-1}k) = a^*\cdot \xi(gk) \quad \text{for all } \xi \in X_r^0 \text{ and } k \in G.\]  

\smallskip We have now defined a map $\rho_0(a): \X^0 \to \X^0$ for each $a\in \A$, and it can be verified without too much trouble that $\rho_0: \A \to \F(\X^0) $ is an $(\A, \varphi)$-pre-action on $\X^0$. For example, let us check that  \[\langle \rho_0(a)\xi, \eta\rangle_\B^0 = \langle \xi, \rho_0(a^*)\eta\rangle_\B^0\]
for all $a\in \A$ and $\xi, \eta \in \X^0$. We first note that if $k, k' \in G$ and $r, r'\in H$, then one readily sees that
\[[a\cdot x , y]_T = [x, a^*\cdot y]_T\]
for 
all $x \in A_k \odot B_{\varphi(k)^{-1}r}$ and $y\in A_{k'} \odot B_{\varphi(k')^{-1}r'}$. Hence we get
\[  \langle \rho_0(a)\xi, \eta\rangle_\B^0 =  \sum_{k, k'\in G} [a\cdot\xi(g^{-1}k), \eta(k')]_T 
 =  \sum_{k, k'\in G} [\xi(g^{-1}k), a^*\cdot \eta(k')]_T\]
 \[= \sum_{l, l'\in G} [\xi(l), a^*\cdot \eta(gl')]_T= \langle \xi, \rho_0(a^*)\eta\rangle_\B^0.\]
 After separation and completion, $\rho_0$ induces an $(\A, \varphi)$-action $\rho$ on $\X$, cf.~Remark \ref{pre-act-def}.
 
\medskip  Finally, for $r \in H, k' \in G, a'\in A_{k'}$ and $b' \in B_{\varphi(k')^{-1}r}$, let $\delta_{(k', a', b')}^r \in X_r^0$ be given by 
  \[ \delta_{(k', a', b')}^r(k):= \begin{cases} a'\odot b' \quad \text{if } k = k',\\ \ \ 0 \quad \text{otherwise}\end{cases}.\]
  Then for any $g \in G$ and $a\in A_g$, we have that \[\rho_0(a)\delta_{(e, 1_{A_e}, 1_{B_e})}^e = \delta_{(g, a, 1_{B_e})}^{\varphi(g)}.\]
  Thus, letting $\xi$ denote the canonical image of $\delta_{(e, 1_{A_e}, 1_{B_e})}^e$ in $X_e$, 
   we get that 
  \begin{align*} \langle \xi, \rho(a)\xi\rangle_\B &= \big\langle \delta_{(e, 1_{A_e}, 1_{B_e})}^e, \,\rho_0(a)\delta_{(e, 1_{A_e}, 1_{B_e})}^{e}\big\rangle_\B^0 = \big\langle \delta_{(e, 1_{A_e}, 1_{B_e})}^e,\, \delta_{(g, a, 1_{B_e})}^{\varphi(g)}\big\rangle_\B^0 \\ 
  &= \sum_{k, k'\in G}[  \delta_{(e, 1_{A_e}, 1_{B_e})}^e(k),\, \delta_{(g, a, 1_{B_e})}^{\varphi(g)}(k')]_T =
  [1_{A_e}\odot 1_{B_e}, a\odot 1_{B_e}]_T\\
  &= 1_{B_e}^* T_{e^{-1}g}(1_{A_e}^*a)1_{B_e}
  = T_g(a),
  \end{align*}
  as desired.

\subsection{Summing up}
By combining Theorem \ref{GR} with some previous results, we get the following:

\begin{theorem} \label{pd-etc}
Let $\A=(A_g)_{g \in G}$ and $\B=(B_h)_{h \in H}$ be Fell bundles and  $\varphi \in {\rm Hom}(G,H)$. Let $T=(T_g)_{g\in G}$ be an $\A$-$\varphi$-$\B$ bundle map.
Consider the following properties:
\begin{itemize}
\item[i$)$] $T$ is positive definite.
\item[ii$)$] $T$ is full and  $\Phi_T: C^*(\A) \to C^*(\B)$ is completely positive.
\item[iii$)$] $T$ is reduced  and $M_T: C^*_r(\A) \to C^*_r(\B)$ is completely positive.
\item[iv$)$] There exists an $(\A, \varphi)$-action $\rho$ on some  Hilbert $\B$-bundle $\X=(X_h)_{h \in H}$ and $x \in X_e$  such that 
\[T_g(a)= \langle x, \rho(a) x\rangle_\B\] for all $g\in G$ and $a\in A_g$.
\end{itemize}
Then we have i$)$ $\Leftrightarrow$ ii$)$, iii$)$ $\Rightarrow$ i$)$ and iv$)$ $\Rightarrow$ i$)$.
Further, if $\ker(\varphi)$ is amenable, then i$)$ 
$\Leftrightarrow$ iii$)$.
Also, if both $\A$ and $\B$ are unital, then i$)$ $\Rightarrow$ iv$)$.\\
In particular, if $\A$ and $\B$ are unital and the kernel of $\varphi \in {\rm Hom}(G,H)$ is amenable, then the conditions i$)$-iv$)$ are all equivalent.
\end{theorem}
\begin{proof}
The equivalence $i) \Leftrightarrow ii)$ is from Theorem \ref{PDCP}. The implication $iii) \Rightarrow i)$ follows from \cite[Proposition 3.9]{BeCo7} (see the comment after Definition 3.10). The implication  $iv) \Rightarrow i)$ is shown in Proposition \ref{action-pd}. 
 If $\ker(\varphi)$ is amenable, then $i) 
\Leftrightarrow iii)$ is stated in Theorem \ref{PDCP}.
Finally, if both $\A$ and $\B$ are unital, then $i) \Rightarrow iv)$ by Theorem \ref{GR}.
\end{proof}

\noindent Of course, if $T$ is an $\A$-$\B$ bundle map (i.e., $G=H$ and $\varphi={\rm id}_G$), then the amenability assumption on the kernel is trivially satisfied.

\smallskip \noindent 
For completeness, we also mention that one can
add to the last statement of Theorem \ref{pd-etc} another equivalent characterization of positive definiteness of $T$ 
that can be expressed in terms of diagonal coefficient maps associated to 
 some $C^*(\A)$-$C^*(\B)$ correspondence (or some $C_r^*(\A)$-$C_r^*(\B)$ correspondence if $\ker(\varphi)$ is amenable), see the discussion in \cite[Example 4.8]{BeCo7}.
 
\section{$C^*$-correspondences associated to actions of Fell bundles}\label{C*-c}

Let $\A=(A_g)_{g\in G}$ and $\B=(B_h)_{h\in H}$ be Fell bundles over $G$ and $H$, and $\varphi \in {\rm Hom}(G, H)$. 
Several papers in the literature deal with the construction of $C^*$-correspondences over $C^*$-crossed products associated to $C^*$-dynamical systems (see e.g.~\cite{Combes, Kasp, Kasp95, EKQR00, EKQR06, BeCo6}), in order to provide functors from the natural equivariant category of  $C^*$-dynamical systems to the category of $C^*$-algebras with morphisms given by isomorphism classes of $C^*$-correspondences. This functorial approach is a fundamental tool in the study of the representation theory and the structure of $C^*$-crossed products, in particular in connection with induction processes and imprimitivity results.  Our aim in this section is to construct a natural $C^*(\A)$-$C^*(\B)$ correspondence $Y=\X\rtimes\B$
associated to an $(\A, \varphi)$-action on a  Hilbert $\B$-bundle $\X$.
We also discuss $C_r^*(\A)$-$C_r^*(\B)$ correspondences. 
It is quite clear that this construction will give rise to a functor from the category of Fell bundles over discrete groups with morphisms consisting of actions to the category of $C^*$-algebras mentioned above, but we will not discuss this here any further.

\medskip We will say that an $\A$-$\varphi$-$\B$ action $\rho$ on some Hilbert $\B$-bundle $\X$  is \emph{nondegenerate}
when for every $k\in H$,

\medskip  Span$\{(\rho(a)w)b: g \in G, a \in A_g, w\in X_e, h\in H, b\in B_h,  \varphi(g)h = k\}$ is dense in $X_{k}$.

\medskip 

\medskip\noindent  It is then obvious that $\rho$ is nondegenerate whenever $\rho$ has a cyclic vector $x\in X_e$ (in the sense of Remark \ref{cyclic1}). 
The main point is that this notion will imply that the associated left action of $C^*(\A)$ on the Hilbert $C^*(\B)$-module  $Y=\X\rtimes\B$ is nondegenerate, i.e., the span of $\{z\cdot y: z \in C^*(\A), y \in Y\}$ is dense in $Y$. It will also imply that the associated left action of $C_r^*(\A)$ on the Hilbert $C_r^*(\B)$-module  $\X\rtimes_r\B$ is nondegenerate whenever it is well-defined. 

\subsection{The main construction} \label{mainconstr}
As a start-up, we  describe how to associate a Hilbert $C^*(\B)$-module $Y$ (resp.~a Hilbert $C_r^*(\B)$-module $Y_r$) to a  given Hilbert $\B$-bundle $\X=(X_h)_{h\in H}$. Let $C_c(\X)$ denote the vector space of finitely supported sections $\xi$  from $H$ to $\X$, i.e., $\xi:H\to \X$ is such that $\xi(h) \in X_h$ for every $h\in H$ and $\xi(h)= 0$ for all but finitely many $h$'s in $H$.
 For  $f\in C_c(\B)$ and $\xi, \eta \in C_c(\X)$, we define 
$\xi\cdot f \in C_c(\X)$ and  $\langle \xi, \eta\rangle_{C_c(\B)} \in C_c(\B)$ by
\begin{align*} 
(\xi\cdot f)(h)  &:= \sum_{k\in H} \xi(k) \cdot f(k^{-1}h),\\
\langle \xi, \eta\rangle_{C_c(\B)}(h) &:= \sum_{k\in H} \, \langle \xi(k), \eta(kh) \rangle_\B
\end{align*}
for all $h\in H$. Then it is not difficult to verify that these operations are well-defined and that the first one turns $C_c(\X)$ into
 a right $C_c(\B)$-module. To simplify notation, we will from now on suppress $\kappa^\B$ (resp.~$\iota^\B$) in our notation, i.e., we will identify $C_c(\B)$  with its canonical copy in $C^*(\B)$ (resp.~in $C_r^*(\B)$). Moreover, we will denote the norm of $f\in C_c(\B)$ in $C^*(\B)$ (resp.~in $C_r^*(\B)$) by $\|f\|_{C^*(\B)}$ (resp.~$\|f\|_{C_r^*(\B)}$), so that $ \|f\|_{C^*(\B)} = \|f\|_{\rm u}$ (resp.~$ \|f\|_{C_r^*(\B)} = \|f\|_{r}$).  
  \begin{proposition} \label{indprod}  $\langle \cdot, \cdot \rangle_{C_c(\B)} $
 is a  $C^*(\B)$-valued inner product on $C_c(\X)$.
  \end{proposition}   
  
  \vspace{-5ex}
  \begin{proof} 
  We check only two of the axioms, leaving to the reader to verify the other ones.  
   Let $\xi, \eta \in C_c(\X), f \in C_c(\B)$ and $h\in H$. Then we have
  \begin{align*} \langle \xi, \eta\cdot f\rangle_{C_c(\B)}(h) &= \sum_{h'\in H} \, \langle \xi(h'), (\eta\cdot f)(h'h) \rangle_\B =  \sum_{h',k\in H} \, \langle \xi(h'), \eta(k) \cdot f(k^{-1}h'h) \rangle_\B\\
   &=  \sum_{h',k\in H} \, \langle \xi(h'), \eta(k) \rangle_\B \, f(k^{-1}h'h)
  =\sum_{h', k\in H}  \langle \xi(h'), \eta(h'k)\rangle_{\B} \,f(k^{-1}h) \\
  &=\sum_{k\in H}  \langle \xi, \eta\rangle_{C_c(\B)}(k) \,f(k^{-1}h)=(\langle \xi, \eta\rangle_{C_c(\B)} \star f) (h),
  \end{align*}
  showing that $\langle \xi, \eta\cdot f\rangle_{C_c(\B)} = \langle \xi, \eta\rangle_{C_c(\B)} \star f$.
  
 Next, let $\pi=(\pi_h)_{h\in H}$ be a $*$-representation of $\B$ on some Hilbert space $\H$ such that 
 its integrated form $\Phi_\pi: C^*(\B)\to \L_{\Complessi}(\H)$ is  faithful. 
  We will show that $\Phi_\pi (\langle \xi, \xi\rangle_{C_c(\B)}) $ is positive in $\L_{\Complessi}(\H)$ for every $\xi \neq 0 \in C_c(\X)$. Let $v \in \H$. Then 
\begin{align*} \langle v, \Phi_\pi (\langle \xi, \xi\rangle_{C_c(\B)})v\rangle_{\H}&= \sum_{h\in H} \big\langle v,  \pi_h\big(\langle \xi, \xi\rangle_{C_c(\B)}(h)\big) v\big\rangle_{\H}\\
&= \sum_{k, h\in H} \big\langle v,  \pi_h\big(\langle \xi(k), \xi(kh)\rangle_\B\big) v\big\rangle_{\H}\\
&=  \sum_{k\in H} \sum_{h\in H} \big\langle v,  \pi_{k^{-1}h}\big(\langle \xi(k), \xi(h)\rangle_\B\big) v\big\rangle_{\H}
\end{align*} 
Let now $h_1, \ldots, h_n$ be an enumeration of the support of $\xi$, and set $\mathbf{h} := (h_1, \ldots, h_n) \in H^n$. 
As $\xi(h_i) \in X_{h_i}$ for each $i$, we get from Lemma \ref{pos-mat} that 
\[\big[\langle \xi(h_i), \xi(h_j)\rangle_\B\big] \in M_{\mathbf{h}}(\B)^+.\]
Since the map $[r_{ij}] \mapsto [\pi_{{h_i}^{-1}h_j}(r_{ij})]$ is  a $*$-homomorphism 
from $M_{\mathbf{h}}(\B)$ into 
$M_n(\L_{\Complessi}(\H))$, this implies that
\[ \big[\pi_{{h_i}^{-1}h_j}\big(\langle \xi(h_i), \xi(h_j)\rangle_\B\big)\big] \in M_n(\L_{\Complessi}(\H))^+,\]
and it readily follows that 
 \begin{equation*} \langle v, \Phi_\pi (\langle \xi, \xi\rangle_{C_c(\B)})v\rangle_{\H} 
 = \sum_{i, j=1}^n \big\langle v,  \pi_{h_{i}^{-1}h_j}\big(\langle \xi(h_i), \xi(h_j)\rangle_\B\big) v\big\rangle_{\H} \geq 0.
 \end{equation*}
 Thus,  $\Phi_\pi (\langle \xi, \xi\rangle_{C_c(\B)}) \in \L_{\Complessi}(\H) \geq 0$, i.e., $\langle \xi, \xi\rangle_{C_c(\B)}$ is positive in $C^*(\B)$ (since $\Phi_\pi$ is faithful).

To show definiteness of $\langle \cdot, \cdot\rangle_{C_c(\B)}$,
 assume that $\langle \xi, \xi\rangle_{C_c(\B)} = 0$. Then \[\sum_{k\in H} \langle \xi(k), \xi(k)\rangle_\B = \langle \xi, \xi\rangle_{C_c(\B)}(e) = 0.\] As $ \langle \xi(k), \xi(k)\rangle_\B \in B_e^+$ for each $k\in H$, we get that $\langle \xi(k), \xi(k)\rangle_\B= 0$, i.e., $\xi(k) = 0$ for every $k\in H$. Thus $\xi =0$, as desired. 
\end{proof}
  Setting \[\|\xi\| := \big\| \langle \xi, \xi\rangle_{C_c(\B)}\big\|_{C^*(\B)}^{1/2}\]
 for every $\xi \in C_c(\X)$, we get a norm on $C_c(\X)$. Taking the completion 
 of the right inner product $C_c(\B)$-module $C_c(\X)$ w.r.t.~this norm,
  we get a  Hilbert $C^*(\B)$-module $Y$ associated to the Hilbert $\B$-bundle $\X$, that we will sometimes denote by $\X \rtimes \B$.
 
Note that by a similar argument as in the proof above, $\langle \cdot, \cdot \rangle_{C_c(\B)} $ may also be considered 
as a $C_r^*(\B)$-valued inner product on $C_c(\B)$.
Hence, setting \[\|\xi\|_r := \big\| \langle \xi, \xi\rangle_{C_c(\B)}\big\|_{C_r^*(\B)}^{1/2}\]
 for every $\xi \in C_c(\X)$, we get another norm on $C_c(\X)$, and taking the completion w.r.t.~$\|\cdot\|_r$,  
  we get a  Hilbert $C_r^*(\B)$-module $Y_r$ associated to the Hilbert $\B$-bundle $\X$,  sometimes denoted by $\X \rtimes_{r} \B$.

 \medskip 
 Next, assume that $\rho$ is an $(\A, \varphi)$-action on the Hilbert $\B$-bundle $\X=(X_h)_{h\in H}$. For  $f\in C_c(\A)$ and $\xi \in C_c(\X)$, we define 
$f\cdot \xi \in C_c(\X)$ by 
\begin{equation} \label{left-action} (f\cdot \xi)(h) := \sum_{g\in G} \rho(f(g)) \xi(\varphi(g)^{-1}h)\end{equation}
for all $h\in H$. 
One easily checks that this turns $C_c(\X)$ into a left $C_c(\A)$-module.  Further, using property (iv) for $\rho$ in Definition \ref{BC-action-def}, it readily follows that
\begin{equation}
f \cdot (\xi \cdot f') = (f \cdot \xi) \cdot f'
\end{equation}
for all $f \in C_c(\A)$, $f'\in C_c(\B)$ and $\xi \in C_c(\X)$, i.e., the left action of $C_c(\A)$ and the right action of $C_c(\B)$ on $C_c(\X)$ commute. 

The left action of $C_c(\A)$ on $C_c(\X)$ is associated to a $*$-representation $\pi^0$ of $\A$ in  $C_c(\X)$  defined as follows. 
For $g\in G$ and $a\in A_g$, let $\pi^0_g(a): C_c(\X) \to C_c(\X)$ be the  linear map given by
\[(\pi^0_g(a)\xi)(h) := \rho(a) \xi(\varphi(g)^{-1}h)\] for all $\xi \in C_c(\X)$ and $h\in H$. In particular, if $a \in A_e$, then
$(\pi_e^0(a) \xi)(h) = \rho(a) \xi(h)$. 
 \begin{lemma}  $\pi^0 := (\pi^0_g)_{g\in G}$ is a $*$-representation of $\A$ in $C_c(\X)$. 
 \end{lemma} 
 \begin{proof} Let $g,g' \in G, k\in H$, $a \in A_g, a'\in A_{g'}$ and $\xi, \eta \in C_c(\X)$. Then we have
\begin{align*} (\pi_g^0(a) \pi_{g'}^0(a') \xi)(k) &= \rho(a) (\pi_{g'}^0(a') \xi)(\varphi(g)^{-1}k)= \rho(a) \rho(a') \xi(\varphi(g')^{-1}\varphi(g)^{-1}k) \\ &= \rho(aa') \xi(\varphi(gg')^{-1}k) =(\pi_{gg'}^0(aa') \xi)(k),\end{align*}
hence  $\pi_g^0(a) \pi_{g'}^0(a') = \pi_{gg'}^0(aa')$. Moreover, for all $h\in H$,
\begin{align*} \big\langle \xi, \pi^0_g(a)\eta \big\rangle_{C_c(\mathcal{\B})}(h) &= \sum_{k\in H} \, \big\langle \xi(k), (\pi^0_g(a)\eta)(kh) \big\rangle_\B = \sum_{k\in H} \, \big\langle \xi(k), \rho(a)\eta(\varphi(g)^{-1}kh) \big\rangle_\B\\ 
& =\sum_{k\in H} \, \big\langle \rho(a^*)\xi(k), \eta(\varphi(g)^{-1}kh) \big\rangle_\B  =\sum_{k\in H} \, \big\langle \rho(a^*)\xi(\varphi(g)k), \eta(kh) \big\rangle_\B \\ &= \sum_{k\in H} \, \big\langle (\pi^0_{g^{-1}}(a^*)\xi)(k), \eta(kh) \big\rangle_\B 
=\big\langle \pi^0_{g^{-1}}(a^*) \xi, \eta \big\rangle_{C_c(\mathcal{\B})}(h),
\end{align*}
i.e., $ \big\langle \xi, \pi^0_g(a)\eta \big\rangle_{C_c(\mathcal{\B})} = \big\langle \pi^0_{g^{-1}}(a^*) \xi, \eta \big\rangle_{C_c(\mathcal{\B})}$. Thus, $\pi^0_g(a)^*=\pi^0_{g^{-1}}(a^*)$. 
\end{proof}
In order to extend $\pi^0$ to a $*$-representation of $\A$ on $Y$ (resp.~$Y_r$), the following lemma will be crucial. 
\begin{lemma} \label{pi_e}
Let $a\in A_e$ and $\xi\in C_c(\X)$. Then 
\[ \|\pi_e^0(a)\xi\| \leq \|a\|\, \|\xi\| \, \text{ and } \, \|\pi_e^0(a)\xi\|_r \leq \|a\| \, \|\xi\|_r.\] 
\end{lemma}

\begin{proof} Let $\phi=(\phi_g)_{g\in G}$ be a $*$-representation of $\B$ 
in $\L_D(\H)$, 
where  $D$ is some $C^*$-algebra and $\H$ is a Hilbert $D$-module.
 Let then $\Phi: C^*(\B) \to  \L_D(\H)$ denote the integrated form of $\phi$.
  We may assume $\xi \neq 0 \in C_c(\X)$. Let $h_1, \ldots, h_n$ be an enumeration of the support of $\xi$ and set $\mathbf{h} := (h_1, \ldots, h_n) \in H^n$.  Then, for each $v$ in the unit ball of $\H$, we have
 \begin{align*}
 \langle v , \Phi\big(\langle \pi_e^0(a)\xi, \pi_e^0(a)\xi\rangle_{C_c(\B)}\big) v\rangle_\H&
 = \langle v , \sum_{h\in H} \phi_h\big(\langle \pi_e^0(a)\xi, \pi_e^0(a)\xi\rangle_{C_c(\B)}(h)\big) v\rangle_{\H}\\
 &= \langle v , \sum_{k, h\in H} \phi_h\big(\langle (\pi_e^0(a)\xi)(k), (\pi_e^0(a)\xi)(kh)\rangle_{\B}\big) v\rangle_{\H}\\
 &= \langle v , \sum_{k, h\in H} \phi_h\big(\langle \rho(a)\xi(k), \rho(a)\xi(kh)\rangle_{\B}\big) v\rangle_{\H}\\
  &= \langle v , \sum_{k, h\in H} \phi_{k^{-1}h}\big(\langle \rho(a)\xi(k), \rho(a)\xi(h)\rangle_{\B}\big) v\rangle_{\H}\\
   &= \langle v , \sum_{i, j=1}^n \phi_{h_i^{-1}h_j}\big(\langle \rho(a)\xi(h_i), \rho(a)\xi(h_j)\rangle_{\B}\big) v\rangle_{\H}\\
     &= \langle v , \sum_{i, j=1}^n \phi_{h_i^{-1}h_j}\big(\langle y_i, y_j\rangle_{\B}\big) v\rangle_{\H}
 \end{align*}
where $y_i:= \rho(a)\xi(h_i) \in X_{h_i}$ for each $i=1, \ldots, n$. Since the matrix $\big[\langle y_i, y_j\rangle_{\B}\big]$ is a positive element in $M_\mathbf{h}(\B)$, cf.~Lemma \ref{pos-mat}, there exists some $C = [c_{ij}] \in M_\mathbf{h}(\B)$ such that $C^*C = \big[\langle y_i, y_j\rangle_{\B}\big]$, i.e., we have
\[  \langle y_i, y_j\rangle_{\B} = \sum_{k=1}^n c_{ki}^*\,c_{kj}\,\]
for all $i,j=1, \ldots, n$.
Thus we get
\begin{align*}
\langle v , \sum_{i, j=1}^n \phi_{h_i^{-1}h_j}\big(\langle y_i, y_j\rangle_{\B}\big) v\rangle_{\H}&=\langle v , \sum_{i, j, k=1}^n \phi_{h_i^{-1}h_k h_k^{-1}h_j}\big( c_{ki}^*\,c_{kj}\big) v\rangle_{\H}\\
&=\langle v , \sum_{i, j, k=1}^n \phi_{h_i^{-1}h_k} \big(c_{ki}^*\big) \,\phi_{h_k^{-1}h_j}\big(c_{kj}\big) v\rangle_{\H}\\
&=\sum_{k=1}^n\Big\langle \sum_{i=1}^n\phi_{h_k^{-1}h_i} \big(c_{ki}\big)v , \sum_{j=1}^n \phi_{h_k^{-1}h_j}\big(c_{kj}\big) v\Big\rangle_{\H}\\
&= \sum_{k=1}^n\big\langle v_k, v_k\big\rangle_{\H} =  \sum_{k=1}^n\big \|v_k\|^2,
\end{align*}
where $v_k:= \sum_{i=1}^n\phi_{h_k^{-1}h_i} \big(c_{ki}\big)v \in \H$ for $k=1, \ldots, n$. 

Now, for each $i, j =1, \ldots, n$, set
$t_{ij} := \phi_{h_i^{-1}h_j}(c_{ij}) \in \L_D(\H)$ and $T:= [t_{ij}] \in M_n(\L_D(\H))$. Letting $M_n(\L_D(\H))$ act on the direct sum $\H^n$ in the standard way, we get that 
\[ \sum_{k=1}^n\big \|v_k\|^2 = \|T(v, \ldots, v)\|^2 = \Big\langle (v, \ldots, v), T^*T(v, \ldots, v) \Big\rangle_{\H^n}.\]
On the other hand, we have that 
\begin{align*} T^*T &= \Big[\sum_{k=1}^n t_{ki}^*\,t_{kj}\Big] = \Big[ \sum_{k=1}^n\phi_{h_i^{-1}h_k}(c_{ki}^*) \phi_{h_k^{-1}h_j}(c_{kj})\Big]\\ &= \Big[ \sum_{k=1}^n\phi_{h_i^{-1}h_j}(c_{ki}^*c_{kj})\Big] =  \Big[ \phi_{h_i^{-1}h_j}\big(\langle y_i, y_j\rangle_\B\big]\big)\Big] \\
&= \phi_\mathbf{h}\Big(\big[\langle y_i, y_j\rangle_\B\big]\Big).
\end{align*}
Hence, altogether,  we get that 
\[  \big\langle v , \Phi\big(\langle \pi_e^0(a)\xi, \pi_e^0(a)\xi\rangle_{C_c(\B)}\big) v\big\rangle_{\H} = \sum_{k=1}^n\big \|v_k\|^2 =  \Big\langle (v, \ldots, v), \phi_\mathbf{h}\Big(\big[\langle y_i, y_j\rangle_\B\big]\Big)(v, \ldots, v) \Big\rangle_{\H^n}.\]
Setting $x_i:= \xi(h_i)$ for $i=1, \ldots, n$, a computation very similar to the one above gives that 
\[\langle v , \Phi\big(\langle \xi, \xi\rangle_{C_c(\B)}\big) v\rangle_{\H} = 
   \Big\langle (v, \ldots, v), \phi_\mathbf{h}\Big(\big[\langle x_i, x_j\rangle_\B\big]\Big) (v, \ldots, v) \Big\rangle_{\H^n}.\]
   Now, by Lemma \ref{Lance-gen}, we have $ \big[\langle y_i, y_j\rangle_\B\big] \leq \|a\|^2\, \big[\langle x_i, x_j\rangle_\B\big]$. Thus we get that
      \[  \big\langle v , \Phi\big(\langle \pi_e^0(a)\xi, \pi_e^0(a)\xi\rangle_{C_c(\B)}\big) v\big\rangle_{\H} \leq \|a\|^2 \langle v , \Phi\big(\langle \xi, \xi\rangle_{C_c(\B)}\big) v\rangle_{\H}.\]
      Since this holds for all $v$ in the unit ball of $\H$, we conclude that 
    \[  \|  \Phi\big(\langle \pi_e^0(a)\xi, \pi_e^0(a)\xi\rangle_{C_c(\B)}\big)\| \leq \|a\|^2 \, \| \Phi\big(\langle \xi, \xi\rangle_{C_c(\B)}\big)\|.\]
\quad Choosing $\phi$ such that $\Phi$ is a faithful $*$-representation of $C^*(\B)$ on a Hilbert space $\H$,
this gives that 
     \[ \|\pi_e^0(a)\xi \|^2 = \| \langle \pi_e^0(a)\xi, \pi_e^0(a)\xi\rangle_{C_c(\B)}\|_{C^*(\B)}\, \leq \|a\|^2 \, \| \langle \xi, \xi\rangle_{C_c(\B)}\|_{C^*(\B)} =  \|a\|^2 \, \|\xi \|^2,\]
    which shows the first inequality stated in the lemma.
    
    Next, choosing $\phi = \lambda^\B$, so that $\Phi = \Lambda^\B: C^*(\B) \to C_r^*(\B) \subseteq \L_{B_e}(\ell^2(\B))$, we get that 
     \[ \|\pi_e^0(a)\xi \|_r^2 = \| \langle \pi_e^0(a)\xi, \pi_e^0(b)\xi\rangle_{C_c(\B)}\|_{C_r^*(\B)}\, \leq \|a\|^2 \, \| \langle \xi, \xi\rangle_{C_c(\B)}\|_{C_r^*(\B)} =  \|a\|^2 \, \|\xi \|_r^2,\]
     which shows the second inequality.
\end{proof}  
     
   Using Lemma \ref{pi_e}, we can extend $\pi_e^0(a)$ for each $a\in A_e$ to a bounded operator $\pi_e(a) $ on $Y$ (resp.~$\pi^r_e(a) $ on $Y_r$), satisfying 
    that  $\pi_e(a) \in \L_{C^*(\B)}(Y)$, with $\pi_e(a)^* = \pi_e(a^*)$ (resp.~$\pi^r_e(a) \in \L_{C_r^*(\B)}(Y_r)$, with $\pi^r_e(a)^* = \pi^r_e(a^*)$). 
   It is then clear that $\pi_e$ (resp.~$\pi^r_e$) is a $*$-homomorphism from $A_e$ into $\L_{C^*(\B)}(Y)$ (resp.~$\L_{C_r^*(\B)}(Y_r)$).

\begin{lemma} \label{b-bounded} Let $g \in G$, $a \in A_g$ and $\xi\in C_c(\X)$. Then we have 
\[ \|\pi_g^0(a)\xi\| \leq \|a\|\, \|\xi\| \, \text{ and } \, \|\pi_g^0(a)\xi\|_r \leq \|a\| \, \|\xi\|_r.\]  
\end{lemma}
\begin{proof}  We first show that $\big\langle \pi^0_{e}(a^*a)\xi, \xi \big\rangle_{C_c(\mathcal{\B})} \leq \|a^*a\|\, \big\langle \xi, \xi \big\rangle_{C_c(\mathcal{\B})}$. Set $d:= a^*a \in A_e^+$.
Then, using \cite[Proposition 1.2]{La1} at the second step, 
we get that
\begin{align*} \big\langle \pi^0_{e}(a^*a)\xi, \xi \big\rangle_{C_c(\B)} &= \big\langle \pi_{e}(d^{1/2})\xi, \pi_{e}(d^{1/2})\xi \big\rangle_{C^*(\B)} \leq \|\pi_e(d^{1/2})\|^2 \, \big\langle \xi, \xi \big\rangle_{C^*(\B)}\\
&= \|\pi_e(d)\| \, \big\langle \xi, \xi \big\rangle_{C_c(\B)} \leq \|a^*a\|\, \big\langle \xi, \xi \big\rangle_{C_c(\B)}.
\end{align*}
Thus, we obtain that
\begin{align*} \|\pi^0_g(a) \xi\|^2 &= \big\|  \big\langle \pi^0_g(a)\xi, \pi^0_g(a)\xi \big\rangle_{C_c(\mathcal{\B})}\big\|_{C^*(\B)} = \big\|  \big\langle \pi^0_{g^{-1}}(a^*)\pi^0_g(a)\xi, \xi \big\rangle_{C_c(\B)}\big\|_{C^*(\B)}\\
& =\big\|  \big\langle \pi^0_{e}(a^*a)\xi, \xi \big\rangle_{C_c(\B)}\big\|_{C^*(\B)}
\leq  \|a^*a\| 
\, \big\|\big\langle \xi, \xi \big\rangle_{C_c(\B)}\big\|_{C^*(\B)} = \|a\|^2 
\, \big\|\xi\big\|^2.
\end{align*}
This shows the first desired inequality. 
The second can be shown in an analogous way by replacing $C^*(\B)$ with $C_r^*(\B)$.
\end{proof}
Using Lemma \ref{b-bounded}, we can extend  $\pi^0_g(a)$ for every $g\in G$ and $a\in A_g$ to a bounded operator $\pi_g(a)$ on $Y$ (resp. $\pi_g^r(a)$ on $Y_r$), which is adjointable with $\pi_g(a)^* = \pi_{g^{-1}}(a^*)$ (resp.~$\pi^r_g(a)^* = \pi^r_{g^{-1}}(a^*)$). It follows that   $\pi=\{\pi_g\}_{g\in G}$
(resp.~$\pi^r=\{\pi_g^r\}_{g\in G}$)
 is a $*$-representation of $\A$ in $\L_{C^*(\B)}(Y)$ (resp. $\L_{C_r^*(\B)}(Y_r)$).

\begin{theorem} \label{correspond} 
Let $\rho$ be an $(\A, \varphi)$-action on the Hilbert $\B$-bundle $\X=(X_h)_{h\in H}$.
Then the left action of $C_c(\A)$ on $C_c(\X)$ defined in (\ref{left-action}) extends to a 
left action $\Phi$ of $C^*(\A)$ on the Hilbert $C^*(\B)$-module 
 $Y= \X\rtimes \B$. Hence $Y$ is a $C^*(\A)$-$C^*(\B)$ correspondence. Moreover, $\Phi$ is nondegenerate whenever $\rho$ is nondegenerate. 
\end{theorem}
\begin{proof} Let $\Phi:=\Phi_\pi:C^*(\A)\to \L_{C^*(\B)}(Y)$ be the integrated form of $\pi$. Then we have
\begin{align*} (\Phi(f) \xi) (h) &= \Big(\sum_{g\in G} \pi_g(f(g))\xi\Big)(h) = \sum_{g\in G} (\pi^0_g(f(g))\xi)(h)\\
&= 
 \sum_{g\in G} \rho(f(g))\xi(\varphi(g)^{-1}h) = (f\cdot \xi)(h)
 \end{align*}
 for all  $f\in C_c(\A)$, $\xi \in C_c(\X)$ and $h\in H$. Hence,  $\Phi(f) \xi = f\cdot \xi$ for all $f\in C_c(\A)$ and $\xi \in C_c(\X)$. This shows that the left action of $C_c(\A)$ on $C_c(\X)$ extends to the left action $\Phi$ of $C^*(\A)$ on $Y$. 
 
 Assume now that $\rho$ is nondegenerate.  In order to show that $\Phi$ is nondegenerate, using that $C_c(\X)$ is dense in $Y$, it suffices to show that  \[C_c(\X) \subseteq \overline{{\rm Span}\big(\{f\cdot \xi: f \in C_c(\A), \xi \in C_c(\X)\}\big)}.\]  
Let $\eta \in C_c(\X)$. As  $\eta = \sum_{j=1}^n x_j \odot k_j$ for some $k_1, \ldots, k_n \in H$ and $x_j \in X_{k_j}$, $j=1, \ldots, n$,
it is clear that we only have to consider $\eta= v\odot k$ where $k \in H$ and $v\in X_k$. 

Let $\varepsilon >0$. Using nondegeneracy of $\rho$, we can 
find $v'\in X_k$ of the form $v'= \sum_{i=1}^m \rho(a_i)(w_ib_i)$, where $g_1, \ldots, g_m \in G$, $h_1, \ldots, h_m \in H$, $w_1, \ldots, w_m \in X_e$, and $a_i \in A_{g_i}, b_i \in B_{h_i}, \varphi(g_i)h_i = k$ for $i=1, \ldots m$,  such that $\|v-v'\| < \varepsilon$.

  Now, observe that for all $w \in X_k$ and $h\in H$, we have 
       \begin{equation}\label{eq-ind} \langle w\odot k, w \odot k\rangle_{C_c(\B)}(h)= \sum_{k'\in H} \langle (w\odot k)(k'), (w\odot k)(k'h)\rangle_\B = \langle w , (w\odot k)(kh)\rangle_\B= (\langle w, w\rangle_\B \odot e)(h).
       \end{equation}
       Hence,
       \[ \|w\odot k\|^2 = \|\langle w\odot k, w \odot k\rangle_{C_c(\B)}\|_{C^*(\B)} =
       \|\langle w, w\rangle_\B \odot e\|_{C^*(\B)}= \|\langle w, w\rangle_\B\| = \|w\|^2,\]
       where we have used that if $b\in B_e$, then $\| b \odot e\|_{C^*(\B)} = \| \widehat j^\B_e(b)\| = \|b\|$ (cf.~\cite[Proposition 17.9 (iv)]{Exel2} for the last equality). 
       
       Setting $\eta':= v'\odot k =\sum_{i=1}^m (\rho(a_i)(w_ib_i)) \odot \varphi(g_i)h_i$,  we get  that   
      \[\|\eta-\eta'\|= \|v\odot k - v'\odot k\| = \|(v-v')\odot k\| = \|v-v'\| < \varepsilon.\]  
Moreover, an easy computation gives that \[(\rho(a_i)(wb_i)) \odot \varphi(g_i)h_i = (a_i \odot g_i) \cdot (wb_i \odot h_i)\] for each $i=1, \ldots, m$, so  $\eta'  \in {\rm Span}\big(\{f\cdot \xi: f \in C_c(\A), \xi \in C_c(\X)\}\big)$.

\smallskip This shows that $\eta \in \overline{{\rm Span}\big(\{f\cdot \xi: f \in C_c(\A), \xi \in C_c(\X)\}\big)}$, as desired.
\end{proof}

\begin{remark} \label{full-compl} Let $x\in X_e$ and $T$ be the $\A$-$\varphi$-$\B$ bundle map given by
 \[ T_g(a) = \langle x, \rho(a)x\rangle_\B \quad \text{for all } g\in G, a\in A_g.\]
 Since $T$ is positive definite, we know from \cite[Theorem 3.10]{BeCo7} (cf.~Theorem \ref{pd-etc}, $ i) \Rightarrow ii)$) that 
 $T$ is full with  $\Phi_T$ completely positive. This can also be deduced from Theorem \ref{correspond} as follows.
 
Define $\xi_x \in C_c(\X)$ by 
 \[\xi_x(h) =\begin{cases} x  \quad \text{if } h=e,\\ 0 \quad \text{if } h\neq e\end{cases}.\] 
 Then we get a completely positive linear map $\Psi: C^*(\A)\to C^*(\B)$ by setting 
  \[\Psi(z) := \langle \xi_x, z\cdot \xi_x \rangle_{C^*(\B)} \quad \text{ for all } z \in C^*(\A).\]
 Note that $\Psi$ maps $C_c(\A)$ into $C_c(\B)$ 
  (since $f\cdot \xi_x \in C_c(\X)$ whenever $f \in C_c(\A)$).
 We claim that $\Psi = \phi_T$ on $C_c(\A)$. 
 
 Indeed, let $f\in C_c(\A)$, $h \in H$, and set $F:=\{g\in G : f(g)\neq 0, \varphi(g) = h\}$, which is a finite subset of $G$. Then  we have
 \[ (f\cdot\xi_x)(h) = \sum_{g\in G} \rho(f(g))\xi_x(\varphi(g)^{-1}h)= \sum_{g \in F} \rho(f(g)) x\,,\]
hence
 \begin{align*} (\Psi(f))(h) &
 =\langle \xi_x, f\cdot \xi_x\rangle_{C_c(\B)}(h) \\
 &= \sum_{k\in H} \langle \xi_x(k), (f\cdot\xi_x)(kh)\rangle_\B = \langle x, (f\cdot\xi_x)(h)\rangle_\B\\
& =  \sum_{g \in F} \langle x, \rho(f(g)) x\rangle_\B = \sum_{g \in F} T_g(f(g)) \\
&= \Big(\sum_{g\in G} j_{\varphi(g)}^\B(T_g(f(g)))\Big)(h) = (\phi_T(f))(h),
\end{align*}
which shows that
$\Psi(f) = \phi_T(f)$, as claimed.
This shows that $T$ is full with $ \Phi_T=\Psi $.
 \end{remark}

\begin{remark} \label{mixed-cor}
Let $\Phi_{\pi^r}:C^*(\A)\to \L_{C^*_r(\B)}(Y_r)$ denote the integrated form of the $*$-representation $\pi^r$ from  $\A$ into $\L_{C^*_r(\B)}(Y_r)$. Then, proceeding as in the proof of Theorem \ref{correspond} with $\Phi_{\pi^r}$ instead of $\Phi_{\pi}$, we can extend the left action of $C_c(\A)$ on $C_c(\X)$ defined in (\ref{left-action}) to a left action of $C^*(\A)$ on the Hilbert $C_r^*(\B)$-module 
 $Y_r= \X\rtimes_r \B$, thereby turning $Y_r$ into a $C^*(\A)$-$C_r^*(\B)$ correspondence.

It is not obvious which conditions will guarantee that $\Phi_{\pi^r}$ factors through a $*$-homomorphism $\Psi^{r}: C_r^*(\A)\to \L_{C^*_r(\B)}(Y_r)$ such that $\Phi_{\pi^r} = \Psi^{r}\circ \Lambda^\A$, 
 which would ensure that the left action of $C_c(\A)$ on $C_c(\X)$ can be extended and turn $Y_r$ into a $C_r^*(\A)$-$C_r^*(\B)$ correspondence.
 We believe that this should be possible at least when $G=H$ and $\varphi$ is the identity map. This is for example known to be true in the case where $\A$ and $\B$ are the Fell bundles over $G$ associated to some $C^*$-dynamical systems, 
 see Example \ref{C*-c-dyn-syst} for further details. 
We will discuss some related results in subsection \ref{C*-c-red}. 
\end{remark}

Despite the difficulties discussed in the remark above, we show right away that
the Fell's absorption principle for Fell bundles \cite{Exel2}
provides a cheap way to build an amplified $C_r^*(\A)$-$C_r^*(\B)$ correspondence which always works. 

\begin{theorem} \label{Fell-abs}
Consider the external tensor product $\widetilde{Y_r}:=\ell^2(G) \otimes Y_r$ as a Hilbert $C_r^*(\B)$-module in the obvious way. Then there exists a $*$-homomorphism $\Psi_r$  from $C_r^*(\A)$ into $\L_{C_r^*(\B)}(\widetilde{Y_r})$ 
such that 
\[ \Psi_r(\lambda^\A_g(a))
= \lambda^G_g \otimes \pi^r_g(a)
\]
for all $g\in G$ and $a\in A_g$, where $\lambda^G$ denotes the left regular representation of $G$ on $\ell^2(G)$. In particular,  $\widetilde{Y_r}$ is a $C_r^*(\A)$-$C_r^*(\B)$ correspondence.
\end{theorem}
\begin{proof} Let $\lambda^G\otimes \pi^r$ be the $*$-representation of $\A$ on $\widetilde{Y_r}$ given by
\[(\lambda^G\otimes \pi^r)_g(a)  = \lambda^G_g \otimes \pi^r_g(a)\] for every $g\in G$ and $a\in A_g$. Then the integrated form of $\lambda^G\otimes \pi^r$ is a $*$-homomorphism $\widetilde\Psi$ from $C^*(\A)$ into $\L_{C_r^*(\B)}(\widetilde{Y_r})$ satisfying $\widetilde\Psi(a\odot g) =  \lambda^G_g \otimes \pi^r_g(a)$ for every $g\in G$ and $a\in A_g$. 
Moreover, adapting the proof of \cite[Proposition 18.4]{Exel2} to $*$-representations of Fell bundles on Hilbert $C^*$-modules, we get that $\widetilde\Psi$
factors through $C_r^*(\A)$. Thus, there exists a $*$-homomorphism $\Psi_r:C_r^*(\A) \to \L_{C_r^*(\B)}(\widetilde{Y_r})$ such that $\widetilde\Psi = \Psi_{r}\circ \Lambda^\A$. We then have 
\[ \Psi_r(\lambda^\A_g(a)) = \Psi_r(\Lambda^\A(a\odot g)) =  \widetilde\Psi(a\odot g)
= \lambda^G_g \otimes \pi^r_g(a)
\]
for all $g\in G$ and $a\in A_g$.
\end{proof}

\begin{remark}
Assume $\X$ is the trivial Hilbert $\B$-bundle $\B$,
and $\rho=\rho_\B$ is the trivial $\B$-action on $\X=\B$, as in Example \ref{trivialBact}. Then $C_c(\X)=C_c(\B)$ is equipped with the right action of $C_c(\B)$ given by $\xi \cdot f = \xi \star f$ and the $C_c(\B)$-valued inner product $\langle \xi, \eta \rangle_{C_c(\B)} = \xi^* \star \eta$, where $f,\xi, \eta \in C_c(\B)$. Further, the left action of $C_c(\B)$ on itself associated to $\rho= \rho_\B$ is simply given by $f \cdot \xi = f \star \xi$.
It readily follows that the associated Hilbert $C^*(\B)$-module $Y = \B \rtimes \B$ is nothing but the trivial $C^*(\B)$-module $C^*(\B)$, equipped with the natural left action of $C^*(\B)$ on itself, i.e., $Y = \B \rtimes \B$ is the trivial $C^*(\B)$-$C^*(\B)$ correspondence. Similarly, $Y_r = \B \rtimes_ r \B$ gives rise to the trivial $C_r^*(\B)$-$C_r^*(\B)$ correspondence.
\end{remark}

\begin{remark} Let $\rho=\widecheck{\rho_\B}$  be the regular $\B$-action  on the regular Hilbert $\B$-bundle $\widecheck \B$ introduced in Example \ref{regB-action}. By unwrapping the relevant definitions, elements in $C_c(\widecheck \B)$ can be approximated by finitely supported functions $\xi$ on $G$ such that, for each $g \in G$, $\xi(g)$ is a finitely supported function from $G$ into $B_g$ or, in other words, setting $\xi(g,t):= (\xi(g))(t)$, by finitely supported functions from $G \times G$ into $\B$ such that $\xi(g,t) \in B_g$, for all $g,t \in G$. Then 
one may translate the various operations in this language to get a handy description of $C_c(\widecheck \B)$ and thus of $Y = \widecheck \B \rtimes \B $ as a $C^*(\B)$-$C^*(\B)$ correspondence. 
When $\B$ is saturated, it can also be shown after some efforts that the left action of $C_c(\B)$ on $C_c(\widecheck\B)$ is bounded w.r.t.~$\|\cdot\|_r$; see  Example \ref{B-saturated}.
\end{remark}

\begin{example} \label{C*-c-dyn-syst} Let $\Sigma=(A, G, \alpha)$ and $\Omega=(B, G, \beta)$ be $C^*$-dynamical systems and $\gamma$ be 
an $\alpha$-$\beta$ compatible action of $G$ on an $A$-$B$ correspondence $X$, cf.~subsection \ref{action-ds}.
We can then form the $\A_\Sigma$-action $\rho_\gamma$ on the Hilbert $\B_\Omega$-bundle $\X$ associated to $X$  and $\Omega$. Since  $C^*(\Sigma)\simeq C^*(\B_\Sigma)$ and $C^*(\Omega) \simeq C^*(\B_\Omega)$, we get from Theorem \ref{correspond} a $C^*(\Sigma)$-$C^*(\Omega)$ correspondence $\X\rtimes \B_\Omega$. 

Identifying $C_c(\X)$ with $C_c(G,X)$, and similarly $C_c(\B_\Omega)$ (resp. $C_c(\A_\Sigma)$) with $C_c(G,B)$ (resp. $C_c(G,A)$),  one checks that the operations on $C_c(\X) = C_c(G,X)$  are given as follows:
\begin{align*} (\xi \cdot f)(h) &= \sum_{k\in G} \xi(k)\beta_k(f(k^{-1}h)) \quad \text{(right $C_c(G, B)$-action)},\\
\langle \xi,\eta \rangle_{C_c(G, B)}(h) &= \sum_{k\in G} \beta_k^{-1}(\langle \xi(k),\eta(kh)\rangle_B) \quad \text{($C_c(G, B)$-valued inner product)},\\
(g \cdot \xi)(h) &= \sum_{k\in G} g(k)\gamma_k(\xi(k^{-1}h)) \quad \text{(left $C_c(G, A)$-action)}
\end{align*}
for $f \in C_c(G, B), \xi, \eta \in C_c(G, A), g \in C_c(G, A)$ and $h\in G$. Hence, we see that $\X\rtimes \B_\Omega$ coincides with the $C^*(\Sigma)$-$C^*(\Omega)$ correspondence previously constructed in \cite[Section 6, Definition 1 and Theorem 1]{Kasp95} and in \cite[Proposition 3.5]{EKQR00} (when assuming that the action of $A$ on $X$ is non-degenerate). Interestingly, it is shown in \cite[Lemma 3.10]{Kasp} and \cite[Proposition 3.2]{EKQR06} that the left action of $C_c(G, A)$  extends to a left action of $C_r^*(\Sigma)$ on $\X\rtimes_r \B_\Omega$. 
\end{example}

 \begin{remark} Assume that $\A=(A_g)_{g\in G}$ and $\B=(B_h)_{h\in H}$ are unital Fell bundles,  $T$ is a positive definite $\A$-$\varphi$-$\B$ bundle map and $\ker(\varphi)$ is not amenable.  Using Theorem \ref{GR}, 
 we can write $T$ as the diagonal coefficient of an $(\A, \varphi)$-action $\rho$ on a  Hilbert $\B$-bundle $\X$, and form the associated Hilbert $C_r^*(\B)$-module $Y_r  = \X \rtimes_r \B$. Now, recall that it is not necessarily true that $T$ is reduced with $M_T$ completely positive,  cf.~the discussion before Theorem 3.10 in \cite{BeCo7}.  
This means that
it is not always possible to extend the left action of $C_c(\A)$ on $C_c(\X)$  to an action of $C_r^*(\A)$ on $Y_r$. Otherwise, we could argue in a similar way as in Remark \ref{full-compl} to deduce that  $T$ has to be reduced with $M_T$ completely positive. 
A simple  example illustrating when such a situation arises is given below.
\end{remark} 

\begin{example} Consider the trivial group bundles $\A= ( \Complessi \times \{g\} )_{g\in G}$ and $\B= ( \Complessi \times \{h\} )_{h\in H}$.  
 Let $\rho$ be the $(\A, \varphi)$-action  on $\X:=\B$ (considered as the trivial Hilbert $\B$-bundle) given by
\[ [\rho(z, g)](w, h)= (zw , \varphi(g) h) \quad \text{for all } z, w \in\Complessi, g \in G, h\in H,\]
and form the associated Hilbert $C_r^*(\B)$-module $Y_r$ (which is nothing but $C_r^*(\B)$, as a Hilbert $C^*$-module over itself). 
Set $x= (1, e) \in X_e= B_e$, and let $T$ the positive definite $\A$-$\varphi$-$\B$ bundle map  given by
\[  T(z, g):= \big\langle x, \rho(z, g) x\big\rangle_\B= \big\langle (1, e), (z, \varphi(g)) \big\rangle_\B = (1, e)^*(z, \varphi(g)) = (z, \varphi(g))\]
for all $z\in \Complessi, g\in G$.      
 Then we claim that if the left action of $C_c(\A)$ on $C_c(\X)$ can be extended to a left action of $C_r^*(\A)$ on $Y_r$, then $\ker(\varphi)$ has to be amenable.  
Indeed, assume that this left action extends.  
Then, arguing similarly as in Remark \ref{full-compl}, we get that $T$ is reduced with $M_T: C_r^*(\A)\to C_r^*(\B)$ completely positive. Now, identifying $C_r^*(\A)$ with $C_r^*(G)$, and $C_r^*(\B)$ with $C_r^*(H)$, it follows that the natural $*$-homomorphism from $C_c(G)$ into $C_c(H)$ induced by $\varphi$ extends to a $*$-homomorphism from $C_r^*(G)$ into $C_r^*(H)$ (mapping each $\lambda_G(g)$ to $\lambda_H(\varphi(g))$). But, as shown in \cite{BedlH}, this is true (if and) only if $\ker(\varphi)$ is amenable, as claimed above.

Thus, if $\ker(\varphi)$ is not amenable, then  the left action of $C_c(\A)$ on $C_c(\X)$ can not be extended. 
\end{example}

 It would be interesting to know whether
 Theorem \ref{Fell-abs} can  be used to show that the left action of $C_c(\A)$ on $C_c(\X)$ always extends to an action of $C_r^*(\A)$ on $Y_r$ whenever $\ker(\varphi)$ is amenable, or at least when $H=G$ and $\varphi={\rm id}_G$.  

 \subsection{An application to Morita equivalence} \label{Morita}
In this subsection, we let $\A$ and $\B$ be two Fell bundles, both over the same discrete group $G$.  We recall \cite{ABF1} that a (right) Hilbert $\B$-bundle $\X$ is called \emph{full} 
when $\overline{{\rm Span}( \langle X_r, X_r\rangle_\B)} = B_e$ for every $r\in G$.
 Equivalently,  $\X$ is full if $B_r= \overline{{\rm Span}(\{ \langle X_s, X_t\rangle_\B : s^{-1}t =r\})}$ for all $r\in G$, 
 cf.~\cite[Remark 2.4]{ABF1}. 
It is a routine exercise to check that if $\X$ is a full right Hilbert $\B$-bundle, then $Y=\X\rtimes \B$, as defined in the previous subsection, is a full right $C^*(\B)$-module. 

We also recall that there is an analogous notion of  \emph{left} Hilbert $\A$-bundle, along with 
a similar notion of fullness for such a bundle (see \cite{ABF1}). As usual, we will denote the $\A$-valued inner product on a left $\A$-bundle by ${}_\A\langle \cdot, \cdot\rangle$. 

According to \cite{ABF1}, $\A$ and $\B$ are said to be \emph{weakly equivalent} 
if there exists a Banach bundle $\X$ over $G$ which is  a  full 
left Hilbert $\A$-bundle and a  full 
right Hilbert $\B$-bundle satisfying that ${}_\A\langle x, y\rangle\, z = x\,\langle y, z\rangle_\B$ for all $x, y, z\in \X$. As shown in \cite[Theorem 4.5]{AF}, $C^*(\A)$ is Morita equivalent to $C^*(\B)$
whenever  $\A$ and $\B$ are weakly equivalent.\footnote{Weak equivalence of Fell bundles is just called equivalence in \cite{AF}. However, in \cite{AF}, fullness is included in the definition of a (left or right) Hilbert $\B$-bundle. In this paper we follow the terminology used in \cite{ABF1}. 
A notion of \emph{strong} equivalence is introduced in \cite{ABF1}, but it is not relevant for our discussion.}

Assume now that $\A$ and $\B$ are weakly equivalent and let $\X$ be a Banach bundle $\X$ over $G$ implementing the weak equivalence between $\A$ and $\B$. We will show below that $Y=\X\rtimes \B$ can be turned into a $C^*(\A)$-$C^*(\B)$ imprimivity bimodule (as defined for example in \cite{EKQR06}), thus providing an alternative and more direct proof of the Morita equivalence between $C^*(\A)$ and $C^*(\B)$.

\begin{proposition} With assumptions as above, $Y=\X\rtimes \B$ is a $C^*(\A)$-$C^*(\B)$ imprimitivity bimodule.
\end{proposition} 

\begin{proof} Let $\rho:\A\to \mathcal{F}(\X)$ be defined by $\rho(a)x = a\cdot x$ for all $a\in \A$ and $x \in \X$.  Using that $\X$ is left Hilbert $\A$-bundle and \cite[Lemma 2.7, (vii) and (ix)]{AF}, we get that $\rho$ is an $\A$-action on the right Hilbert $\B$-bundle $\X$. Theorem \ref{correspond} now gives that 
$Y$ is a $C^*(\A)$-$C^*(\B)$ $C^*$-correspondence.  Since $\X$ is a full as a right Hilbert $\B$-bundle, it is clear that $Y$ is a full right $C^*(\B)$-module. 
Now, for $f\in C_c(\A)$ and $\xi \in C_c(\X)$, recall that $f\cdot \xi \in C_c(\X)$ is given by
\[ (f\cdot\xi)(h) = \sum_{k\in G} \rho(f(k)) \xi(k^{-1}h)=  \sum_{k\in G} f(k) \cdot \xi(k^{-1}h)\] 
for all $h\in G$. Define ${}_{C_c(\A)}\langle \cdot, \cdot\rangle:  C_c(\X)\times C_c(\X) \to C_c(\A) \subseteq C^*(\A)$ by 
\[{}_{C_c(\A)}\langle \xi, \eta\rangle(h) := \sum_{k\in G} \, {}_{\A}\langle \xi(hk), \eta(k) \rangle\]
for all $\xi, \eta \in C_c(\X)$ and $h\in G$.
In particular, for $g, g' \in G$, $x\in X_g$ and $y \in X_{g'}$, we get that
\[ {}_{C_c(\A)}\langle x\odot g, y \odot g'\rangle(h) =  \sum_{k\in G} \, {}_{\A}\langle (x\odot g)(hk), (y \odot g')(k) \rangle = \big({}_{\A}\langle x, y\rangle  \odot (gg'^{-1})\big)(h)\]
for all $h \in G$, hence that ${}_{C_c(\A)}\langle x\odot g, y \odot g'\rangle = {}_{\A}\langle x, y\rangle  \odot (gg'^{-1})$.

It is straightforward to verify that \[({}_{C_c(\A)}\langle \xi, \eta\rangle)^* = {}_{C_c(\A)}\langle \eta, \xi\rangle,  \quad{}_{C_c(\A)}\langle f \cdot \xi, \eta\rangle = f \star {}_{C_c(\A)}\langle \xi, \eta\rangle\] for all $\xi, \eta\in C_c(\X)$ and $f \in C_c(\A)$. Further, arguing similarly as in the proof of Proposition \ref{indprod}, one gets 
that ${}_{C_c(\A)}\langle \cdot, \cdot\rangle$ is a $C^*(\A)$-valued inner product on the left  $C_c(\A)$-module $C_c(\X)$. 

We also have  that 
\begin{equation} \label{eq-A-B}
{}_{C_c(\A)}\langle \xi, \eta\rangle\,\cdot \zeta = \xi \cdot \langle \eta, \zeta\rangle_{C_c(\B)}
\end{equation} for all $\xi, \eta, \zeta\in C_c(\X)$.
For this it suffices to check that 
\[{}_{C_c(\A)}\langle x\odot g, y\odot h\rangle\,\cdot (z\odot k) = (x\odot g) \cdot \langle y\odot h, z\odot k\rangle_{C_c(\B)}\] for all $g, h, k \in G$ and $x\in X_g, y\in X_h , z \in X_k$. We compute:
\begin{align*} {}_{C_c(\A)}\langle x\odot g, y\odot h\rangle\,\cdot (z\odot k) &= \big( {}_{\A}\langle x, y\rangle  \odot (gh^{-1})\big) \cdot (z\odot k)
= \big({}_{\A}\langle x, y\rangle \cdot z\big) \odot (gh^{-1}k) \\ 
&= \big( x\,\langle y, z\rangle_\B \big) \odot (gh^{-1}k) =
 (x\odot g) \cdot \big(\langle y, z\rangle_\B \odot h^{-1}k\big)\\
&=(x\odot g) \cdot \langle y\odot h, z\odot k\rangle_{C_c(\B)}
\end{align*}
Now, for every $\xi \in C_c(\X)$ and $f \in C_c(\B)$, we observe that 
\begin{equation}\label{eq-left} \|{}_{C_c(\A)}\langle \xi\cdot f, \xi\cdot f\rangle \|_{C^*(\A)}  \leq \| {}_{C_c(\A)}\langle \xi, \xi\rangle\|_{C^*(\A)}\;  \|f\|^2_{C^*(\B)}.
\end{equation}
The point here is 
that we can proceed in a perfectly analogous way as we did in Section \ref{mainconstr}, now with $C_c(\B)$ acting on the right of the inner product $C_c(\A)$-module $C_c(\X)$ via $\xi\mapsto \xi\cdot f$ to construct a $*$-representation of $C^*(\B)$ by adjointable (hence bounded) operators on the left $C^*(\A)$-module $Y'$ obtained by completion of   $C_c(\X)$ w.r.t.~the norm given 
 by $\|\xi\|': = \|{}_{C_c(\A)}\langle \xi, \xi\rangle \|^{1/2}_{C^*(\A)}$. 

Applying 
the inequality (\ref{eq-left}) in the form $\|\xi\cdot f\|' \leq \|\xi\|'\, \|f\|_{C^*(\B)}$ 
with $f={}_{C_c(\A)}\langle \xi, \xi\rangle$
we get 
\begin{align*}
 \|{}_{C_c(\A)}\langle \xi, \xi\rangle\|_{C^*(\A)}^2 
  &= \|{}_{C_c(\A)}\langle \xi, \xi\rangle \; {}_{C_c(\A)}\langle \xi, \xi\rangle\|_{C^*(\A)}\\
  &=\| {}_{C_c(\A)}\big\langle {}_{C_c(\A)}\langle \xi, \xi\rangle\cdot\xi, \xi\big\rangle\|_{C^*(\A)}\\
  & = \| {}_{C_c(\A)}\big\langle \xi\cdot \langle \xi, \xi\rangle_{C_c(\B)}, \xi\big\rangle\|_{C^*(\A)}\\
  &\leq \|\xi\cdot \langle \xi, \xi\rangle_{C_c(\B)}\|' \; \|\xi\|'\\
&\leq \|\xi\|' \, \| \langle \xi, \xi\rangle_{C_c(\B)}\|_{C^*(\B)} \, \|\xi\|' \\
&= \| \langle \xi, \xi\rangle_{C_c(\B)}\|_{C^*(\B)} \;
 \|{}_{C_c(\A)}\langle \xi, \xi\rangle\|_{C^*(\A)},
 \end{align*}
 i.e., $\|{}_{C_c(\A)}\langle \xi, \xi\rangle\|_{C^*(\A)} \leq \| \langle \xi, \xi\rangle_{C_c(\B)}\|_{C^*(\B)}$. Hence  
 \[\|{}_{C_c(\A)}\langle \xi, \xi\rangle\|_{C^*(\A)} = \| \langle \xi, \xi\rangle_{C_c(\B)}\|_{C^*(\B)}\] by symmetry. This means that both the associated norms on $C_c(\X)$ agree, so $Y'=Y$. 
 In particular,  ${}_{C_c(\A)}\langle \cdot, \cdot\rangle$ 
 extends to a  $C^*(\A)$-valued inner product on $Y$,  
 and it follows  from (\ref{eq-A-B}) by continuity that   
 \[{}_{C^*(\A)}\langle F, G\rangle\,\cdot H = F \cdot \langle G, H\rangle_{C^*(\B)}\] for all $F, G, H\in Y$.
 As $\X$ is also assumed to be full as a left Hilbert $\A$-bundle, 
 we get that $Y$ is full as 
 a left Hilbert $C^*(\A)$-module. All in all, this shows that $Y$ is a $C^*(\A)$-$C^*(\B)$ imprimitivity bimodule, 
 as desired.
\end{proof}

\begin{remark}
Abadie and Ferraro also show that $C_r^*(\A)$ and $C_r^*(\B)$ are Morita equivalent whenever $\A$ and $\B$ are weakly equivalent, cf.~\cite[Proposition 4.13]{AF}. 
It is conceivable that the weak equivalence of  $\A$ and $\B$ implies that  the right Hilbert $C_r^*(\B)$-module  $Y_r = \X\rtimes_r\B$  can be turned into a $C_r^*(\A)$-$C_r^*(\B)$ imprimitivity bimodule. It is not clear to us whether the proof given above in the full case can be adapted to handle the reduced one.
\end{remark}
 
\subsection{Results in the reduced case} \label{C*-c-red}
In this subsection, we assume that the Fell bundles $\A=(A_g)_{g\in G}$ and $\B=(B_h)_{h\in H}$ are unital and $\varphi \in {\rm Hom}(G,H)$.
 We first note that if $T$ is a positive definite $\A$-$\varphi$-$\B$ bundle map, 
then proceeding essentially as in the beginning of section \ref{GR-1}, 
 we can use Paschke's theorem to construct  
 a $C^*(\A)$-$C^*(\B)$ correspondence $Z^T$ associated to the completely positive map $\Phi_T: C^*(\A)\to C^*(\B)$. 
 Similarly, if  ${\rm ker}(\varphi)$ is supposed to be amenable, we can also construct a $C_r^*(\A)$-$C_r^*(\B)$ correspondence $Z^T_r$ associated to $M_T$.

Let now $\rho$ be an $(\A, \varphi)$-action on a  Hilbert $\B$-bundle $\X=(X_h)_{h\in H}$. Let $x \in X_e$ and consider  the positive definite $\A$-$\varphi$-$\B$ bundle map $T^{\rho,x}$ given by 
\[T^{\rho, x}_g(a) = \langle x, \rho(a)x\rangle_\B, \quad g\in G, a \in A_g.\] 
We can then form the $C^*(\A)$-$C^*(\B)$ correspondence $Z^{\rho,x}:= Z^{T^{\rho, x}}$ associated to $\Phi_{T^{\rho, x}}$. Similarly, we can form the $C_r^*(\A)$-$C_r^*(\B)$ correspondence $Z_r^{\rho,x}$ associated to $M_{T^{\rho, x}}$ whenever $\ker(\varphi)$ is amenable. 

We will first discuss how $Z^{\rho, x}$ is related  to the $C^*(\A)$-$C^*(\B)$ correspondence 
$Y= \X\rtimes \B$ associated to $\rho$ in Theorem \ref{correspond}. For $k\in H$ and $y\in X_k$, we let $y\odot k \in C_c(\X)$ be defined by 
\[ (y\odot k)(h) = \begin{cases} y & \quad \text{if } h=k,\\ 0 &\quad \text{otherwise}.\end{cases}\]
 
\begin{proposition} 
Assume $\A$ and $\B$ are unital and $x \in X_e$. Let $Z^{\rho, x}$  be the $C^*(\A)$-$C^*(\B)$ correspondence defined above. Then
$Z^{\rho, x}$ can be identified  with the 
sub-$C^*$-correspondence $Y^{\rho,x}$ of $Y$ 
  defined by
  \[ Y^{\rho, x} := \overline{{\rm span}\{ (\rho(a)x)b \odot \varphi(g)h : g \in G, h\in H, a \in A_g, b\in B_h\}}^{\,\|\cdot\|_Y}.\] 
\end{proposition}
\begin{proof}  By Paschke's construction, $Z^{\rho, x}$ arises as the Hausdorff completion of the  $C^*(\A)$-$C^*(\B)$ bimodule 
  $Z_0^{\rho, x} = 
  C^*(\A)\odot C^*(\B)$ equipped with a $C^*(\B)$-valued semi-inner product $[\cdot, \cdot]$  satisfying 
  \begin{align*}
\Big[ \widehat j^\A_g(a) \odot \widehat j^\B_h(b), \,& \widehat j^\A_{g'}(a')\odot \widehat j^\B_{h'}(b')\Big]  = \widehat j^\B_h(b)^*\, \Phi_{T^{\rho, x}}\big(\,\widehat j^\A_g(a)^* \widehat j^\A_{g'}(a')\big) \,\widehat j^\B_{h'}(b')\\
   &= \widehat j^\B_h(b)^*\, \Phi_{T^{\rho, x}}\big(\,\widehat j^\A_{g^{-1}g'}(a^*a')\big) \,\widehat j^\B_{h'}(b')\\
  & = \widehat j^\B_h(b)^*\,\widehat j^\B_{\varphi(g^{-1}g')}\big(T^{\rho,x}_{g^{-1}g'}(a^*a')\big) \,\widehat j^\B_{h'}(b')\\
   & = \widehat j^\B_h(b)^*\,\widehat j^\B_{\varphi(g^{-1}g')}\big(\langle x, \rho(a^*a')x\rangle_\B\big) \,\widehat j^\B_{h'}(b')\\
  & = \widehat j^\B_{h^{-1}}(b^*)\,\widehat j^\B_{\varphi(g^{-1}g')}\big(\langle \rho(a)x, \rho(a')x\rangle_\B\big) \,\widehat j^\B_{h'}(b')\\
  & = \widehat j^\B_{h^{-1}\varphi(g)^{-1}\varphi(g')h'}\big(b^*\langle \rho(a)x, \rho(a')x\rangle_\B \,b'\big) \\
  & = \widehat j^\B_{(\varphi(g)h)^{-1}\varphi(g')h'}\big(\langle (\rho(a)x)b, (\rho(a')x)b'\rangle_\B \big)
  \end{align*}
  for all $g, g'\in G$, $h, h\in H$, $a\in A_g, b \in B_h, a'\in A_{g'}$ and $b'\in B_{h'}$.  We will denote the associated seminorm on $Z^{\rho, x}_0$ by $\|\cdot\|_0$, i.e.,  $\|z\|_0:= \|[z,z]\|^{1/2}$ for $z \in Z^{\rho, x}_0$. 
  
  On the other hand, identifying $C_c(\X)$ with its canonical copy in $Y$, the $C^*(\B)$-valued inner product on $Y$ satisfies 
     \[\langle y_1 \odot h_1, y_2\odot h_2\rangle_{C^*(\B)} 
    = \widehat j^\B_{h_1^{-1}h_2}( \langle y_1, y_2\rangle_\B )\]
    where   $h_j \in H$ and $y_j \in X_{h_j}$, $j=1,2$.
    
    For $x \in X_e$, define now a linear map $ V_{00}^x: C_c(\A)\odot C_c(\B) \to C_c(\X) \subseteq Y$ by
    \[ V_{00}^x\Big(\sum_{i=1}^m \widehat j^\A_{g_i}(a_i) \odot \sum_{j=1}^n \widehat j^\B_{h_j}(b_j)\Big) = \sum_{i=1}^m\sum_{j=1}^n (\rho(a_i)x) b_j \odot \varphi(g_i)h_j,\]
 where $g_1, \ldots, g_m\in G$, $h_1, \ldots, h_n \in H$, $a_i \in A_{g_i}, i=1, \ldots, m$, $ b_j \in B_{h_j}, j=1, \ldots, n$.   
 By adding some zeros if necessary, we can assume that $m=n$.   
    Then we have that
    \begin{align*}\big\|\sum_{i,j=1}^n (\rho(a_i)x) b_j \odot \varphi(g_i)h_j\big\|_{Y}^2 &= \big\| \langle \sum_{i,j=1}^n (\rho(a_i)x) b_j \odot \varphi(g_i)h_j, 
    \sum_{k,l=1}^n (\rho(a_k)x) b_l \odot \varphi(g_k)h_l\rangle_{C^*(\B)}\big\|\\
    &= \big\|  \sum_{i,j,k, l=1}^n \langle (\rho(a_i)x) b_j \odot \varphi(g_i)h_j, 
     (\rho(a_k)x) b_l \odot \varphi(g_k)h_l\rangle_{C^*(\B)}\big\| \\
     &=  \big\|  \sum_{i,j,k, l=1}^n \widehat j^\B_{(\varphi(g_i)h_j)^{-1}\varphi(g_k)h_l}\big( \langle \rho(a_i)x) b_j, (\rho(a_k)x) b_l\rangle_\B\big) \big\|\\
     &=   \big\|  \sum_{i,j,k, l=1}^n \big[ \,\widehat j^\A_{g_i}(a_i) \odot \widehat j^\B_{h_j}(b_j), \widehat j^\A_{g_k}(a_k)\odot \widehat j^\B_{h_l}(b_l)\big] \big\|\\
     &=   \big\|  \, \Big[ \sum_{i,j=1}^n\widehat j^\A_{g_i}(a_i) \odot \widehat j^\B_{h_j}(b_j), \sum_{k, l=1}^n\widehat j^\A_{g_k}(a_k)\odot \widehat j^\B_{h_l}(b_l)\Big]\, \big\|\\
     &= \big\| \sum_{i,j=1}^n\widehat j^\A_{g_i}(a_i)  \odot \widehat j^\B_{h_j}(b_j) \big\|^2_0.
\end{align*}
Thus we get that the map $V_{00}^x$ is well-defined, and extends to a linear map \[V_{0}^x:C^*(\A)\odot C^*(\B) \to Y\]
satisfying that $\|V_0^x(z)\|_Y = \|z\|_0$ for all $z\in Z^{\rho,x}_0$. Setting $N=\{z\in Z^{\rho,x}_0: [z,z]=0\}$, the map $z+N\mapsto V_0^x(z)$ is then a well-defined linear isometry from $Z^{\rho,x}_0/N$ into $Y$, which  extends to a  linear isometry $ V^x: Z^{\rho, x} \to Y$. 

It is not difficult to check  that $V^x$ is  $C^*(\B)$-linear, i.e., $V^x(zw) = (V^x(z))w$ for all $z \in Z^{\rho,x}$ and $w \in C^*(\B)$. We claim that it also intertwines the respective left actions of $C^*(\A)$ on $Z^{\rho, x}$ and $Y$. 
To verify this, we first   
recall  that $(f\cdot \xi)(h)= \sum_{k\in G} \rho(f(k)) \xi(\varphi(k)^{-1}h)$ and $\Phi_\pi(f) \xi = f\cdot \xi$ for all $f\in C_c(\A), \xi\in C_c(\X)$ and $h\in H$.   
We then have that 
\[ [\Phi_\pi(\widehat j^\A_g(a))](x' \odot h) =  (a\odot g)\cdot (x' \odot h) = \rho(a)x' \odot \varphi(g)h\]
for all $g\in G, h\in H, a\in A_g$ and $x'\in X_h$. Thus we get that
\begin{align*} V_{00}^x\Big(\widehat j^\A_g(a)\cdot(\widehat j^\A_{g'}(a')\odot \widehat j^\B_{h'}(b'))\Big) &= 
V_{00}^x\Big(\widehat j^\A_{gg'}(aa')\odot \widehat j^\B_{h'}(b')\Big)\\
&= (\rho(aa')x)b' \odot \varphi(gg')h' \\ 
&= (\rho(a)\rho(a')x)b' \odot \varphi(g)\varphi(g')h' \\ 
&= \rho(a)((\rho(a')x)b') \odot \varphi(g)\varphi(g')h' \\ 
&= [\Phi_\pi \big(\widehat j^\A_g(a)\big)]((\rho(a')x)b' \odot \varphi(g')h')\\ 
&= [\Phi_\pi \big(\widehat j^\A_g(a)\big)] V_{00}^x\Big(\widehat j^\A_{g'}(a')\odot \widehat j^\B_{h'}(b')\Big)
\end{align*}
for all $g, g'\in G$, $a\in A_g$, $a'\in A_{g'}$, $h'\in H$ and $b\in B_{h'}$. It thus follows that $V_{00}^x$ intertwines the actions of $C_c(\A)$ on $C_c(\A)\odot C_c(\B)$ and $C_c(\X)$, and it is now routine to see that our claim holds. 

Hence we get that $Z^{\rho, x}$ can be identified via $V^x$ with the 
sub-$C^*$-correspondence $Y^{\rho,x}$ of $Y$ 
  defined by $Y^{\rho, x} = \overline{{\rm span}\{ (\rho(a)x)b \odot \varphi(g)h : g \in G, h\in H, a \in A_g, b\in B_h\}}^{\,\|\cdot\|_Y}$, as 
  desired.
  \end{proof}
It follows from this proposition that if $x \in X_e$ can be chosen such that $Y^{\rho,x} = Y$, that is, such that ${\rm span}\{ (\rho(a)x)b \odot \varphi(g)h : g \in G, h\in H, a \in A_g, b\in B_h\}$ is dense in $Y$, then $Z^{\rho, x}$ will be unitarily equivalent to $Y$. 
For instance, if $x\in X_e$ is cyclic for $\rho$ (in the sense of Remark \ref{cyclic1}), 
 then it can be shown that $Y^{\rho,x} = Y$. As we will prove a similar statement in the reduced case (see Proposition \ref{Y-rho-r-2}), we skip the proof. 
   Since $C_c(\X) = {\rm span} \{ y \odot k : k\in H, y\in X_k\}$ is dense in $Y$, this will for example happen when $x \in X_e$ 
   is such that
   ${\rm span}\{(\rho(a)x)b: a \in A_e, b\in B_r\}$ is dense in $X_r$
   for every $r\in H$.
 
 We now turn our attention to the reduced case and therefore assume that $\ker(\varphi)$ is amenable. Let $x \in X_e$. We can then consider the 
 $C_r^*(\A)$-$C_r^*(\B)$ correspondence $Z_r^{\rho,x}$ associated to the completely positive map  $M_{T^{\rho, x}}:C_r^*(\A) \to C_r^*(\B)$. 
 
Analogously to the full case,  $Z_r^{\rho, x}$ arises as the Hausdorff completion of the  $C_r^*(\A)$-$C_r^*(\B)$ bimodule 
  $Z_{0, r}^{\rho, x} := 
  C_r^*(\A)\odot C_r^*(\B)$ equipped with a $C_r^*(\B)$-valued semi-inner product $[\cdot, \cdot]_r$  satisfying 
  \begin{align*}
\Big[\lambda^\A_g(a) \odot \lambda^\B_h(b), \,& \lambda^\A_{g'}(a')\odot \lambda^\B_{h'}(b')\Big]_r  
   = \lambda^\B_{(\varphi(g)h)^{-1}\varphi(g')h'}\big(\langle (\rho(a)x)b, (\rho(a')x)b'\rangle_\B \big)
  \end{align*}
  for all $g, g'\in G$, $h, h\in H$, $a\in A_g, b \in B_h, a'\in A_{g'}$ and $b'\in B_{h'}$.  We will denote the associated seminorm on $Z^{\rho, x}_{0,r}$ by $\|\cdot\|_{0,r}$, i.e.,  $\|z\|_{0,r}:= \|[z,z]_r\|_r^{1/2}$ for $z \in Z^{\rho, x}_{0,r}$. 
  Moreover, we will identify $C_c(\X)$ with its canonical copy in $Y_r$, so the $C_r^*(\B)$-valued inner product on $Y_r$ satisfies 
     \[\langle x_1 \odot h_1, x_2\odot h_2\rangle_{C_r^*(\B)} = \lambda^\B_{h_1^{-1}h_2}( \langle x_1, x_2\rangle_\B ).\]
Then, essentially as above, we get a linear map $V_{0,r}^x:Z^{\rho,x}_{0,r} \to Y_r$ determined by 
 \[ V_{0,r}^x\Big(\sum_{i=1}^m \lambda^\A_{g_i}(a_i) \odot \sum_{j=1}^n \lambda^\B_{h_j}(b_j)\Big) = \sum_{i=1}^m\sum_{j=1}^n (\rho(a_i)x) b_j \odot \varphi(g_i)h_j\]
 whenever $g_1, \ldots, g_m\in G$, $h_1, \ldots, h_n$, $a_i \in A_{g_i}, i=1, \ldots, m$, $ b_j \in B_{h_j}, j=1, \ldots, n$, which satisfies 
that $\|V_{0,r}^x(z)\|_{Y_r} = \|z\|_{0,r}$ for all $z\in Z^{\rho,x}_{0,r}$. 

Setting $N_r=\{z\in Z^{\rho,x}_{0,r}: [z,z]_r=0\}$, the map $z+N_r\mapsto V_{0,r}^x(z)$ is then a well-defined linear isometry from $Z^{\rho,x}_{0,r}/N_r$ into $Y_r$, which  extends to a  linear isometry $ V_r^x: Z_r^{\rho, x} \to Y_r$.

 Further, setting
  \[ Y_r^{\rho, x} := V_r^x\big(Z_r^{\rho, x}\big) = \overline{{\rm span}\{ (\rho(a)x)b \odot \varphi(g)h : g \in G, h\in H, a \in A_g, b\in B_h\}}^{\,\|\cdot\|_{r}},\] we can consider $V_r^x$ as a unitary operator between the Hilbert $C_r^*(\B)$-modules $Z_r^{\rho,x}$ and $Y_r^{\rho, x}$, i.e., 
    $Z_r^{\rho,x}$ and  $Y_r^{\rho, x}$ are isomorphic as Hilbert  $C_r^*(\B)$-modules via $V_r^x$.  We can now transport the left action of $C_r^*(\A)$ on $Z_r^{\rho,x}$ to a left action of $C_r^*(\A)$ on $Y_r^{\rho, x}$ by setting
   \begin{equation}\label{bullet}
    v\bullet y := V_r^x(v\cdot(V_r^x)^*y))\end{equation}
   for all $y\in Y_r^{\rho, x}$ and $v \in C_r^*(\A)$. We then have that  $Z_r^{\rho,x}$ and $Y_r^{\rho, x}$ are unitarily equivalent as $C_r^*(\A)$-$C_r^*(\B)$ correspondences. 
   
   Let $g', h' \in G, a' \in A_{g'}, b'\in B_{h'}$ and set $y:=(\rho(a')x)b' \odot \varphi(g')h' \in Y_r^{\rho,x}$. Then for all $g\in G, a\in A_g$, we get
  \begin{align*} (a\odot g)\bullet  y&= (a\odot g)\bullet \big((\rho(a')x)b' \odot \varphi(g')h' \big)\\
 &= V_r^x\big((a\odot g)\cdot ((V_r^x)^*(\rho(a')x)b' \odot \varphi(g')h' )\big)\\
 &= V_r^x\big((a\odot g)\cdot \big((\lambda_{g'}^\A(a') \odot \lambda^\B_{h'}(b') ) + N_r\big)\big)\\
  &= V_r^x\big(\big((\lambda_{g}^\A(a)\lambda_{g'}^\A(a') \odot \lambda^\B_{h'}(b') ) + N_r\big)\big)\\
  &=V_r^x\big(\big((\lambda_{gg'}^\A(aa') \odot \lambda^\B_{h'}(b') ) + N_r\big)\big)\\
   &=(\rho(aa')x)b' \odot \varphi(gg')h' \\
     &=\rho(a) \big((\rho(a')x))b'\big) \odot \varphi(g)\varphi(g')h' \\
  &= (a\odot g)\cdot y
  \end{align*}
  It follows that the left action of $C_c(\A) \subseteq C^*_r(\A)$ on $Y_r^{\rho, x}$ defined in (\ref{bullet}) coincides on \[W:= {\rm span}\{ (\rho(a)x)b \odot \varphi(g)h : g \in G, h\in H, a \in A_g, b\in B_h\} \, \subseteq \, C_c(\X)\] with the left action of $C_c(\A)$ defined in (\ref{left-action}).
  This can be used to show the following result.
 
 \begin{theorem} \label{Y-rho-r} 
 Let $\rho$ be an $(\A, \varphi)$-action on a  Hilbert $\B$-bundle $\X=(X_h)_{h\in H}$ and let $Y_r = \X \rtimes_{r} \B$  be the
 Hilbert $C_r^*(\B)$-module obtained by taking the completion of $C_c(\X)$  w.r.t.~$\|\xi\|_r = \big\| \langle \xi, \xi\rangle_{C_c(\B)}\big\|_{C_r^*(\B)}^{1/2} , \xi \in C_c(\X)$.
 Further, assume that $\A, \B$ are unital, $\ker(\varphi)$ is amenable and $Y^{\rho,x}_r = Y_r$ for some $x\in X_e$. Then the left action of $C_c(\A)$ on $C_c(\X)$ defined in $(\ref{left-action})$ extends to a left action of $C_r^*(\A)$ on the Hilbert $C_r^*(\B)$-module $Y_r$, i.e.,
 $Y_r$ is a $C_r^*(\A)$-$C_r^*(\B)$ correspondence.  
 \end{theorem} 
 
 \begin{proof}  By assumption, we have $W \subseteq C_c(\X) \subseteq Y_r = Y^{\rho,x}_r$ and 
 $W$ is dense in  $Y_r$ w.r.t.~$\|\cdot\|_r$.  
We will show that the restriction to $C_c(\A)$ of the left action of $C_r^*(\A)$ on $Y_r$ defined in (\ref{bullet}) coincides on $C_c(\X) \subseteq Y_r$ with the action of $C_c(\A)$ defined in  $(\ref{left-action})$, from which the result will follow.

 Let $f\in C_c(\A)$ and $\xi \in C_c(\X)$. Pick $\xi_n \in W$ such that $\|\xi-\xi_n\|_r \to 0$ as $n\to \infty$. Then, using the above computation, we have
 \[f \bullet \xi = \lim_n f \bullet \xi_n = \lim_n f \cdot \xi_n = f \cdot \xi,\]
where the last equality follows from Remark \ref{mixed-cor}, according to which
 \[ \|f \cdot (\xi -\xi_n)\|_r \leq \|f\|_{C^*(\A)} \, \|\xi-\xi_n\|_r \to 0 \text{ as } n \to \infty.\]
 
 \vspace{-2ex} \end{proof}
 
 The assumption that $Y^{\rho,x}_r = Y_r$ for some $x\in X_e$ is always satisfied when $x$ is cyclic for $\rho$ (in the sense of Remark \ref{cyclic1}):
  \begin{proposition} \label{Y-rho-r-2} Assume $\A, \B$ are unital and $\ker(\varphi)$ is amenable. Moreover, assume that  $\rho$ has a cyclic element $x\in X_e$. Then  $Y^{\rho,x}_r = Y_r$. Hence, Theorem \ref{Y-rho-r} gives that $ Y_r =\X\rtimes_r \B$ is a $C_r^*(\A)$-$C_r^*(\B)$ correspondence. 
 \end{proposition}
 \begin{proof}
 By definition of cyclicity of $x$ for $\rho$, we have that
    for each $k \in H$, \[{\rm span}\{(\rho(a)x)b: g\in G, a \in A_g, h\in H, b\in B_h, \varphi(g)h=k\}\] is dense in $X_k$. Now, 
      \[ Y_r^{\rho, x} = \overline{{\rm span}\{ (\rho(a)x)b \odot \varphi(g)h : g \in G, h\in H, a \in A_g, b\in B_h\}}^{\,\|\cdot\|_{r}} \subseteq Y_r,\] so we have to  show that $Y_r \subseteq Y_r^{\rho, x} $. Since 
      \[Y_r = \overline{{\rm span}\{ x \odot k : k\in H, x \in X_k\}}^{\|\cdot\|_r},\]
      it is enough to prove that $ x \odot k \in Y_r^{\rho, x} $ for each $k \in H$ and $x\in X_k$.  
      
      Let $k \in H$, $x\in X_k$ and $\varepsilon >0$. 
      We can then find $x' \in X_k$ 
   of the form $x'= \sum_{i=1}^n (\rho(a_i)x)b_i$ for some  $g_1, \ldots, g_n \in G$, $h_1, \ldots, h_n\in H$, $a_i \in A_{g_i}$, $b_i \in B_{h_i}$ and $\varphi(g_i)h_i = k$ for $i=1, \ldots n$ such that $\|x-x'\| < \varepsilon$.
       
    Now, consider $w\in X_k$. Using equation (\ref{eq-ind}), we get that
       \[ \|w\odot k\|_r^2 = \|\langle w\odot k, w \odot k\rangle_{C_c(\B)}\|_{C_r^*(\B)} =
       \|\langle w, w\rangle_\B \odot e\|_{C_r^*(\B)}= \|\langle w, w\rangle_\B\| = \|w\|^2,\]
       where we have now used that if $b\in B_e$, then $\| b \odot e\|_{C_r^*(\B)} = \| \lambda^\B_e(b)\| = \|b\|$ (cf.~\cite[Proposition 17.9 (v)]{Exel2} for the last equality). 
       
       This implies that   
      $\|x\odot k - x'\odot k\|_r = \|(x-x')\odot k\|_r = \|x-x'\| < \varepsilon$.  Since
       \[ x'\odot k = \Big(\sum_{i=1}^n (\rho(a_i)x)b_i\Big) \odot k = 
       \sum_{i=1}^n \big( (\rho(a_i)x)b_i \odot \varphi(g_i)h_i\big) \in Y_r^{\rho, x},\]
       the desired assertion follows.
\end{proof}
 
 A natural way to produce  actions having a cyclic element is to start with a positive definite bundle map and apply our Gelfand-Raikov type construction. 
   \begin{corollary} \label{Y-rho-r-3} Assume $\A, \B$ are unital and $\ker(\varphi)$ is amenable. Let $T$ be an $\A$-$\varphi$-$\B$ positive definite bundle map, and write $T = T^{\rho, \xi}$ where $\rho$ is an $(\A, \varphi)$-action on a Hilbert $\B$-bundle $\X$ and $\xi \in X_e$,
    as in Theorem \ref{GR}.
     Then $Y_r = \X \rtimes_r \B $ is a $C_r^*(\A)$-$C_r^*(\B)$ correspondence. 
    \end{corollary} 
    
    \begin{proof} The statement  follows from Proposition \ref{Y-rho-r-2} since, as pointed out in Remark \ref{cyclic1},  $\xi\in X_e$ is always cyclic for $\rho$, by construction. 
    \end{proof}

 We end this section with some illustrative examples.
 
 \begin{example} Consider the trivial action $\rho = \rho_\B$ of a unital Fell bundle $\B = (B_g)_{g\in G}$ on $\X=\B$. In this case, $Y_r = C_r^*(\B)$ is the trivial $C_r^*(\B)$-$C_r^*(\B)$ correspondence. For $x=1_{B_e}$, we have 
 $Y^{\rho,x}_r = \overline{{\rm span}\{ ab \odot gh : g, h \in G, a \in B_g, b\in B_h\}}^{\,\|\cdot\|_{r}}$. 
 It is then clear that  $C_c(\B) $ is contained in $Y^{\rho,x}_r$, so it follows that  $Y^{\rho,x}_r = C_r^*(\B) =Y_r$.
 \end{example}
 
 \begin{example} Again, let $\B$ be a unital Fell bundle. For $i=1, \ldots, n$, set $\X^i = \B$ and let $\rho_i$ be the trivial action of $\B$ on $\X^i$. Then we may consider the direct sum $\rho=\oplus_{i=1}^n \rho_i$, which is an action of $\B$ on $\X:= \oplus_{i=1}^n\X^i$.\footnote{The direct sums involved here are defined in the obvious way.
 A general approach to  direct sums of actions will be discussed at length elsewhere.} Note that $Y_r$ 
 is then easily seen to be the direct sum of $n$ copies of the trivial $C_r^*(\B)$-$C_r^*(\B)$ correspondence, i.e., $Y_r = C_r^*(\B)^n$. 
 For each $j=1, \ldots, n$, let $x_j \in \X$ be defined $x_j(i)=\begin{cases} 1_{B_e} &\text{ if } i=j,\\ 0_{B_e} &\text{otherwise}. \end{cases}$. Then 
 $Y^{\rho, x_j}_r$ is the canonical copy of $C_r^*(\B)$ in $Y_r=C_r^*(\B)^n$ at the $j$-th place. 
 
 For each $i=1, \ldots, n$, let $p_i: \X\to \X^i$ denote the canonical $i$-th projection. Let now $\xi\in C_c(\X)$ and set $\xi_i:= p_i\circ \xi \in C_c(\X^i)$. 
Let $f\in C_c(\B)$. Then 
\[\langle f\cdot \xi , f\cdot \xi\rangle_{C_c(\B)} = \sum_{i=1}^n \langle f \cdot \xi_i,  f \cdot \xi_i \rangle_{C_c(\B)}\, \leq\, \|f\|^2_{C_r^*(\B)}  \,\sum_{i=1}^n \langle \xi_i,   \xi_i \rangle_{C_c(\B)} =\|f\|^2_{C_r^*(\B)}  \, \langle \xi,   \xi \rangle_{C_c(\B)} \]
This implies that the left action of $C_c(\B)$ on $C_c(\X)$ can be extended to an action of $C_r^*(\B)$ on $Y_r$, as it was pretty obvious from the start.
The main point with this example is that there is no $x\in \X$ such that $Y^{\rho, x}= Y$, illustrating that  this assumption in Theorem \ref{Y-rho-r} is not necessary for the conclusion to be true.     
Also, this example opens the way for a discussion of more general situations involving direct sums.
  \end{example} 

\begin{example} \label{B-saturated} Assume $\B=(B_g)_{g\in G}$ be a saturated unital Fell bundle. We recall that $\B$ is said to be \emph{saturated} when $B_{gh}$ is the closed span of $B_gB_h$ for each $g,h$ in $G$ (cf.~\cite[Section 16]{Exel2}).
 Consider the regular $\B$-action $\rho:= \widecheck{\rho_\B}$ 
   on the regular Hilbert $\B$-bundle $\widecheck\B$, cf.~Example \ref{regB-action}.     
  On the underlying inner product $C_c(\B)$-bundle $\big(C_c(G,B_r)\big)_{r\in G}$,  $\rho$ is given by
 \[ ( \rho(a) \xi) (s) = a\,\xi(t^{-1}s)\]
for all  $r, s, t \in G$, $a\in B_t$ and $\xi \in C_c(G, B_r)$. 

 Set $x:= 1_{B_e} \odot e \in C_c(G, B_e) \subseteq (\widecheck\B)_e$. Then we claim that $x$ is cyclic for $\rho$. It will then follow from Proposition \ref{Y-rho-r-2} that the left action of $C_c(\B)$ on $C_c(\widecheck\B)$ associated to  $\rho$ extends to an action of $C_r^*(\B)$ on  the Hilbert $C_r^*(\B)$-module  $\widecheck\B \rtimes_r \B$,
hence that  $\widecheck \B\rtimes_r \B$ becomes a $C^*$-correspondence over $C_r^*(\B)$.

 To check that $x$ is cyclic, let $k\in G$. As  $(\rho(a)x)b= (a \odot g)b = ab \odot g $ for all $g, h\in G$,\\ $a\in B_g, b\in B_h$, we have to show that \[Z_k:=
 {\rm span}\{ab\odot g: g, h\in G, a \in A_g, b\in B_h, gh=k\}\] is dense in $(\widecheck\B)_k$.
 Since $C_c(G, B_k)$ is dense in $(\widecheck\B)_k$, it suffices to show that $C_c(G, B_k)$ lies in the closure of $Z_k$ in  $(\widecheck\B)_k$, hence that  $c\odot g$ belongs to $\overline{Z_k}$ for every $c \in B_k$ and $g \in G$. Let $c \in B_k$, $g \in G$ and $\varepsilon > 0$. Then set $h:= g^{-1}k$, so $gh=k$. By saturation, we can find $a_1, \ldots, a_n \in B_g, b_1, \ldots, b_n \in B_h$ such that $\|c - \sum_{i=1}^n a_ib_i\| < \varepsilon$. 
 Then 
 \[\| c\odot g - \sum_{i=1}^n (a_ib_i \odot g) \| = \big\|(c-\sum_{i=1}^n a_ib_i )\odot g\big\| = \|c - \sum_{i=1}^n a_ib_i\| < \varepsilon,\] 
 where we have used that for $d \in B_k$, so $d\odot g \in (\widecheck\B)_k$, we have   
 \[\|d\odot g\|^2 = \|\langle d\odot g, d\odot g\rangle_\B\| = \big\|\sum_{t\in G} \langle(d\odot g)(t),(d\odot g)(t)\rangle_\B\big\| = \|d^*d\| = \|d\|^2.\] 
 Since $ \sum_{i=1}^n (a_ib_i \odot g)  \in Z_k$, this proves our claim. 
 
 In the same vein, we may also consider the  $\B$-action $\rho$ on the canonical $\ell^2$-bundle $\mathcal{L}^2\B$
associated to $\B$, cf.~Examples \ref{ell2} and \ref{ell2-2}. It is not difficult to check that  $x:=(1_{B_e}\odot e, e)$ is an element in the unit fiber of $\mathcal{L}^2\B$ which is cyclic for $\rho$. Proposition \ref{Y-rho-r-2} may therefore be applied, giving that $\mathcal{L}^2\B \rtimes_r \B$ is a $C^*$-correspondence over $C_r^*(\B)$.

  \end{example}

\section{Some final comments} \label{final}

\noindent{\bf 7.1}\,
An important class of Fell bundles over discrete groups,  
which is thoroughly discussed in \cite{Exel2}, 
 arise by considering partial actions of discrete groups on $C^*$-algebras. 
 If $\B$ denotes a Fell bundle in this class, one may use our present work to form the trivial action of $\B$ on itself, the regular action of $\B$ on $\widecheck \B$ and the canonical action of $\B$ on $\mathcal{L}^2\B$. One can also create actions associated with $*$-representations of $\B$ on Hilbert $C^*$-modules. It would be interesting to investigate if one can define some other type of actions of $\B$, similar to those associated with equivariant representations (or with compatible actions) of the system when the partial action is a genuine action. 

\medskip
\noindent{\bf 7.2}\,
In subsection \ref{action-ds} we only consider usual (discrete) $C^*$-dynamical systems and their associated Fell bundles. There is little doubt that everything we have done may be enhanced to deal with the case of twisted (discrete) $C^*$-dynamical systems. We have skipped including a twist in our exposition for the sake of clarity, and invite the interested reader to work out the details. 

\medskip      
\noindent{\bf 7.3}\,
When $\A$ and $\B$ are Fell bundles over locally compact groups as in \cite{FD2}, Hilbert $\B$-bundles already appear in \cite{AF,ABF1}, and our definition of an $\A$-$\varphi$-$\B$ action may be extended by adding necessary continuity conditions. 
We expect that most of our results can be adapted to this more general setting, but the proofs will certainly present significant technical challenges.

\medskip
\noindent{\bf 7.4}\,  Finally, we 
mention a couple of potentially interesting 
problems about actions of Fell bundles
that fall outside the scope of the present paper but 
that we feel 
should deserve some attention. 
First, is it possible to define a viable notion of amenability for $\B$-actions on Hilbert $\B$-bundles ? 

\noindent
Second, the space of reduced $\B$-bundle maps whose associated maps on $C_r^*(\B)$ are completely bounded might be thought of as the analogue of the space of Herz-Schur multipliers for a $C^*$-dynamical system, cf.~\cite{MTT, MSTT, MT, MPTT, HTT, He} (see also \cite[Example 7.3]{BKMS} for a broader setting). Examples of such Herz-Schur multipliers on $\B$ associated to coefficient maps of  $\B$-actions are described in Remark \ref{pre-act-pd}. A more comprehensive study of Herz-Schur multipliers on $\B$ is a natural task for further investigations.

\bigskip
{\parindent=0pt Addresses of the authors:\\

\smallskip Erik B\'edos, Department of Mathematics, University of
Oslo, \\
P.B. 1053 Blindern, N-0316 Oslo, Norway.\\ E-mail: bedos@math.uio.no. \\

\smallskip \noindent
Roberto Conti, Dipartimento SBAI,
Sapienza Universit\`a di Roma \\
Via A. Scarpa 16,
I-00161 Roma, Italy.
\\ E-mail: roberto.conti@sbai.uniroma1.it
\par}

\end{document}